\documentclass[a4paper,12pt]{article}
\title{
From Hyperbolic to Parabolic Parameters \\
along Internal Rays
}
\author{
Yi-Chiuan Chen 
and
Tomoki Kawahira
}
\usepackage{enumerate}
\usepackage{graphicx}
\usepackage{amsmath}
\usepackage{amssymb}
\usepackage{amsthm}
\usepackage{color}

\usepackage{tocloft}
\numberwithin{equation}{section}
\theoremstyle{plain}
\newtheorem{thm}{Theorem}[section]
\newtheorem{prop}[thm]{Proposition}

\newtheorem{cor}[thm]{Corollary}
\theoremstyle{remark}
\newtheorem{remark}[thm]{Remark}
\theoremstyle{definition}
\newtheorem{eg}{Example}
\newcommand{\C}{\mathbb{C}}
\newcommand{\Cbar}{\overline{\C}}
\newcommand{\R}{\mathbb{R}}
\newcommand{\D}{\mathbb{D}}
\newcommand{\Dbar}{\overline{\D}}

\newcommand{\N}{\mathbb{N}}
\newcommand{\M}{\mathbb{M}}

\newcommand{\X}{\mathbb{X}}
\newcommand{\Xbar}{\overline{\X}}
\newcommand{\abs}[1]{{\left| #1 \right|}}
\newcommand{\paren}[1]{{\left( #1 \right)}}
\newcommand{\brac}[1]{{\left\{ #1 \right\}}}

\newcommand{\cI}{{\mathcal{I}}}

\newcommand{\cS}{{\mathcal{S}}}
\newcommand{\cV}{{\mathcal{V}}}

\newcommand{\st}{\,:\,}
\newcommand{\QED}{\hfill $\blacksquare$}

\newcommand{\e}{\epsilon}
\newcommand{\dist}{\mathrm{dist}\,}
\newcommand{\diam}{\mathrm{diam}\,}
\newcommand{\lam}{\lambda}
\newcommand{\Lam}{\Lambda}
\newcommand{\sS}{{\sf S}}

\newcommand{\bhat}{{\hat{b}}}
\newcommand{\chat}{{\hat{c}}}
\newcommand{\lamhat}{{\hat{\lambda}}}

\renewcommand{\Re}{{\mathrm{Re}\,}}
\renewcommand{\Im}{{\mathrm{Im}\,}}

\begin{document}

\maketitle

\begin{abstract}
For the quadratic family $f_{c}(z) = z^2+c$ with $c$ 
in a hyperbolic component  of the Mandelbrot set, 
it is known that every point in the Julia set moves holomorphically. 
In this paper we give a uniform derivative estimate of such a motion when the parameter $c$ converges to a parabolic parameter $\chat$ radially; 
in other words, it stays within a bounded Poincar\'e distance from the internal ray that lands on $\chat$. 
We also show that the motion of each point in the Julia set is uniformly one-sided H\"older continuous at $\chat$ 
with exponent depending only on the petal number. 

This paper is a parabolic counterpart of the authors' paper
``From Cantor to semi-hyperbolic parameters along external rays"
 ({\it Trans. Amer. Math. Soc.} {\bf 372} (2019) pp. 7959--7992).    

\end{abstract}

\tableofcontents

\section{Introduction and main results}

\paragraph{Hyperbolic components.}
Let $\M$ be the {\it Mandelbrot set},
the connectedness locus of the quadratic family
$$
\brac{f_c:z \mapsto z^2 + c}_{c \, \in \, \C}.
$$ 
That is, the Julia set $J(f_c)$ is connected if and only if 
$c \in \M$.
A parameter $c \in \M$ is called {\it hyperbolic} if $f_c$ has a (super-)attracting periodic point. 
Equivalently, there exist positive numbers $\gamma_c$ and $\varepsilon_c$
such that $|Df_c^n(z)| \ge \gamma_c(1 + \varepsilon_c)^n$ 
for any $n \ge 0$ and $z \in J(f_c)$.
The set of hyperbolic parameters in $\M$ is an open subset 
and its connected components are called 
{\it hyperbolic components} of the Mandelbrot set.
(The complement of $\M$ is also called a hyperbolic component, 
but in this paper we only consider those contained in the Mandelbrot set.)

Let $\D$ be the unit disk in $\C$ and 
$\X$ a hyperbolic component of $\M$.
Sullivan and Douady-Hubbard 
(see \cite[Expos\'es XIV \& XIX]{DH} and \cite[Thm.6.5]{Mi1})
gave a {\it uniformization} of $\X$,
which is a canonical homeomorphism
$\Phi=\Phi_{\X}:\Dbar \to \Xbar$ such that 
$\Phi|_{\D}:\D \to \X$ is an conformal isomorphism;
and for any $c=\Phi(\mu)$ with $\mu \in \Dbar-\{1\}$,
the map $f_c$ has a periodic point of multiplier $\mu$
with a common period.
The parameters $\Phi(0)$ and $\Phi(1)$ in $\Xbar$ are called 
the {\it center} and the {\it root} of $\X$ respectively. 
Note that the {\it Poincar\'e distance} in $\X$
is defined by pulling-back the Poincar\'e (hyperbolic) 
metric $|dz|/(1-|z|^2)$ on $\D$ by the isomorphism $\Phi$.

\paragraph{Internal rays and thick internal rays.}
For a given hyperbolic component $\X$ and a given real number $\theta$,
we define the {\it internal ray} $I(\theta)$ of angle $\theta$ by
$$
I(\theta)=I_\X(\theta)
:=\brac{\Phi( r e^{2 \pi i \theta}) \in \X \st 0 \le r <1}. 
$$
The point 
$$
\chat
:=\Phi(e^{2 \pi i \theta}) \in \partial \X
$$ 
is called the {\it landing point} of $I(\theta)$. 
For a given $\delta \ge 0$, 
we define the {\it $\delta$-thick internal ray} 
$\cI(\theta, \delta)$ of angle $\theta$ 
by the closed $\delta$-neighborhood of $I(\theta)$ in $\X$ with respect to the Poincar\'e distance.
 We say the parameter $c$ {\it tends to $\chat$ 
along a thick internal ray} 
if there exists a $\delta \ge 0$ such that
$c$ stays in the $\delta$-thick internal ray $\cI(\theta, \delta)$
while it tends to $\chat$.
It is rather common to say that such a $c$ converges to 
$\chat$ {\it radially} (after McMullen \cite{Mc2}) or {\it non-tangentially}. 
Indeed, for any angle $A_0 \in [0, \pi/2)$,
if $c$ stays in the $\delta(A_0)$-thick internal ray 
with
\begin{equation}
\delta(A_0)
=\frac{1}{2}\log \frac{1+\tan (A_0/2)}{1-\tan (A_0/2)} \in [0,\infty),
\label{eq_A_0}\end{equation}
then by letting $\mu_c:=\Phi^{-1}(c)$ and $\hat{\mu}:=\Phi^{-1}(\chat)$ 
we have 
$$
\abs{
\arg \paren{1-\frac{\mu_c}{\hat{\mu}}}
} \le A_0 
$$
for $c$ sufficiently close to $\chat$.
In other words, $\mu_c$ stays in the Stolz angle at $\hat{\mu} \in \partial \D$ with opening angle $2A_0$. 
(See \cite[p.7]{P} and Figure \ref{fig_thick_ray} in the next section.)

\paragraph{Holomorphic motion of the hyperbolic Julia sets.}
It is well-known that there exists 
a {\it holomorphic motion} 
(\cite{BR, L, Mc1, MSS}) of the Julia sets
over any hyperbolic component $\X$ of $\M$.
Indeed, we have an equivariant holomorphic motion as follows. 
For any base point $\sigma \in \X$,
there exists a unique map $H:\X \times J(f_{\sigma})  \to \C$ such that
\begin{enumerate}[(1)]
\item
$H(\sigma, z) = z$ for any $z \in J(f_{\sigma})$.
\item
For any $ c \in \X$,
the map $z \mapsto {H}(c,z)$ 
is injective on $J(f_{\sigma})$.
\item
For any $z \in J(f_{\sigma})$,
the map $c \mapsto H(c,{z})$ is holomorphic on $\X$.
\item
For any $c \in \X$,
the map $h_c(z):=H(c,z)$ satisfies $h_c(J(f_{\sigma}))=J(f_c)$ and 
$f_c \circ h_c = h_c \circ f_{\sigma}$ on $J(f_{\sigma})$.
\end{enumerate}
See \cite[\S 4]{Mc1} for more details.
In this paper we choose the center $\sigma:=\Phi(0)$ of $\X$ as the base point of the motion. 
We are concerned with boundary behavior of
such an equivariant holomorphic motion of the Julia set $J(f_{\sigma})$
when $c \in \X$ tends to some $\chat \in \partial \X$
along a thick internal ray.

\paragraph{Parabolic parameters.}
Now suppose that the angle $\theta$ of the internal ray 
$I(\theta)$ is a rational number.
Then for the landing point $\chat$ of $I(\theta)$,
$f_\chat$ has a parabolic periodic point, 
that is, a periodic point 
whose multiplier is a root of unity.
We say such a parameter $\chat$ is {\it parabolic},
and a parabolic periodic point $\bhat$ of $f_\chat$
has {\it $q$ petals} if 
the local dynamics of $f_\chat^{k}$ 
near $\bhat$ is of the form
$\zeta \mapsto \zeta+\zeta^{q+1}+O(\zeta^{q+2})$
for some $k$ in an appropriate local coordinate.
In our setting, it is known that $\bhat$ has $q$ petals 
if and only if the multiplier of $\bhat$
is a primitive $q$-th root of unity 
(since $f_\chat$ is quadratic and has only one critical point in $\C$).

\begin{figure}[htbp]
\begin{center}
\includegraphics[width=.6\textwidth]{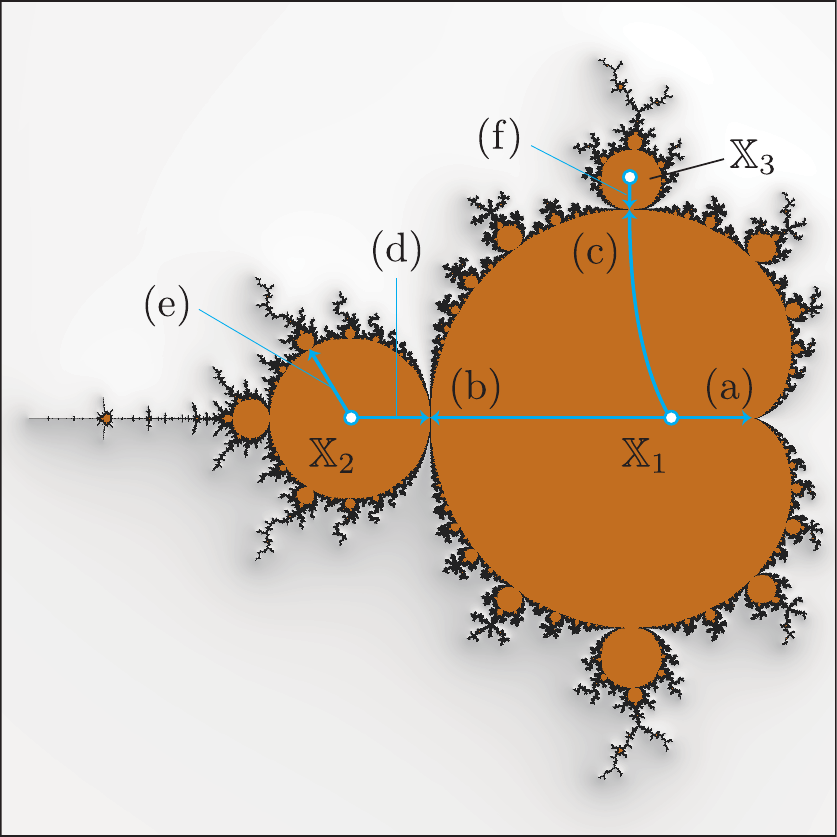}
\end{center}
\caption{Internal rays (a) -- (f) of the Mandelbrot set.}
\label{fig_mandelbrot}
\end{figure}

\begin{figure}[htbp]
\begin{center}
\includegraphics[height=60mm, bb =0 0 1146 894]{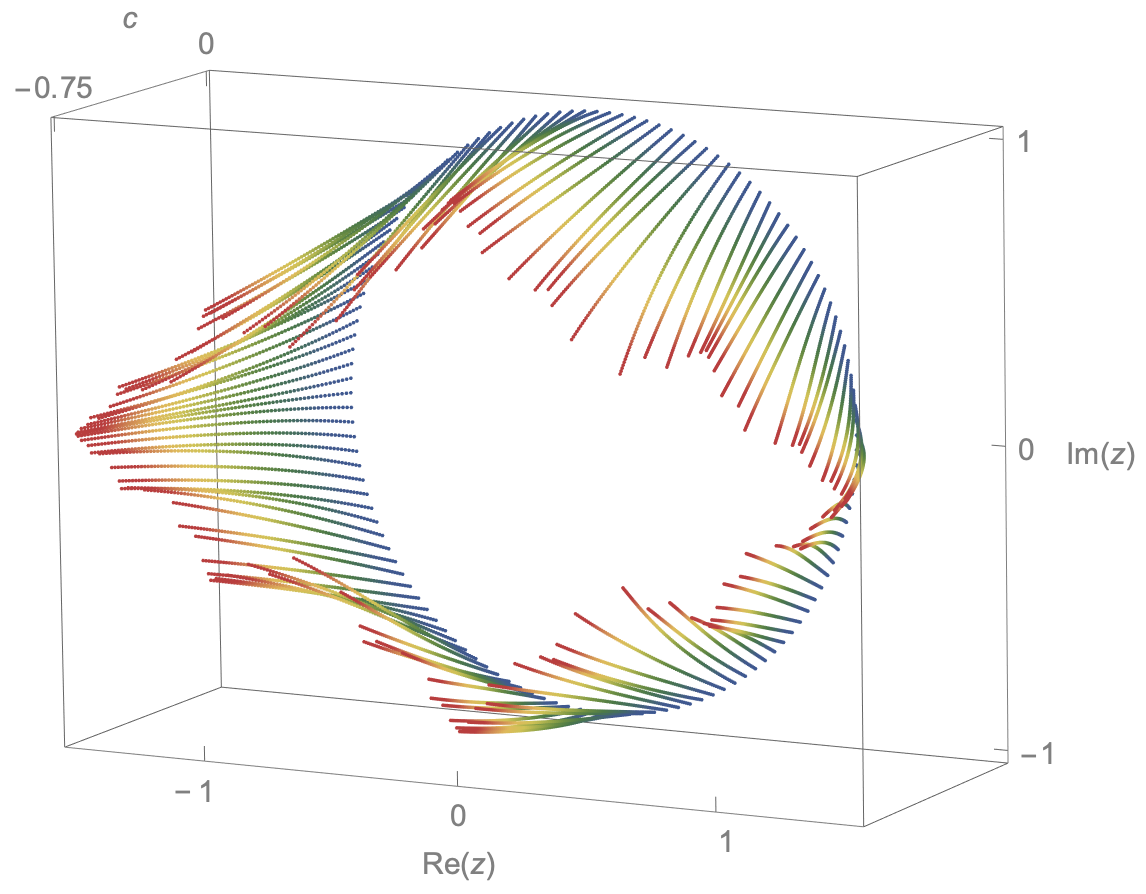}
\includegraphics[height=60mm, bb =0 0 1344 1362]{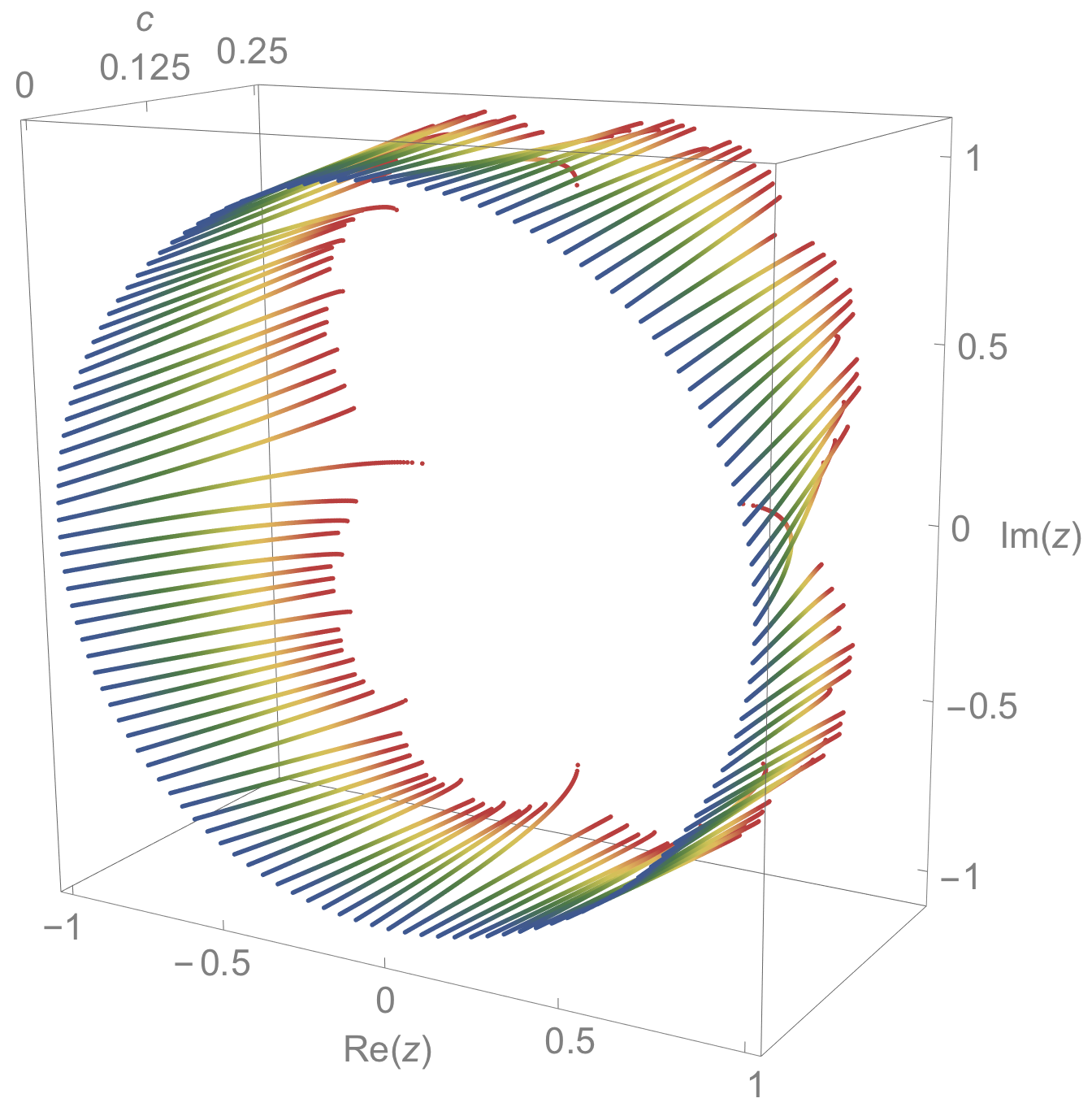}
\end{center}
\caption{Real analytic motion of the preimages 
of the repelling fixed point along the internal rays $I(1/2)$ (left) and $I(0)$
(right) in $\X_1$.}
\label{fig_motion}
\end{figure}

\begin{figure}[htbp]
\begin{center}
\fbox{\includegraphics[width=.18\textwidth, bb = 0 0 600 600]{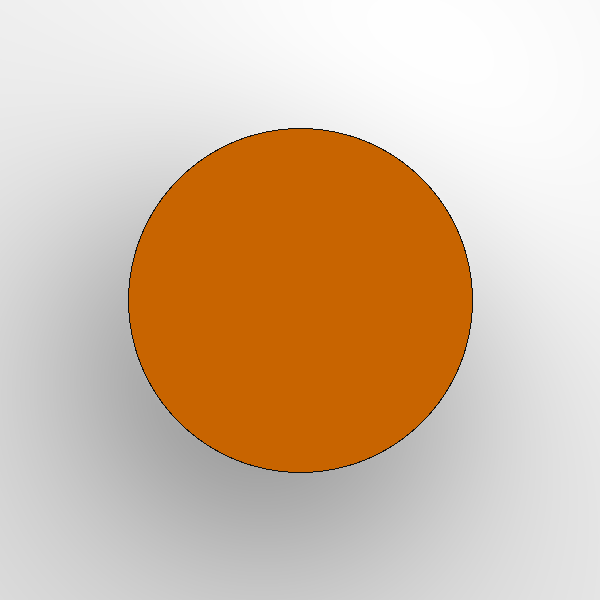}}
\fbox{\includegraphics[width=.18\textwidth, bb = 0 0 600 600]{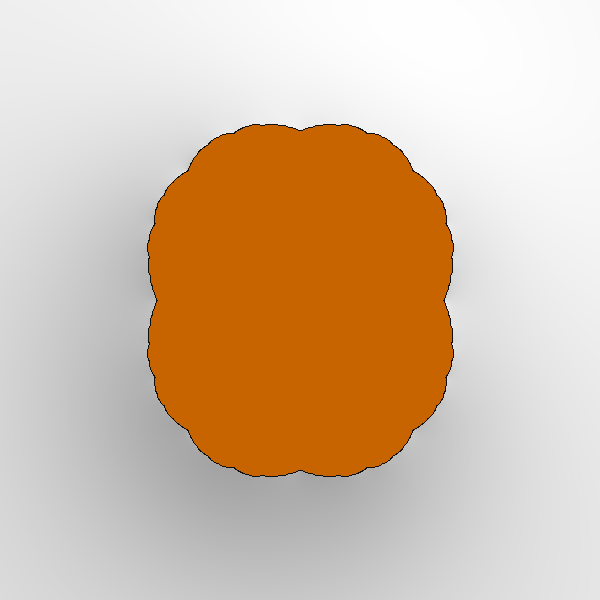}}
\fbox{\includegraphics[width=.18\textwidth, bb = 0 0 600 600]{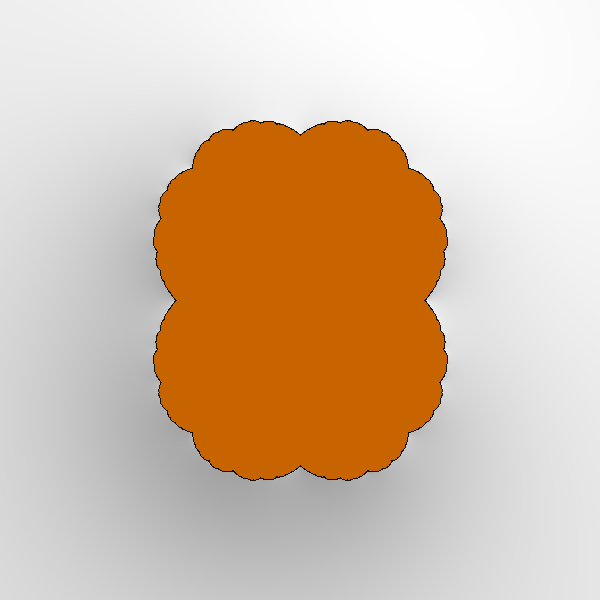}}
\fbox{\includegraphics[width=.18\textwidth, bb = 0 0 600 600]{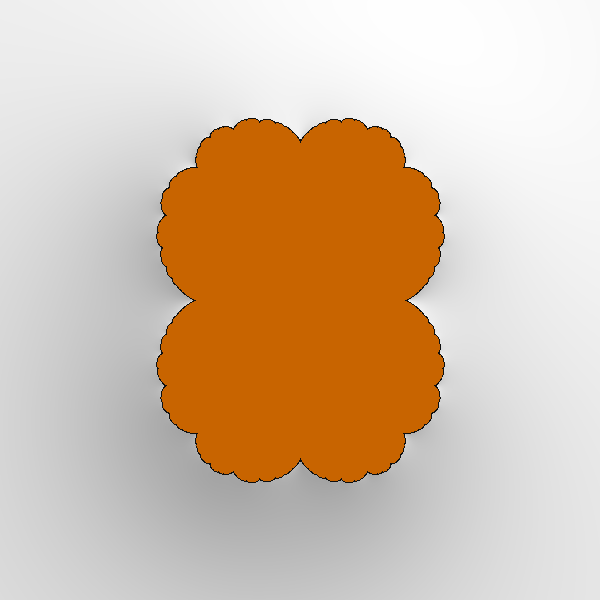}}
\fbox{\includegraphics[width=.18\textwidth, bb = 0 0 600 600]{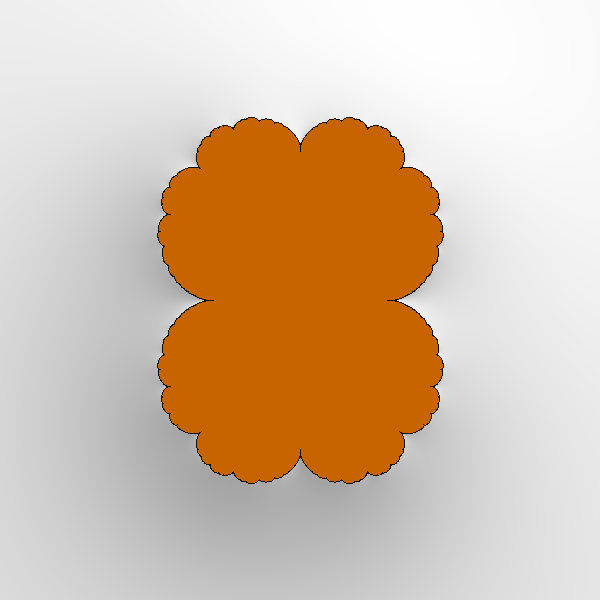}}
\\[.5em]
(a): $I(0)$ in $\X_1$\\[.5em]
\fbox{\includegraphics[width=.18\textwidth, bb = 0 0 600 600]{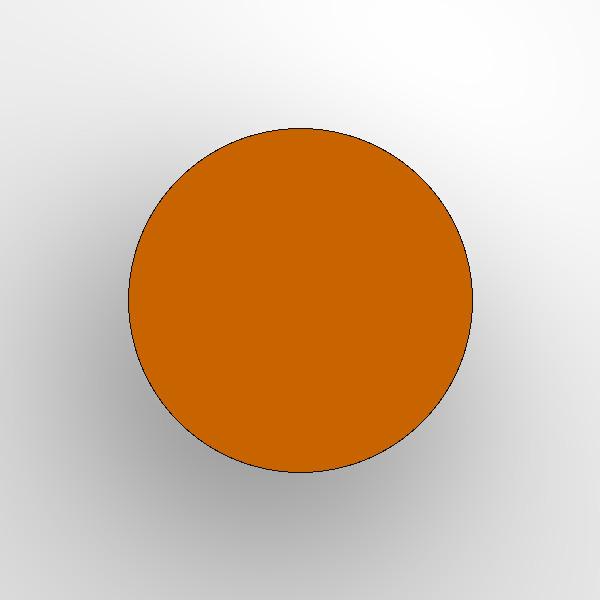}}
\fbox{\includegraphics[width=.18\textwidth, bb = 0 0 600 600]{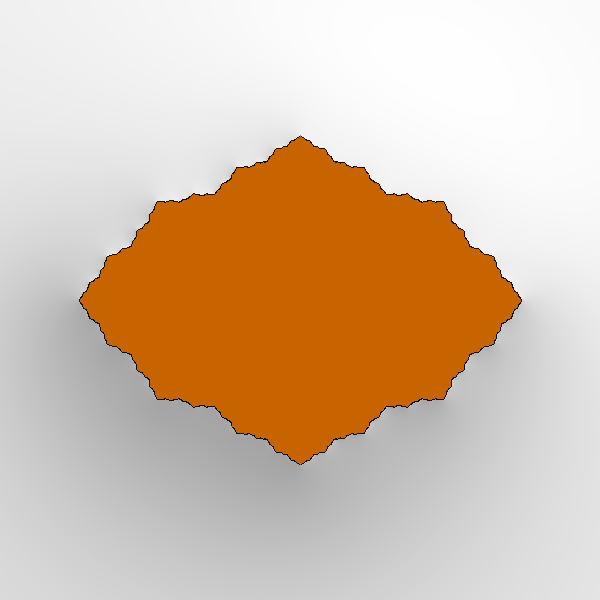}}
\fbox{\includegraphics[width=.18\textwidth, bb = 0 0 600 600]{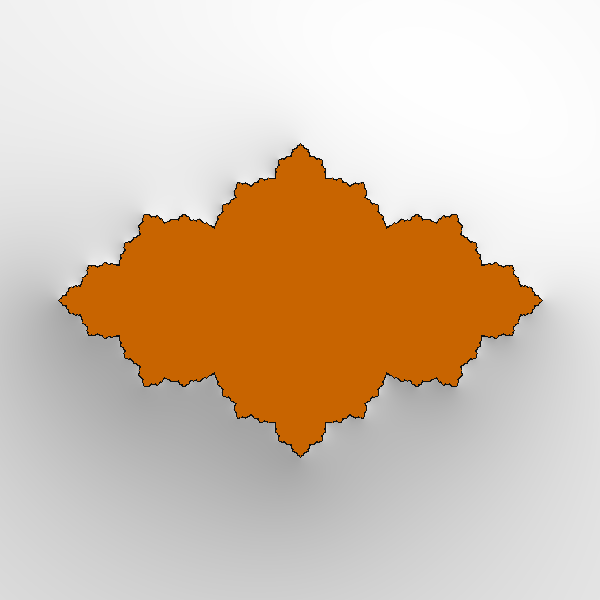}}
\fbox{\includegraphics[width=.18\textwidth, bb = 0 0 600 600]{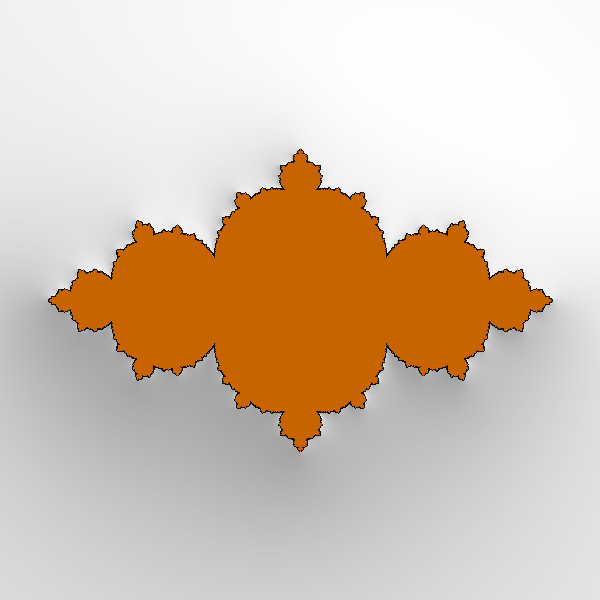}}
\fbox{\includegraphics[width=.18\textwidth, bb = 0 0 600 600]{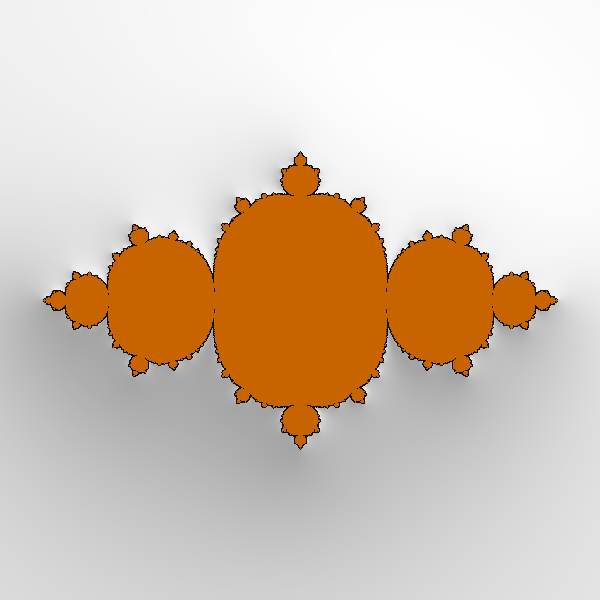}}
\\[.5em]
(b): $I(1/2)$ in $\X_1$\\[.5em]
\fbox{\includegraphics[width=.18\textwidth, bb = 0 0 600 600]{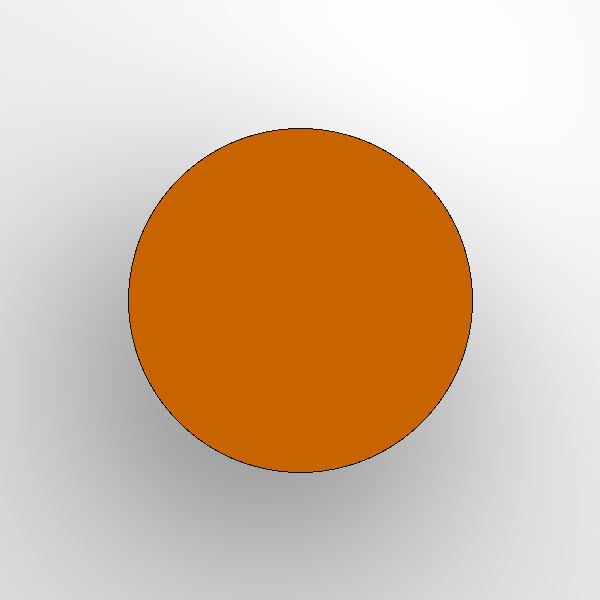}}
\fbox{\includegraphics[width=.18\textwidth, bb = 0 0 600 600]{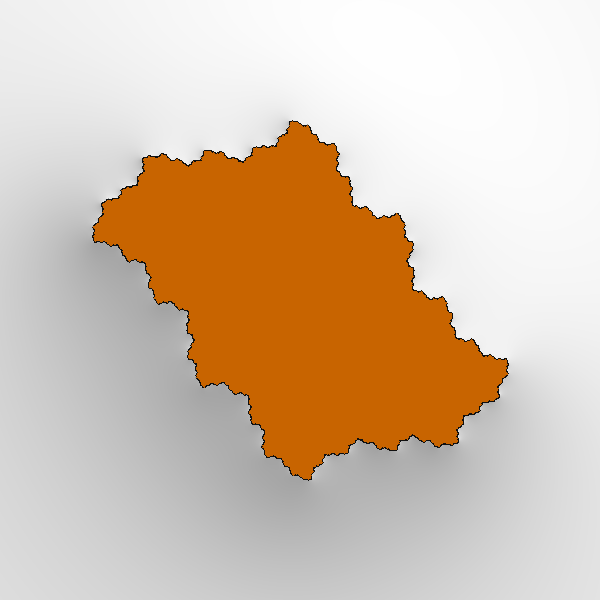}}
\fbox{\includegraphics[width=.18\textwidth, bb = 0 0 600 600]{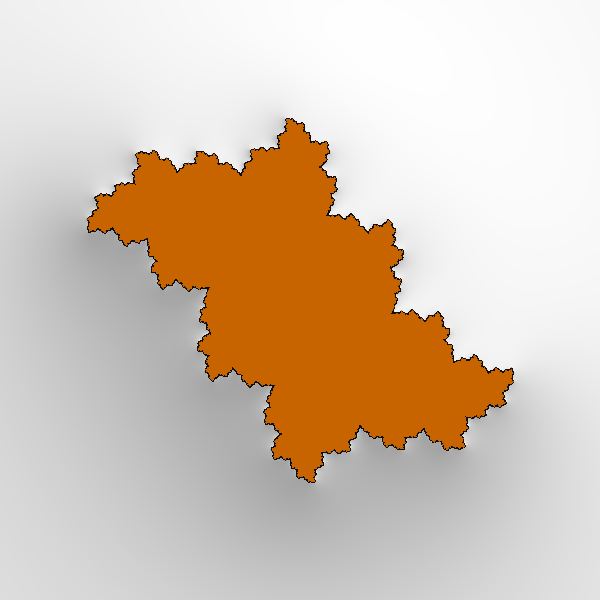}}
\fbox{\includegraphics[width=.18\textwidth, bb = 0 0 600 600]{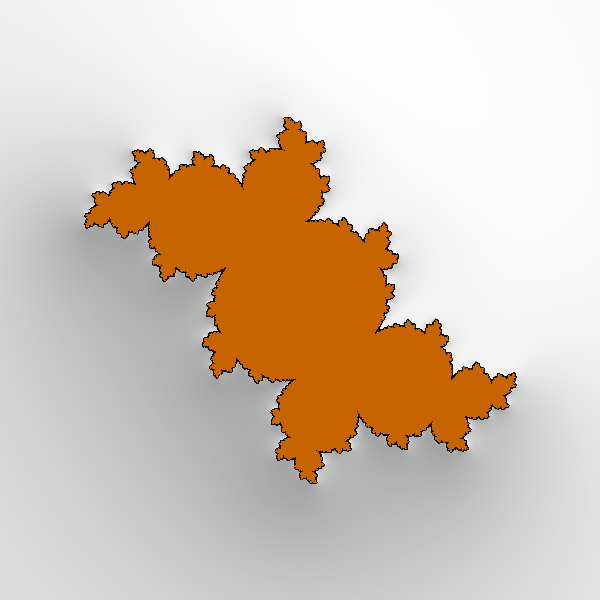}}
\fbox{\includegraphics[width=.18\textwidth, bb = 0 0 600 600]{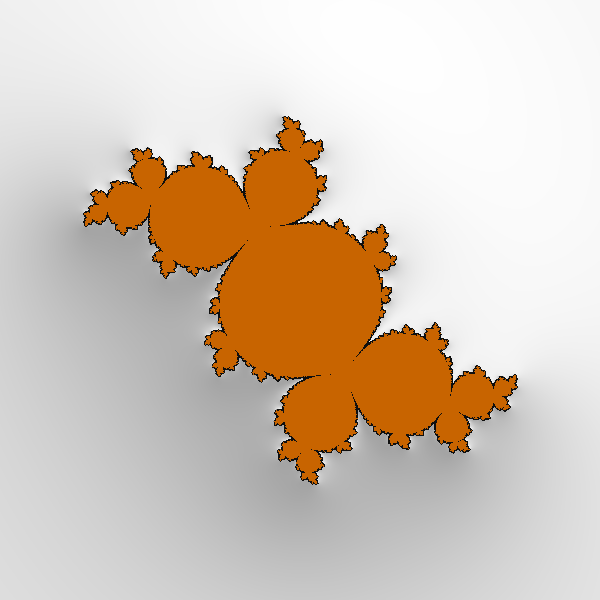}}
\\[.5em]
(c): $I(1/3)$ in $\X_1$\\[.5em]
\fbox{\includegraphics[width=.18\textwidth, bb = 0 0 600 600]{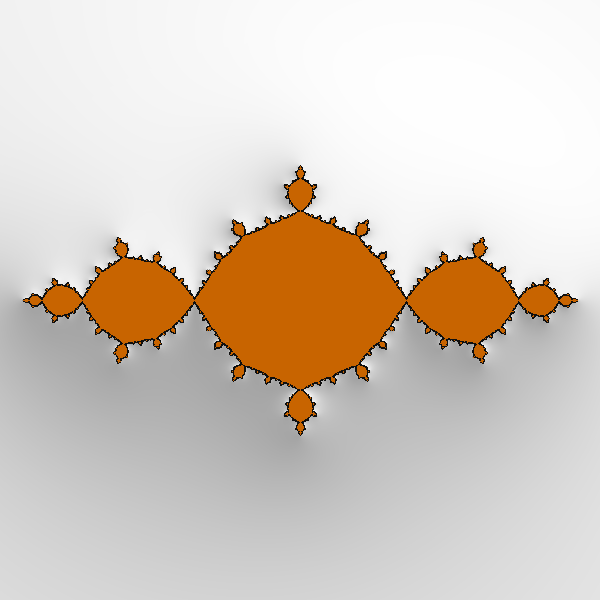}}
\fbox{\includegraphics[width=.18\textwidth, bb = 0 0 600 600]{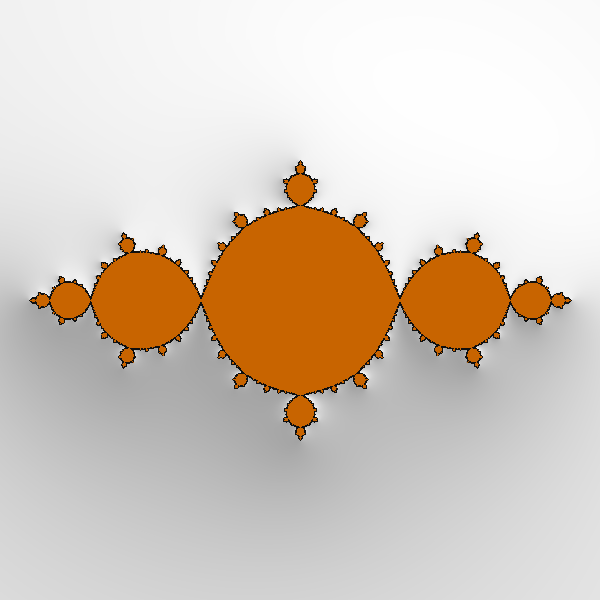}}
\fbox{\includegraphics[width=.18\textwidth, bb = 0 0 600 600]{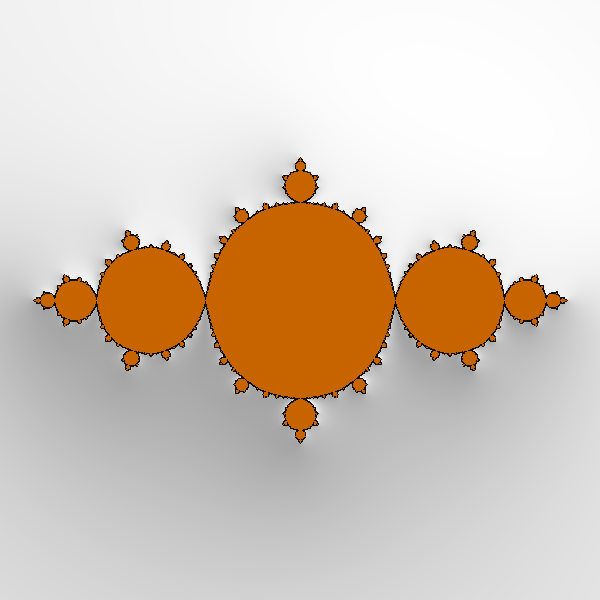}}
\fbox{\includegraphics[width=.18\textwidth, bb = 0 0 600 600]{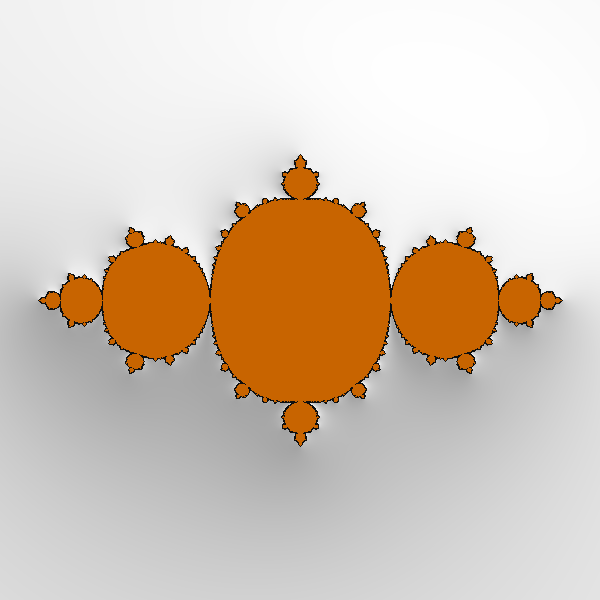}}
\fbox{\includegraphics[width=.18\textwidth, bb = 0 0 600 600]{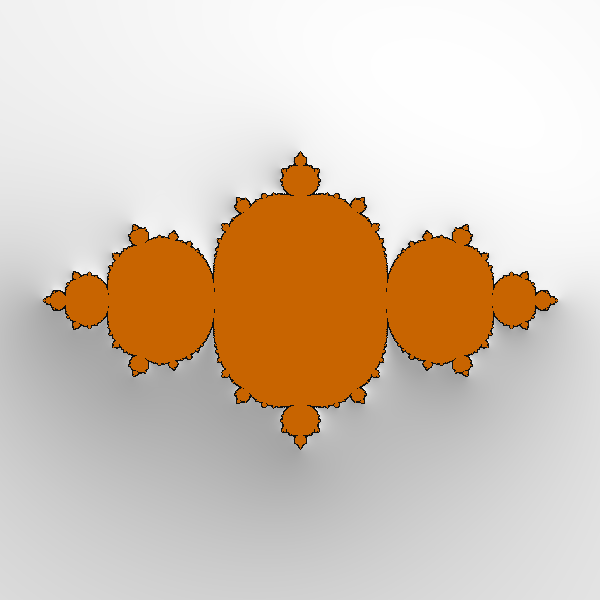}}
\\[.5em]
(d): $I(0)$ in $\X_2$\\[.5em]
\fbox{\includegraphics[width=.18\textwidth, bb = 0 0 600 600]{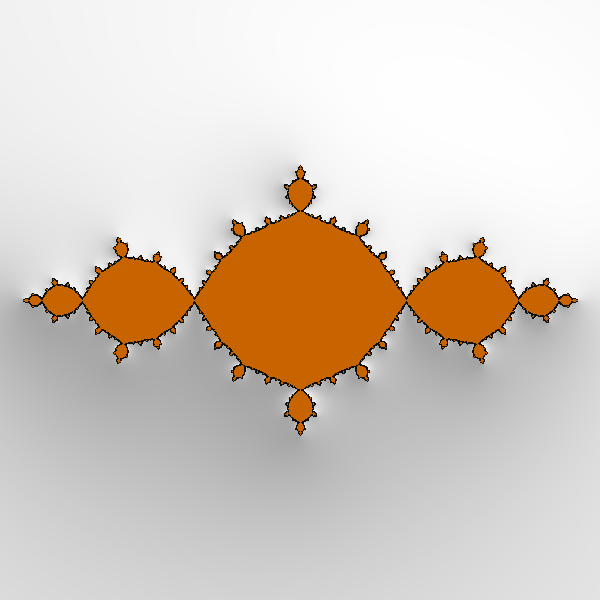}}
\fbox{\includegraphics[width=.18\textwidth, bb = 0 0 600 600]{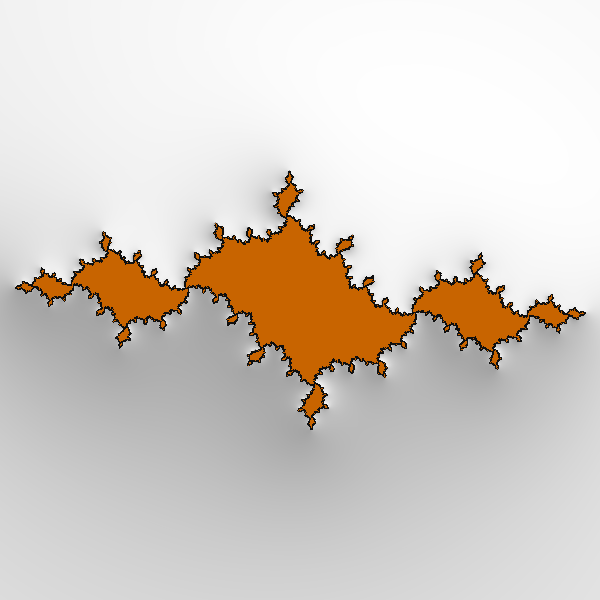}}
\fbox{\includegraphics[width=.18\textwidth, bb = 0 0 600 600]{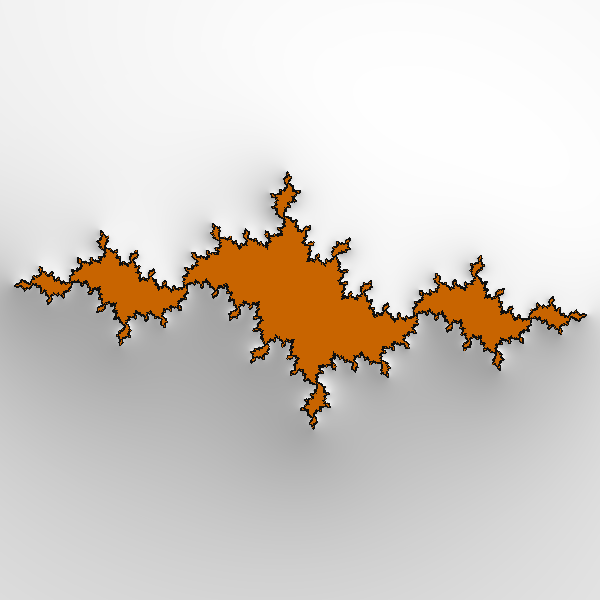}}
\fbox{\includegraphics[width=.18\textwidth, bb = 0 0 600 600]{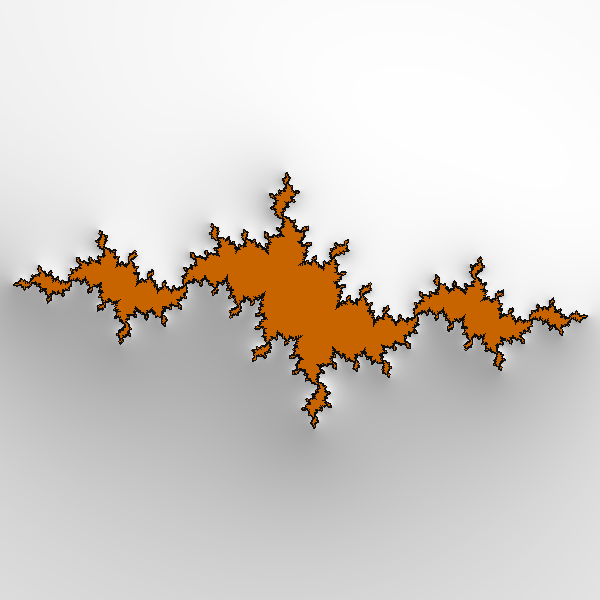}}
\fbox{\includegraphics[width=.18\textwidth, bb = 0 0 600 600]{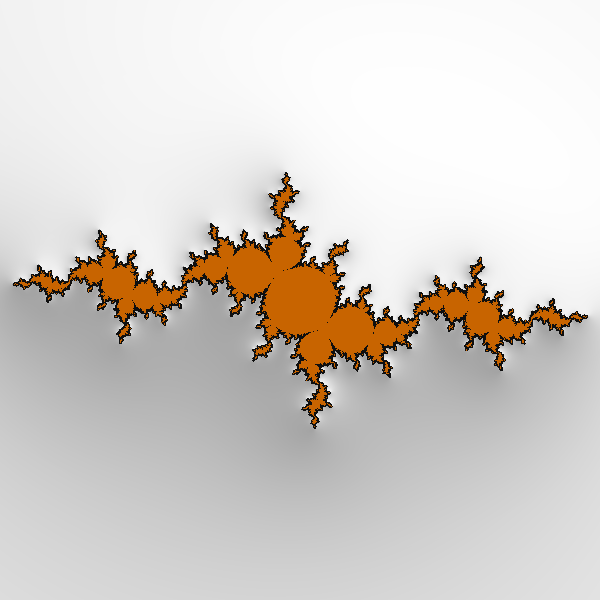}}
\\[.5em]
(e): $I(1/3)$ in $\X_2$\\[.5em]
\fbox{\includegraphics[width=.18\textwidth, bb = 0 0 600 600]{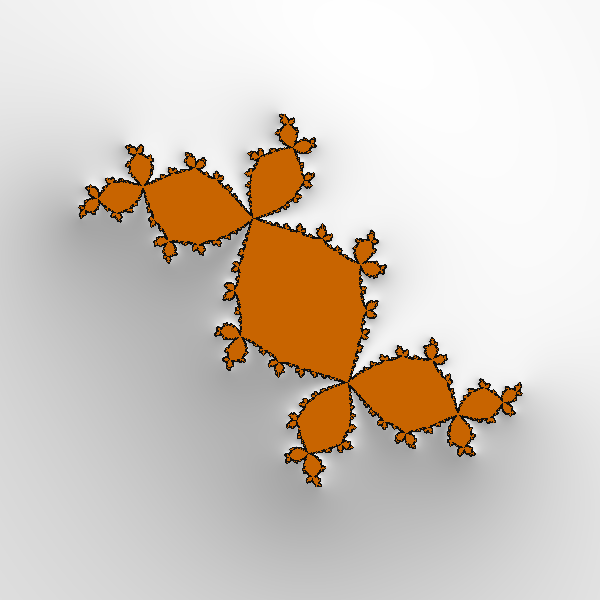}}
\fbox{\includegraphics[width=.18\textwidth, bb = 0 0 600 600]{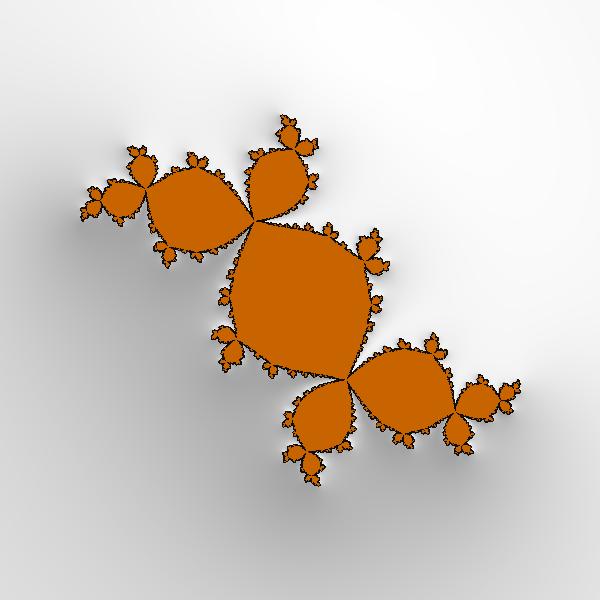}}
\fbox{\includegraphics[width=.18\textwidth, bb = 0 0 600 600]{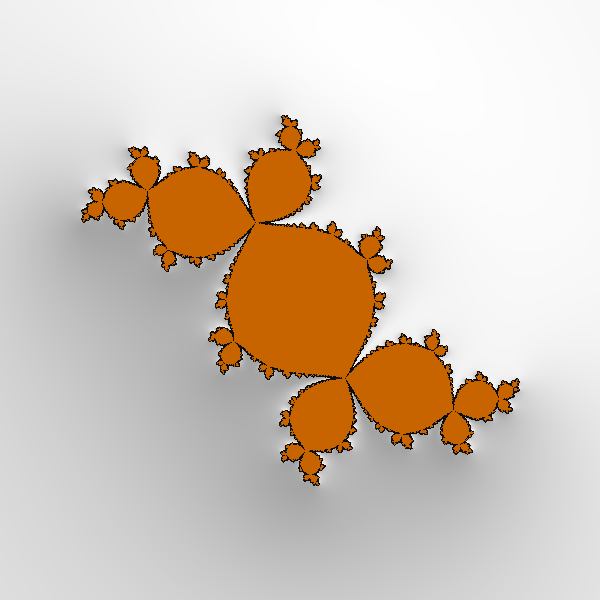}}
\fbox{\includegraphics[width=.18\textwidth, bb = 0 0 600 600]{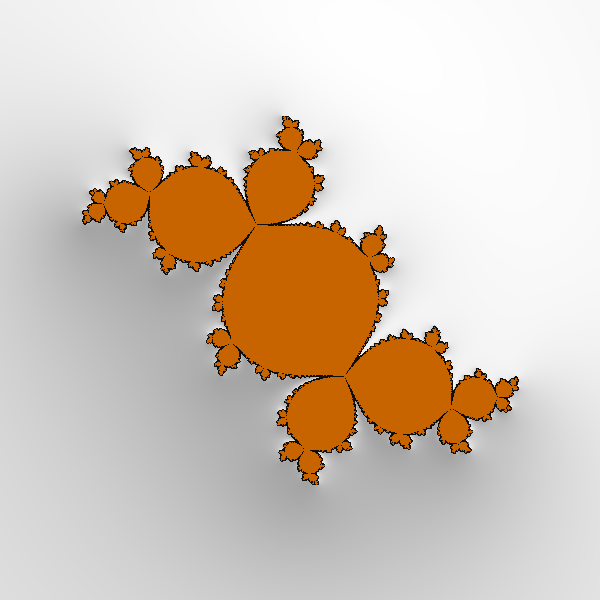}}
\fbox{\includegraphics[width=.18\textwidth, bb = 0 0 600 600]{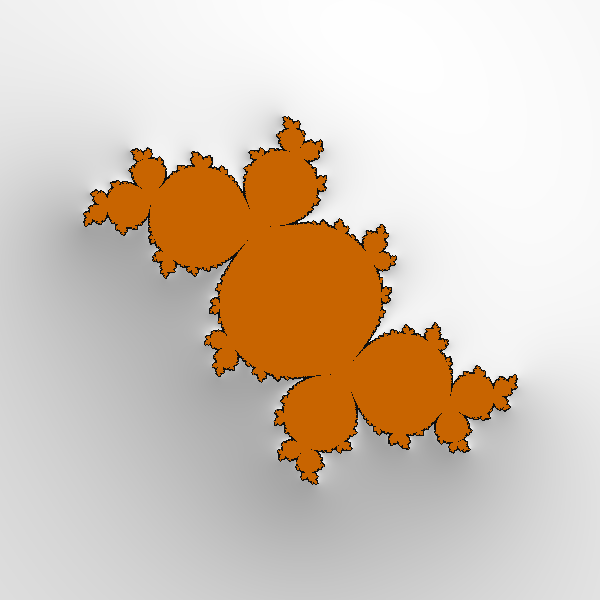}}
\\[.5em]
(f): $I(0)$ in $\X_3$
\end{center}
\caption{
Holomorphic motions along the internal rays depicted in Figure \ref{fig_mandelbrot}.}
\label{fig_tables}
\end{figure}

\begin{eg}[Period one, the main cardioid]
For the hyperbolic component $\X_1$ containing $0$ 
({\it the main cardioid}), 
the map $\Phi=\Phi_{\X_1}$ is explicitly given by 
$\Phi(\mu):= -\mu^2/4+\mu/2~(\mu \in \Dbar)$
and $f_{\Phi(\mu)}$ has a fixed point with multiplier $\mu$.
The internal ray of angle $\theta$ is given by 
$
I_{\X_1}(\theta)
=\brac{-r^2 e^{4 \pi \theta i}/4+r e^{2 \pi \theta i}/2 \st
0 \le r<1}.
$ 
In Figure \ref{fig_mandelbrot}, $I_{\X_1}(\theta)$ 
for $\theta=0, 1/2$ and $1/3$ are 
depicted as paths (a), (b), and (c) respectively. The corresponding holomorphic motions along $I_{\X_1}(0)$ and $I_{\X_1}(1/2)$ are illustrated in Figure \ref{fig_motion} and in Figures \ref{fig_tables}(a) and \ref{fig_tables}(b).  The motion  along $I_{\X_1}(1/3)$ is depicted in Figure  \ref{fig_tables}(c).
\end{eg}

\begin{eg}[Period two]
Similarly the hyperbolic component $\X_2$ 
containing $-1$ consists of hyperbolic parameters $c$
such that $f_c$ has an attracting cycle of period two.
The map $\Phi=\Phi_{\X_2}$ is explicitly given by 
$\Phi(\mu):= \mu/4-1$ for $\mu \in \Dbar$,
and $f_{\Phi(\mu)}$ with $\mu \in \Dbar-\{1\}$ 
has a periodic point of period two with multiplier $\mu$.
The internal ray of angle $\theta$ is given by 
$
I_{\X_2}(\theta)
=\brac{r e^{2 \pi \theta i}/4-1 \st 0 \le r<1}.
$
In Figure \ref{fig_mandelbrot}, $I_{\X_2}(\theta)$ 
for $\theta=0$ and $1/3$ are 
depicted as paths (d) and (e) respectively.
The internal ray $I_{\X_2}(0)$ lands 
at the root $\chat=\Phi(1)=-3/4$ of $\X_2$,
where the map $f_\chat$ has a parabolic fixed point of 
multiplier $-1$ that has two petals.
Note that $\chat$ is the landing point of 
another internal ray $I_{\X_1}(1/2)$ of $\X_1$. See Figures \ref{fig_tables}(d) and \ref{fig_tables}(e) for the motions  along $I_{\X_2}(0)$ and $I_{\X_2}(1/3)$.
\end{eg}

\begin{eg}[Period three]
There is a hyperbolic component $\X_3$ attached to the main cardioid $\X_1$
that consists of parameters $c$ with $\mathrm{Im}\, c >0$
for which $f_c$ has an attracting cycle of period three.
The internal ray $I_{\X_3}(0)$
(depicted as path (f) in Figure \ref{fig_mandelbrot})
joins the center (so-called ``rabbit") 
and the root $\chat=\Phi(1)$ (``the fat rabbit"),
where
the map $f_\chat$ has a parabolic fixed point of 
multiplier $e^{2 \pi i /3}$ that has three petals.  (See Figure \ref{fig_tables}(f)  for the holomorphic motion  along $I_{\X_3}(0)$.) 
Again $\chat$ is the landing point of 
another internal ray $I_{\X_1}(1/3)$ of $\X_1$.
\end{eg}

\paragraph{Main results.}
Let $\X$ be a hyperbolic component in the Mandelbrot set and
$\sigma$ be its center.
For any $z_\ast$ in $J(f_{\sigma})$, 
the map $c \mapsto z(c):=H(c,z_\ast)$ is holomorphic over $\X$. 
Let $\chat \in \partial \X$ be the landing point of the 
internal ray $I_{\X}(\theta)$ of rational angle $\theta$.
Our main theorem states that 
the speed of $z(c)=H(c,z_\ast)$ is uniformly bounded by a function of $|c-\chat|$ as $c$ tends to $\chat$ along a thick internal ray: 

\begin{thm}[Main Theorem] \label{thm_main}
Suppose that $f_\chat$ has a parabolic periodic point with $q$ petals
and $c$ tends to $\chat$ along a thick internal ray $\cI(\theta,\delta)$ in $\X$.
Then there exists a constant $K>0$ depending only on $\chat$ and $\delta$
such that for any $z = z(c)\in J(f_c)$, 
the point $z(c)$ moves holomorphically with 
$$
\abs{\frac{d}{dc}z(c)}
\le 
\frac{K}{\abs{\, c-\chat \,}^{1-1/Q}},
$$
where $Q=\max\{2, q\}$. 

\end{thm}

\medskip
By this theorem we obtain one-sided H\"older continuity of the holomorphic motion along thick internal rays landing on parabolic parameters:

\begin{thm}[One-sided H\"older Continuity]
\label{thm_HolderContinuity}
Under the same assumption as Theorem \ref{thm_main} above,
the point $z=z(c)$ in $J(f_c)$ tends to a limit $z(\chat)$ in $J(f_\chat)$
as $c$ tends to $\chat$ along the thick internal ray $\cI(\theta,\delta)$. 
Moreover, there exists a constant $K'$ 
depending only on $\chat$ and $\delta$ such that 
\begin{equation}
\label{eq_thm1.2}
|z(c)-z(\chat)| \le K' |\, c-\chat \, |^{1/Q}
\end{equation}
for any $c$ in $\cI(\theta,\delta)$.
\end{thm}
\medskip
As an immediate consequence the holomorphic motion of 
each point $z(c) \in J(f_c)$ lands when $c$ moves along 
the internal ray $I(\theta)=\cI(\theta,0)$. 
This theorem yields a precise description of 
the degeneration of the dynamics on the Julia sets
along the internal rays of rational angles:

\begin{thm}[Pinching Semiconjugacy]\label{thm_H2P}
Under the same assumption as the theorems above,
the conjugacy $H(c,\cdot)=h_c: J(f_{\sigma}) \to J(f_c)$ 
converges uniformly to a semiconjugacy $h_{\chat}:J(f_{\sigma}) \to J(f_\chat)$
from $f_{\sigma}$ to $f_{\chat}$ as $c$ tends to $\chat$ along the thick internal ray $I(\theta, \delta)$.
Moreover, $h_\chat$ satisfies the following:
\begin{enumerate}[\rm (1)]
\item 
If $\chat$ is the root of $\X$,
then $h_{\chat}$ is injective and thus a conjugacy.
\item 
If $\chat$ is not the root of $\X$ (hence $q \ge 2)$,
then the preimage $h_\chat^{-1}(\{w\})$ of any $w \in J(f_\chat)$
consists of one or $q$ distinct points, and 
the latter holds if and only if $w$ eventually lands on a 
parabolic periodic point of $f_\chat$. 
\item
The semiconjugacy $\eta_c:=h_\chat \circ h_c^{-1}:J(f_c) \to J(f_\chat)$
satisfies 
\begin{equation}\label{eq_H2P}
|\eta_c(z)-z| \le K' |c-\chat|^{1/Q}
\end{equation}
for any $c$ in the thick internal ray $\cI(\theta,\delta)$.
\end{enumerate}
\end{thm}

\medskip
By (3) of this theorem we obtain:

\begin{cor}[Hausdorff Convergence]\label{cor_HausdorffConvergence}
The Hausdorff distance between $J(f_c)$ and $J(f_\chat)$ 
is $O(|c-\chat|^{1/Q})$ as $c$ tends to $\chat$ along a thick internal ray.
\end{cor}

\begin{remark} 
{}~
\begin{itemize}
\item
These results are parabolic counterparts of the authors' results in
\cite{CK1} about parameter rays (external rays) landing on semi-hyperbolic parameters of the Mandelbrot set.
\item
For any $c\in \X$ and $z =z(c) \in J(f_c)$, 
we have
\begin{equation}\label{eq_CK1_Prop3.1}
\abs{\frac{d}{dc}{z}(c)} 
\le 
\frac{\, 1+\sqrt{1+6\,|c|\,}\,}{\dist(c, \partial \X)}
\end{equation}
by Proposition 3.1 in \cite{CK1}. 
However, this will only give 
$$
\abs{\frac{d}{dc}{z}(c)}=O\paren{\frac{1}{|c-\chat|}}
$$
as $c$ tends to $\chat$.
\item
In \cite{CK2}, the authors showed that 
for any $c \in [0, 1/4)=I_{\X_1}(0)$ and $z =z(c) \in J(f_c)$, 
we have an optimal estimate
$$
\abs{\frac{d}{dc}{z}(c)} 
\le 
\frac{1}{2\sqrt{1/4-c}}.
$$
In particular, the Hausdorff distance between 
$J(f_c)$ and $J(f_{1/4})$ is exactly $\sqrt{1/4-c}$.
\item
The existence of the semiconjugacy in Theorem \ref{thm_H2P} 
and the Hausdorff convergence of the Julia sets in 
Corollary \ref{cor_HausdorffConvergence}
are previously shown in a more general context 
by the second author \cite{K1} and McMullen \cite{Mc2} respectively. 
The novel part of our results 
is the quantitative estimate $O(|c-\chat|^{1/Q})$.
\end{itemize}
\end{remark}

\paragraph{Structure of the paper.}
In Section \ref{sec:access_cycle}, we define a parametrization of $c$ with a complex parameter $t\in \C$ such that $c=c_t$ converges to $\chat$ along a thick internal ray as $t\to 0$.  Also in Section \ref{sec:access_cycle}, we state three propositions,  Propositions \ref{prop_c_and_epsilon}, \ref{prop_near_parabolic} and \ref{prop_leaving_U_0}, which concern the local dynamics of $f_c$ in a neighborhood $U_0$ of a parabolic point of $f_\chat$ when $c$ is near $\chat$, and will be employed to prove Theorems \ref{thm_HolderContinuity} and \ref{thm_H2P} 
as well as some lemmas in the paper. Then, we introduce the notion of ``S-cycle" to describe how an orbit of $f_c$ repeatedly (infinite or finite times or never) enters and leaves a fixed subset of $U_0$. In Section \ref{sec:Proof of the main theorem}, by assuming Lemmas A, B and C, we prove our main theorem, Theorem \ref{thm_main}. It is well-known (for example \cite[\S 3.2]{Mc1}) that the Julia set is expanding with respect to the hyperbolic metric on $\C-P(f_c)$, where $P(f_c)$ denotes the postcritical set of $f_c$. In order to estimate the expansion of the Julia set with respect to the Euclidean metric, we give an estimate of the distance  between $z\in J(f_c)$ and $P(f_\chat)$ in Lemma D in Section \ref{sec:Hyperbolic metrics} for $z$ not too close to the parabolic cycle of $f_\chat$. Lemma A is proved in Section \ref{sec: proof_of_A} by assuming another two lemmas, Lemmas G and H. Both lemmas 
  rely on local dynamics of perturbed parabolic cycle. We prove Lemma B in Section \ref{sec:Proof_B}, and Lemma G in Section \ref{sec:Proof_G}.  Section \ref{sec_parabolic} is devoted to the proofs of  Propositions \ref{prop_c_and_epsilon} and
\ref{prop_near_parabolic}. We use a branched coordinate  to prove Proposition \ref{prop_leaving_U_0} in Section \ref{sec:Proof_of_Leaving}. We also employ the branched coordinate to prove Lemma H in Section \ref{sec:Proof_H}. Then, using some results presented in Section \ref{sec:Proof_H}, we are able to prove Lemma D in Section \ref{sec:Proof_of_Lemma_D}. Some arguments in the proofs of Lemmas H and D are used to prove Lemma C in Section \ref{sec:Proof_Of_Lemma_C}. Finally, in Section \ref{sec:together} we prove Theorems \ref{thm_HolderContinuity} and \ref{thm_H2P} simultaneously.

\fboxsep=0pt
\fboxrule=.5pt

\section{Radial access condition and S-cycles} \label{sec:access_cycle}

In this section we introduce the notion of S-cycles for a given orbit in the Julia set. 
The idea of S-cycles was introduced in \cite{CK1} to describe orbits 
that repeatedly come close to the postcritical set.
Here we present a modified version where the postcritical set is replaced by the parabolic cycle. 

\paragraph{Notation.}
We start with some notation and the terminology that will be used in what follows.

\begin{itemize}
\item
Let $\N$ denote the set of positive integers. We denote the set of non-negative integers by $\N_0:=\{0\} \cup \N$.
\item
Let $\D(a,r)$ denote the disk in $\C$ centered at $a$ and of radius $r>0$.
When $a=0$ we denote it by $\D(r)$. 
\item
For non-negative variables $X$ and $Y$, by $X \asymp Y$ we mean there exists 
an implicit constant $C>1$ independent of $X$ and $Y$ such that $X/C \le Y \le CX$.  
\item 
When we say ``for any $X \ll 1$" it means that 
``for any sufficiently small $X>0$". 
More precisely,  we mean there exists an implicit constant $C>0$ such that $0<X <C$. 
\end{itemize}

\paragraph{Hyperbolic components and internal rays.}
Let $\chat \in \partial \M$ be a parabolic parameter having
a parabolic periodic point $\bhat$ of period exactly $p$.
Let $\hat{\lambda}:=Df_{\chat}^p(\bhat)$ 
be the multiplier of this cycle, 
and assume that it is a primitive $q$-th root of unity. 
We specify an internal ray $I(\theta)=I_\X(\theta)$ of 
the hyperbolic component $\X$ that lands at $\chat$ as follows.
(See \cite{DH} or \cite{Mi1} for details on the hyperbolic components
of $\M$.) 

\paragraph{Case 1.}
If $q =1$, then there is only one hyperbolic component $\X$
such that $\chat = \Phi_{\X}(1) \in \partial \X$,
where $\Phi_{\X}$ is the uniformizing map of $\X$. 
Hence by letting $\theta=0$ the internal ray $I(\theta)$ 
of $\X$ lands at $\chat$.

\paragraph{Case 2.}
If $q \ge 2$, then there are exactly two hyperbolic components
$\X^-$ and $\X^+$ such that
\begin{itemize}
\item
$\partial \X^- \cap \partial \X^+=\{\chat \}$.
\item
$\chat = \Phi_{\X^-}(\hat{\lambda})$ and $\chat = \Phi_{\X^+}(1)$, 
where $\Phi_{\X^\pm}$ is the uniformizing map of $\X^{\pm}$.
\end{itemize} 
Hence $\X$ can be either $\X^-$ or $\X^+$,
and the case of $q \ge 2$ is divided into two sub-cases.
\begin{itemize}
\item
{\bf Case 2$^-$:}
If $\X=\X^-$, then we let $\theta= (\arg \hat{\lambda})/(2 \pi)$; and
\item
{\bf Case 2$^+$:}
If $\X=\X^+$, then we let $\theta=0$
\end{itemize}
in such a way that the internal ray 
$I(\theta)=I_\X(\theta)$ of $\X$ lands at $\chat$.

Note that $\chat$ is the root of $\X$ if and only if it is as Case 1 or Case 2$^+$. Hence, 
$$
\theta=\left\{\begin{array}{ll}0 
& \text{for {Case 1} or {Case 2$^+$},}\\
[.7ex] (\arg \hat{\lambda})/(2 \pi)
 & \text{for {Case 2$^-$}.}\end{array}\right.
$$

\begin{eg}\label{eg_X}
Case 1 holds when $\chat=1/4$ with $\X=\X_1$ (Example 1),
or $\chat=-7/4$ with $\X$ whose center is a unique real parameter 
$\sigma<0$ with $f_{\sigma}^3(0)=0$ (``the airplane").
If $\hat{c}=-3/4$, 
Case 2$^-$ holds when $\X^-=\X_1$, and
Case 2$^+$ holds when $\X^+=\X_2$ (Example 2).  
\end{eg}

\begin{figure}[htbp]
\begin{center}
\includegraphics{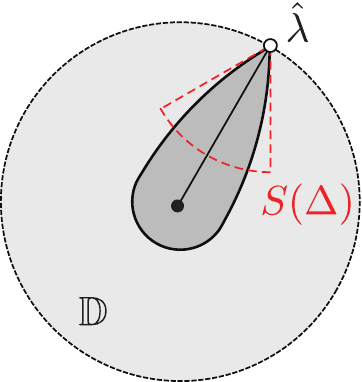}
\end{center}
\caption{The Stolz angle $S(\Delta)$ with opening angle $2A_0$ equal to 
$\pi/3$ at $\hat{\lam} \in \partial \D$. The dark gray region is the image of $\delta$-thick internal ray $\cI(\theta, \delta)$ under $\Phi^{-1}$ with $\delta = \delta(A_0)$.} 
\label{fig_thick_ray}
\end{figure}

\paragraph{Radial convergence and thick internal rays.}
Let $A_0$ and $T_0$ be constants with 
$0 \le A_0 <\pi/2$ and $0<T_0 < 2\cos A_0$, and let
$$
\Delta =\Delta(A_0, T_0):=\{ t \in \C \st 0< |t| \le T_0,~|\arg t\,| \le A_0\}. 
$$
The {\it Stolz angle} at $\hat{\lambda}$
with opening angle $2 A_0$ is given by
\begin{equation}
S(\Delta) :=\{ \mu =(1-t) \hat{\lambda} 
\in \D \st t \in \Delta\} \subset \D. \label{Stolz_mu_t}
\end{equation}
Let $\Phi:=\Phi_\X$.
If the parameter $c \in \X$ tends to $\chat$ 
satisfying $\Phi(c) \in S(\Delta)$,
we say $c \to \chat$ {\it radially} after McMullen \cite{Mc2}.

For a given $\delta$-thick internal ray $\cI(\theta,\delta) \subset \X$ of angle $\theta$,
one can easily check that 
$$
\Phi^{-1}(\cI(\theta,\delta)) \cap \D(\hat{\lambda}, T_0) 
\subset S(\Delta)
$$
if $\delta \le \delta(A_0)$, which is given in \eqref{eq_A_0}. See Figure \ref{fig_thick_ray}.
Hence in what follows it is enough to consider the parameters 
of the form $c=\Phi(\mu)$ with $\mu \in S(\Delta) \cup \{ \hat{\lambda} \}$.

\begin{remark} 
Conversely, the Stolz angle $S(\Delta)$ is contained in 
$\Phi^{-1}(\cI(\theta,\delta))$ with $\delta=\delta(A_0)+1$ 
by taking a sufficiently small $T_0$. 
This implies that the convergence in a thick internal ray 
is equivalent to the radial convergence for some angle. 
\end{remark}

\paragraph{Parametrization and notation.}
For a technical reason, instead of \eqref{Stolz_mu_t},  it is convenient to re-parametrize 
$\mu \in S(\Delta)$ by
$$
\mu=
\mu_t:=\left\{\begin{array}{ll}
1-q\,t
& \text{in {Case 1} or {Case 2$^+$}}\\[.7ex] 
(1-t/q )\hat{\lambda}
& \text{in {Case 2$^-$}}\end{array}\right.
$$
for $t \in \Delta=\Delta(A_0, T_0)$
with sufficiently small $T_0$.

\paragraph{The radial access condition.}
In what follows, by 
$$
c \approx \chat \quad \text{or} \quad c=c_t \approx \chat
$$
we mean the parameter $c$ is of the form 
$$
c=c_t=\Phi(\mu_t)
$$
for some $t \in \Delta=\Delta(A_0,T_0)$, 
where we take a smaller $T_0$ 
in the definition of $\Delta$ if necessary. 
We say such a parameter $c$ {\it satisfies 
the radial access condition} 
or {\it $c$ is in a thick internal ray}.

\paragraph{Perturbation of parabolic points.}
Let $c_0:=\chat$ and $\lam_0:=\hat{\lambda}$.
In Section \ref{sec_parabolic}, we will show the following two propositions under the radial access condition:

\begin{prop}\label{prop_c_and_epsilon}
The function $t \mapsto c=c_t$
is holomorphic on $\Delta$,
and there exists a constant $B_0 \neq 0$ such that
\begin{itemize}
\item {\bf Case 1 ($q =1$):} 
$c_t=\chat + B_0 t^2 +O(t^3)$ 
\item {\bf Case 2$^\pm$ ($q \ge 2$):} 
 $c_t=\chat + B_0 t +O(t^2)$ 
\end{itemize}
as $t \in \Delta$ tends to $0$.
In particular, 
we have $\sqrt{|c_t-\chat|}\asymp |t|$ 
or $|c_t-\chat|\asymp |t|$ 
for $t \in \Delta$
according to Case 1 or Case 2$^\pm$.
\end{prop}

\begin{prop}\label{prop_near_parabolic}
There exists a continuous map $t \mapsto b_t$ defined for $t \in \Delta \cup \{ 0 \}$ such that
\begin{enumerate}[\rm (1)]
\item
$b_0=\hat{b}$ and $b_t$ is a periodic point of $f_{c_t}$ with 
the same period $p$ as $\bhat$. 
\item
Let $\lambda_t:= Df_{c_t}^p(b_t)$.
Then there exist two families of holomorphic local coordinates 
$\brac{\zeta=\varphi_{t}(z)}_{t \,\in \,\Delta \cup \{ 0 \}}$ 
and 
$\brac{w=\psi_{t}(z)}_{t \,\in \,\Delta \cup \{ 0 \}}$ 
defined on a disk $\hat{U}:=\D(\bhat, \hat{R})$
such that: For each $z \in \hat{U}$, 
$\varphi_{t}(z)$ and $\psi_{t}(z)$ 
are holomorphic in $t \in \Delta$ and continuous at $t=0$; 
$\varphi_{t}(b_{t})=\psi_{t}(b_{t})=0$; and 
\begin{align}
\varphi_{t} \circ f_{c_t}^{p} \circ \varphi_{t}^{-1}(\zeta)
&= \lambda_t \zeta +\zeta^{q+1} +O(\zeta^{2q+1}), ~\text{and}
\label{eq_local_coordinates1}
\\[.5em]
\psi_{t} \circ f_{c_t}^{pq} \circ \psi_{t}^{-1}(w)
&= \lambda_t^q w \,(1+w^q +O(w^{2q})).
\label{eq_local_coordinates1}
\end{align}
In particular, both $D\varphi_{t}^{-1}(0)$ and $D\psi_{t}^{-1}(0)$ 
are uniformly bounded away from zero.
\item
In Case 1 or Case 2$^+$, 
$b_t$ is repelling for $t \in \Delta$, and the multiplier satisfies
$$
\lambda_t=(1+t/q+o(t))\, \hat{\lambda}.
$$
Moreover, there are $q$ distinct attracting fixed points 
$\alpha_t^1,\, \cdots,\, \alpha_t^q$ of $f_{c_t}^{pq}$
satisfying $Df_{c_t}^{pq}(\alpha_t^{j})=\mu_t=1-qt$ 
and $\left(\psi_{t}(\alpha_t^{j})\right)^q = -t+o(t)$ for $j=1,\cdots, q$.
\item
In Case 2$^-$, $b_t$ is attracting for $t \in \Delta$ and the multiplier satisfies 
$$
\lambda_t=\mu_t=(1-t/q) \,\hat{\lambda}.
$$
Moreover, there are $q$ distinct repelling fixed points 
$\beta_t^1,\, \cdots,\, \beta_t^q$ of $f_{c_t}^{pq}$
satisfying $Df_{c_t}^{pq}(\beta_t^{j})=1+q t+o(t)$ 
and $\left(\psi_{t}(\beta_t^{j})\right)^q= t +o(t)$ for $j=1,\cdots, q$.
\end{enumerate}
\end{prop}

\begin{figure}[htbp]
\includegraphics[width=.95\textwidth]{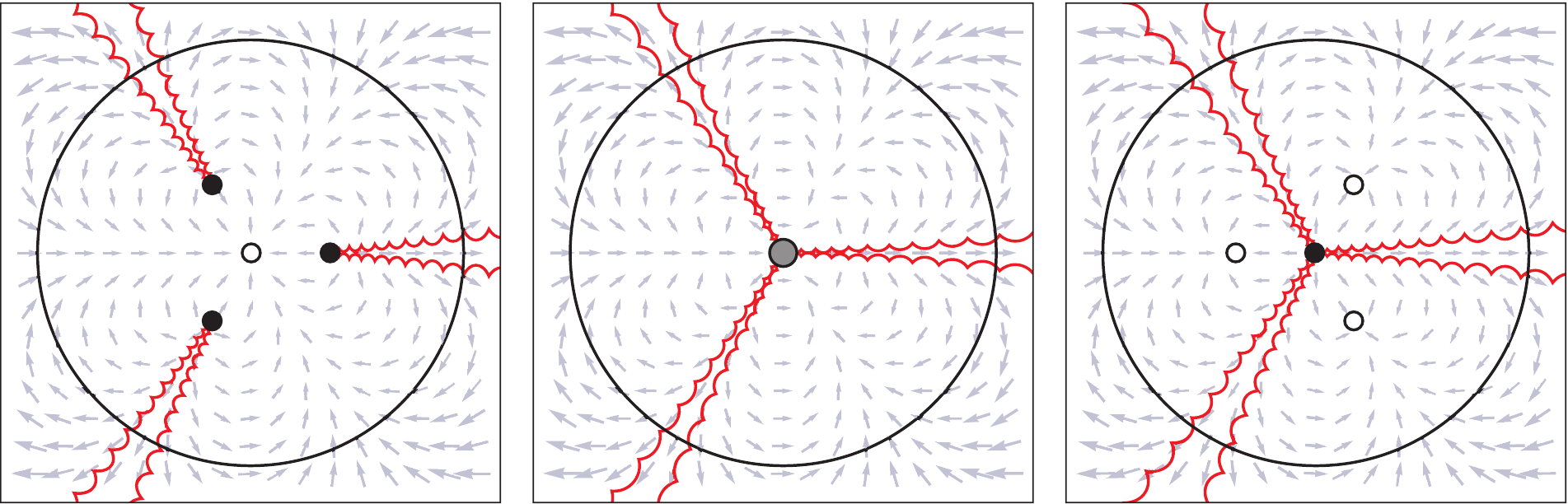}
\caption{
The middle frame depicts the local dynamics
of $f_{\chat}^{pq}$ near $\hat{U}$ (indicated by the big circle)
that can be mildly perturbed in two ways:
Case 2$^-$ on the left, and Case 2$^+$ on the right.
The white, black, and gray dots indicate attracting, repelling, and parabolic fixed points respectively. The red spiky curves indicate the Julia sets.}
\label{fig_perturbation}
\end{figure}

The local dynamics of $f_{c_t}^{pq}$ observed as 
\eqref{eq_local_coordinates1} behaves quite similar 
to that of $f_{\chat}^{pq}$. 
See Figure \ref{fig_perturbation}.
In particular, 
in the domain $\hat{U}=\D(\bhat, \hat{R})$ 
of $\psi_{t}$, 
the map $f_{c_t}^{pq}$ has exactly
\begin{itemize}
\item
one repelling fixed point $b_t$ in Case 1 or Case 2$^+$; and
\item
$q$ repelling fixed points in Case 2$^-$ that are symmetrically arrayed
near $b_t$.    
\end{itemize}

\paragraph{Definition of $U_0$.}
We fix a small $R_0 \in (0, \hat{R})$ such that the disk
$$
U_0:=\D(\bhat, R_0) \quad (\Subset \hat{U})
$$
possesses the  property that  $f_\chat^p:U' \to U_0$ is univalent for each connected component $U'$ of $f_\chat^{-p}(U_0)$.  
Such an $R_0$ exists because the orbit $f_\chat^n(0)$ of $0$ 
(the critical point) 
keeps a definite distance from the parabolic cycle 
for $0 \le n \le p$.
The next proposition will be proved in Section \ref{sec:Proof_of_Leaving}:

\begin{prop}\label{prop_leaving_U_0}
For each parameter $c_t$ with $t \in \Delta$ 
if $z_0 \in J(f_{c_t}) \cap U_0$ 
is none of the repelling fixed points of $f_{c_t}^{pq}$ described as above, then the orbit $z_{kpq}:=f_{c_t}^{kpq}(z_0)~(k \in \N_0)$
leaves $U_0$ for some $k >0$.
\end{prop}

We choose some $\xi \in (0, 1]$ such that 
$$
\dist(0, f_\chat^{k}(U_0)) \ge \frac{3\xi}{4}
$$ 
for any integer $k$ with $-p \le k \le p$.
Note that by continuity, we have
$$
\dist(0, f_c^{k}(U_0)) \ge \frac{\xi}{2}
$$ 
for any $c \approx \chat$ (taking a smaller $T_0$ 
in the definition of $\Delta$ if necessary) 
and any $k$ with $-p \le k \le p$.

\medskip

\begin{remark}\label{rem_xi} 
We will frequently use the following property: 
{\it If $z \in f_c^{k}(U_0)$ for some 
$c \approx \chat$  and 
$k$ with $-p \le k \le p$, then $|Df_c(z)|=2|z| \ge \xi$.}
\end{remark}

\paragraph{Definition of $V_0$ and $\cV(c)$.}
Now we take a small $\nu \in (0,R_0)$ and let
$$
V_0:=\D(\bhat, \nu) \quad(\Subset U_0).
$$ 
We also define
$$
\cV(c):=\bigcup_{j=0}^{p-1} f_c^{-j}(V_0)
$$
for each $c \approx \chat$. 
By Remark \ref{rem_xi} above, $f_c(z) \in \cV(c)$ implies $|Df_c(z)| \ge \xi$ for $c \approx \chat$.

\paragraph{S-cycles.}
For $c \approx \chat$, let $z_0$ be any point in the Julia set 
$J(f_c)$.
The orbit $z_n:=f_c^n(z_0)~(n \in \N_0)$
may land on $V_0$ ($\Subset U_0$),
and leave $U_0$ by Proposition \ref{prop_leaving_U_0}
(unless it lands exactly on the repelling cycle), 
then it may come back to $V_0$ again. 
To describe the behavior of such an orbit, 
we introduce the notion of ``S-cycle" for the orbit of $z_0$,
where ``S" indicates that orbit stays near the ``singularity" 
of the hyperbolic metric $\rho(z)|dz|$ on the complement of the postcritical set of $f_\chat$ to be defined 
in Section \ref{sec:Hyperbolic metrics}.

\paragraph{Definition {\normalfont (S-cycle)}.}
A {\it finite S-cycle} of the orbit $z_n=f_c^n(z_0)~(n \in \N_0)$ 
is a finite subset of $\N_0$ of the form
$$
\sS=\brac{n \in \N_0 \st M \le n <M'}=[M,M') \cap \N_0
$$
with the following properties:
\begin{itemize}
\item[(S1)]
$z_M \in V_0$, and if $M>0$ then $z_{M-1} \notin V_0$.
\item[(S2)]
There exists a minimal $m \ge 1$ such that 
for $n=M+mpq$, 
$z_{n-pq} \in U_0$ but $z_n \notin U_0$.
\item[(S3)]
$M'=M + mpq +L$ for some $L \in [1,\infty)$
such that $z_n \notin V_0$ 
for $n=M+mpq+i~(0 \le i < L)$
and $z_{M'} \in V_0$.
\end{itemize}

An {\it infinite S-cycle} $\sS$ 
of the orbit $z_n=f_c^n(z_0)~(n \in \N_0)$ 
is an infinite subset of $\N_0$ of the form
$$
\sS=\brac{n \in \N_0 \st M \le n <\infty}=[M,\infty) \cap \N_0, 
$$
satisfying either 
\begin{itemize}
\item 
Type {(I):}
(S1), (S2), and
\begin{itemize}
\item[(S3)']
$z_n \notin V_0$ for all $n \ge M+mpq$; 
\end{itemize}
\end{itemize}
or
\begin{itemize}
\item Type {(II):}
(S1) and
\begin{itemize}
\item[(S2)']
$z_{M+kpq} \in V_0$ for any $k \in \N$. 
Equivalently, $z_M$ is a repelling periodic point of period $p$ or period $pq$ in $V_0$ (by Proposition \ref{prop_leaving_U_0}).
\end{itemize}
\end{itemize}
By an {\it S-cycle} we mean a finite or infinite S-cycle.
In both cases, we denote them by
$\sS=[M, M')$ or $\sS=[M, \infty)$ 
for brevity.

\medskip

\begin{figure}[htbp]
\begin{center}
\includegraphics{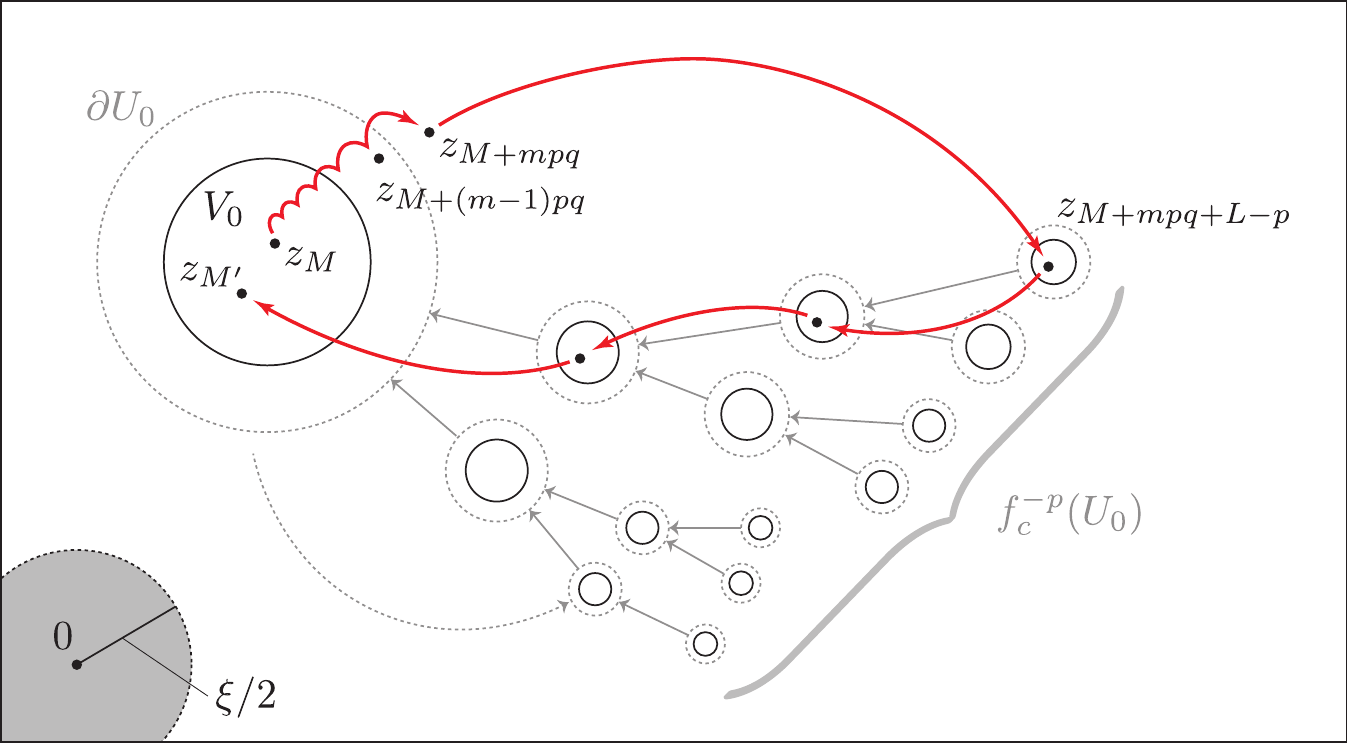}
\end{center}
\caption{A finite S-cycle $[M, M')$ with $M'=M+mpq+L$, $p=3$.
The dotted arrow indicates there is another connected component 
of $f_c^{-p}(U_0)$ intersecting with $U_0$ which is not drawn.
}\label{fig_S-cycle}
\end{figure}

\begin{remark} \label{rem_L}
We may assume without loss of generality 
that $L$ of the finite S-cycle is at least $p$
by shrinking the radius of the disk $U_0$.
Indeed, after the orbit $z_n~(n\in \N_0)$ 
leaves $U_0$ when $n=M+mpq$,
the orbit follows the parabolic cycle for a while 
and it cannot land immediately on $V_0$ by the local dynamics
near the perturbed cycle. 
See Figure \ref{fig_S-cycle}.
\end{remark}

\paragraph{Decomposition of the orbit by S-cycles.}
For a given orbit $z_n=f_c^n(z_0)~(n \in \N_0)$ of $z_0 \in J(f_c)$,
the set $\N_0$ of indices is uniquely decomposed 
by using finite or infinite S-cycles in one of the following three types:
\begin{itemize}
\item
The first type is of the form
\begin{equation}\label{eq_decomposition}
\N_0=[0,M_1) \sqcup \sS_1 \sqcup \sS_2 \sqcup \cdots,
\end{equation}
where $z_n \notin V_0$ for $n \in [0,M_1)$ and  
$\sS_k:=[M_k, M_{k+1})$ is a finite S-cycle for each $k \ge 1$.
\item
The second type is of the form
\begin{equation}\label{eq_decomposition2}
\N_0=[0,M_1) \sqcup \sS_1 \sqcup \sS_2 \sqcup 
\cdots \sqcup \sS_{k_0} \sqcup \emptyset \sqcup \emptyset \sqcup \cdots
\end{equation}
with $k_0 \ge 1$, where  $z_n \notin V_0$ for $n \in [0,M_1)$;  
$\sS_k:=[M_k, M_{k+1})$ is a finite S-cycle for each $1\le k \le k_0-1$;
$\sS_{k_0}=[M_{k_0}, \infty)$ is an infinite S-cycle;
and $\sS_k=\emptyset$ for $k \ge k_0+1$.
\item
The third type is  
\begin{equation}\label{eq_decomposition3}
\N_0=[0,M_1) \sqcup \emptyset \sqcup \emptyset \sqcup \cdots
\end{equation}
 with $M_1=\infty$
and $\sS_k=\emptyset$ for $k \ge 1$,
where $z_n \notin V_0$ for all $n \in \N$. 
\end{itemize}
In the first and second types it is possible that $M_1=0$ and $[0, M_1)$ is empty.

\section{Proof of the main theorem assuming three lemmas}
\label{sec:Proof of the main theorem}

\paragraph{The derivative formula.}
Let $z_\ast$ be any point in $J(f_{\sigma})$ (where $\sigma$ is the center of $\X$), and consider its motion
$z=z(c)=H(c,z_\ast) \in J(f_c)$ for $c \in \X$. 
The estimate of the main theorem is based on the following formula:

\begin{prop}[The Derivative Formula]
\label{prop_Derivative_Formula}
For any $c \in \X$ and $z =z(c) \in J(f_c)$, 
we have 
$$
\frac{d}{dc}{z}(c) = 
-\sum_{n = 1}^\infty\frac{1}{Df_c^{{n}}({z}(c))}.
$$
\end{prop}
See \cite[Proposition 3.2]{CK1} for the proof.
Since $f_c$ is hyperbolic, the convergence of the series above is absolute and it is enough to show 
$$
\left|\sum_{n = 1}^\infty\frac{1}{Df_c^{{n}}({z}(c))}\right| 
\le \frac{K}{ |c-\chat|^{1-1/Q}}
$$
for some constant $K>0$ independent of $c$ 
in a thick internal ray, 
where $Q:= \max\{2,\,q \}$ and $q$ is the petal number of 
the parabolic fixed point $\bhat$.

Now we present three principal lemmas about S-cycles
that are valid for sufficiently small $\nu$ 
(the radius of $V_0$) under the radial access condition 
$c =c_t \approx \chat$.
That is, we only consider 
$c=c_t$ with $t \in \Delta=\Delta(A_0,T_0)$
as in the previous section.

\paragraph{Lemma A.}
{\it 
There exists a constant $K_{\mathrm A}>0$
such that for any $c \approx \chat$,
any $z_0 \in J(f_c)$, and for any 
S-cycle $\sS=[M,M')$ of the orbit 
$z_n=f_c^n(z_0)~(n \in \N_0)$, 
we have 
\begin{equation}
\left|\sum_{i = 1}^{M'-M} \frac{1}{Df_c^i(z_M)}\right|
\le 
\frac{K_{\mathrm A}}{ |c-\chat|^{1-1/Q}},
\end{equation}
where we set $M'-M:=\infty$ if $M'=\infty$. 
}

\paragraph{Lemma B.}
{\it 
There exists a constant $K_{\mathrm B} > 0$
such that for any $c \approx \chat$
and any $M \le \infty$,
if $z_0 \in J(f_c)$ satisfies 
$z_n \notin V_0$ for any $n \in [0, M)$,
then we have \begin{equation}
\sum_{i = 1}^{M} \frac{1}{|Df_c^i(z_0)|}
\le K_{\mathrm B}.
\end{equation}
}

\medskip
An immediate consequence of Lemma B is:

\begin{cor}
\label{cor_of_Lemma B}
For $c \approx \chat$, 
if the orbit of $z=z(c) \in J(f_c)$ by $f_c$ 
never lands on $V_0$,
then the derivative satisfies
\begin{equation}
\abs{\frac{d}{dc}z(c)}
\le 
\sum_{n = 1}^{\infty} 
\frac{1}{|Df_c^n(z(c))|}
\le K_{\mathrm B}.
\end{equation}
\end{cor}

\paragraph{Lemma C {\normalfont (S-cycles Expand Uniformly)}.}
{\it 
There exists a constant $\Lambda> 1$ 
such that for any $c \approx \chat$,
any $z_0 \in J(f_c)$, and for any 
finite S-cycle $\sS=[M,M')$ of the orbit $z_n=f_c^n(z_0)~(n \in \N_0)$, 
we have 
\begin{equation}
|Df_c^{M'-M}(z_M)| \ge \Lambda. 
\end{equation}
}

The constants 
$K_{\mathrm A}$, 
$K_{\mathrm B}$, and
$\Lam$
above depends only on the choice of 
$\chat$, $\nu$, and the thickness $\delta$ of $\cI(\theta, \delta)$
(equivalently, the angle $A_0$ of $\Delta=\Delta(A_0,T_0)$). 
The proofs of these lemmas will be given later.

By assuming these three lemmas, we can give a proof of the main theorem:

\paragraph{Proof of the main theorem assuming Lemmas A, B, and C.}
It is enough to show the theorem for $c \approx \chat$.
(Indeed, if $c$ stays a uniform distance away
from $\partial \X$, the derivative is bounded above
by the inequality \eqref{eq_CK1_Prop3.1}.)

For a given $c \approx \chat$
and $z_\ast \in J(f_\sigma)$, 
let $z_0=z(c)=H(c, z_\ast) \in J(f_c)$. 
We consider the decomposition
$\N_0=[0, M_1)\sqcup \sS_1\sqcup\sS_2\sqcup\cdots$ 
as in (\ref{eq_decomposition}), (\ref{eq_decomposition2}) or (\ref{eq_decomposition3}).
Then we have
\begin{align*}
\abs{\frac{d}{dc}z(c)}
&=\left|\sum_{n=1}^\infty \frac{1}{Df_c^n(z_0)}\right|
\le 
\sum_{n=1}^{M_{1}}
\frac{1}{|Df_c^{n}(z_0)|}
+
\left| \sum_{k \ge 1}\sum_{n \in \sS_k}\frac{1}{Df_c^{n+1}(z_0)}\right| \\
&\le 
\sum_{n=1}^{M_{1}}
\frac{1}{|Df_c^{n}(z_0)|}
+
\sum_{k \ge 1, \sS_k \neq \emptyset} \frac{1}{|Df_c^{M_k}(z_{0})|}
\left|\sum_{i=1}^{M_{k+1}-M_k}
\frac{1}{Df_c^i(z_{M_k})}\right|.
\end{align*}
By Lemma B, we obviously have $1/|Df_c^{M_1}(z_{0})| \le K_{\mathrm B}$.
By Lemma C, we have 
$$
|Df_c^{M_k}(z_{0})| =
|Df_c^{M_k-M_{k-1}}(z_{M_{k-1}})| 
\cdots
|Df_c^{M_2-M_{1}}(z_{M_{1}})| ~
|Df_c^{M_{1}}(z_{0})| 
\ge \Lam^{k-1}/K_{\mathrm B}
$$
as long as $\sS_k \neq \emptyset$.
Hence by Lemma A, we have 
$$
\left|\sum_{n=1}^\infty \frac{1}{Df_c^n(z_0)}\right|
\le
K_{\mathrm B} + \sum_{k \ge 1} \frac{K_{\mathrm B}}{\Lam^{k-1}} 
\cdot \frac{K_{\mathrm A}}{|c-\chat|^{1-1/Q}}
=
K_{\mathrm B} + 
\frac{K_{\mathrm B}\Lam}{\Lam-1} 
\cdot \frac{K_{\mathrm A}}{|c-\chat|^{1-1/Q}}.
$$
We may assume that ${|c-\chat|} \le 1$ 
such that $K_{\mathrm B} \le K_{\mathrm B}/|c-\chat|^{1-1/Q}$.
Hence by setting 
$K:= K_{\mathrm B} + \dfrac{K_{\mathrm B}  K_{\mathrm A} \Lam}{\Lam -1}$, 
we have 
$
\abs{\dfrac{dz}{dc}} \le \dfrac{K}{|c-\chat|^{1-1/Q}}
$
for any $c \approx \chat$.
\QED\\

\section{Postcritical set and hyperbolic metric} \label{sec:Hyperbolic metrics}

In this section we show some properties of 
the hyperbolic metric on the complement of the postcritical set
of $f_\chat$. 

\paragraph{Postcritical sets.}
The {\it postcritical set} $P(f_c) \subset \C$ of the polynomial $f_c(z)=z^2+c$ is
defined by
$$
P(f_c):=
\overline{
\brac{f_c(0), \,f_c^2(0),\, f_c^3(0),\, \cdots}
}.  
$$
In this paper we only consider $P(f_\chat)$,
which is a countable set and accumulates
only on the parabolic cycle of $f_\chat$.
In particular, the universal covering of $\C-P(f_\chat)$ 
is the unit disk $\D$.

Let $\dist(z, P(f_\chat))$ be the Euclidean distance between 
$z$ and $P(f_\chat)$.
The following lemma provides an estimate
of $\dist(z, P(f_\chat))$ for 
$z \in J(f_c) - \cV(c)$ for $c \approx \chat$,
and plays a very important role. 

\paragraph{Lemma D.}
{\it
There exists a constant $K_{\mathrm D}  \in (0,1]$ 
such that for any sufficiently small $\nu \in (0,R_0)$,
we have  
$\dist(z, P(f_\chat)) \ge K_{\mathrm D} \nu$ 
for any $z \in J(f_c) - \cV(c)$ with $c=c_t \approx \chat~(t \in \Delta)$
by taking a sufficiently small $T_0$ (which depends on $\nu$) in the definition of 
$\Delta=\Delta(A_0,T_0)$.
}

\medskip
\noindent
The proof will be given later in Section \ref{sec:Proof_of_Lemma_D}.

\paragraph{Hyperbolic metric.}
For $c=\chat$, let $\rho(z)|dz|$ denote 
the hyperbolic metric of $\C-P(f_\chat)$,
which is induced by the metric $|dz|/(1-|z|^2)$ 
of constant curvature $-4$ on the unit disk. 
The metric $\rho(z)|dz|$ has the following properties: 
\begin{enumerate}[(i)]
\item
$\rho: \C-P(f_\chat) \to \R_+$ is real analytic
and diverges on $P(f_\chat) \cup \{\infty\}$.
\item
if both $z$ and $f_\chat(z)$ are in $\C-P(f_\chat)$, 
we have 
$$
\frac{\rho(f_\chat(z))}{\rho(z)}|Df_\chat(z)| >1.
$$
\end{enumerate}
See \cite[\S 3.2]{Mc1} for example.

There is an important observation: 
Fix any compact subset $\Gamma$ of $\C-P(f_\chat)$,
and suppose that both $z$ and $f_c(z)$ are 
contained in $\Gamma$ for some $c \approx \chat$.
Then by continuity of 
$\rho(f_c(z))|Df_c(z)|$ with respect to $c$,
we have 
$$
\frac{\rho(f_c(z))}{\rho(z)}|Df_c(z)| \ge A, 
\quad\mbox{or equvalently~}
|Df_c(z)| \ge \frac{\rho(z)}{\rho(f_c(z))} \cdot A,
$$
for some uniform constant $A=A_\Gamma>1$ independent of $c \approx \chat$.

To give some estimates of the function $\rho(z)$, 
we need the following:

\begin{prop} \label{previous_Lemma_W}
 
The hyperbolic metric $\rho(z)|dz|$ of $\C-P(f_\chat)$ satisfies
$$
\rho(z) \le \frac{1}{\dist(z,P(f_\chat))}.
$$ 
\end{prop}

\medskip
This is just an application of a standard fact 
on hyperbolic metrics. 
See \cite[Theorems 1.10 \& 1.11]{Ah} for example.
\paragraph{}
Now we are ready to show:

\paragraph{Lemma E.}
{\it If the constant $\nu$ is sufficiently small,
there exists a constant 
$K_{\mathrm E} \asymp \nu$ with the following property:
For any $c \approx \chat$, we have
$$
\frac{\rho(z)}{\rho(\zeta)} \ge K_{\mathrm E}
$$
if $z , \zeta \in J(f_c)-\cV(c)$.
}

\paragraph{Proof.}
Note that $J(f_c) \subset \overline{\D(2)}$ for $c \in \M$.
Since $\rho$ diverges only at the postcritical set $P(f_\chat)$ 
in $\overline{\D(3)}$,
there exists a constant $C_1>0$ such that $\rho(w) \ge C_1$
for any $w \in \overline{\D(3)} - P(f_{\chat})$.
In particular, we have $\rho(z) \ge C_1$. 
By Lemma D, for $\nu \ll 1$ and $c \approx \chat$ 
we have 
$$
\dist (\zeta, P(f_\chat)) \ge K_{\mathrm D}  \nu
$$
and thus Propositon \ref{previous_Lemma_W} implies $\rho(\zeta) \le 1/(K_{\mathrm D} \nu)$
for sufficiently small $\nu$. 
Now we have $\rho(z)/\rho(\zeta) \ge C_1 K_{\mathrm D} \nu =:K_{\mathrm E}$.\QED

\medskip

By Lemma E we obtain a kind of uniform expansion of $f_c$
with respect to $\rho$:

\paragraph{Lemma F.}
{\it There exists a constant $A>1$ such that 
for $c \approx \chat$, if $z, f_c(z), \ldots ,f_c^n(z)$ 
are all contained in $J(f_c)-\cV(c)$, we have 
$$
|Df_c^n(z)| \ge K_{\mathrm E}A^n.
$$
}
\paragraph{Proof.}
By Lemma D, $J(f_c)-\cV(c)$ is contained in a compact set 
$\Gamma$ in $\C-P(f_\chat)$ independent of $c \approx \chat$.
As mentioned in the observation above, 
there exists a constant $A=A_\Gamma>1$ such that
for any $c \approx \chat$,
$$
\frac{\rho(f_c(w))}{\rho(w)}|Df_c(w)| \ge A
$$
if both $w, \, f_c(w) \in J(f_c)-\cV(c)$.

By the chain rule, we have
\begin{equation}
|Df_c^n(z)| 
= 
\prod_{i=0}^{n-1}|Df_c(f_c^i(z))|
\ge  
\prod_{i=0}^{n-1}\frac{\rho(f_c^{i}(z))}{\rho(f_c^{i+1}(z))}A
\ge
\frac{\rho(z)}{\rho(f_c^n(z))}A^n. \label{ineqInF}
\end{equation}
By applying Lemma E with $\zeta=f_c^n(z)$,
 we obtain the desired inequality. \QED

\section{Proof of Lemma A assuming Lemmas G and H} \label{sec: proof_of_A}
For any $c =c_t \approx \chat$, 
we choose an arbitrary $z_0 \in J(f_c)$
and let $z_n:=f_c^n(z_0)~(n \in \N_0)$. 
For a given S-cycle $\sS=[M,M')$ of this orbit,
we may assume that $M=0$ without loss of generality.
We divide the proof into two cases.

First we suppose that $\sS$ is either a finite S-cycle 
or an infinite S-cycle of type (I). 
Then there exist $m \in \N$ 
and $L \in \N \cup \brac{\infty}$ such that 
\begin{itemize}
\item
$z=z_0 \in V_0$; 
\item
$z_{(m-1)pq} \in U_0$ but $z_{mpq} \notin U_0$;
\item
$z_{mpq+i} \notin V_0$ if $0 \le i <L$; and
\item
$M'<\infty$ iff $L<\infty$ and $M'=mpq+L$.
\end{itemize}
As in Remark \ref{rem_L}, 
by shrinking $U_0$ (that is, by taking a smaller $R_0$) if necessary,
each $z_{mpq+i}$ with $0 \le i \le p$ remains near $f^i(U_0)$ and 
we may assume that $L > p$ for any S-cycle.

Recall that under the radial access condition $c=c_t \approx \chat$, 
we have a periodic point $b_t$ of period $p$ in $U_0$
such that $b_t \to \bhat$ as $t \in \Delta$ tends to $0$ 
(Proposition \ref{prop_near_parabolic}). 
For $c=c_t \approx \chat$,
let $b_0(c), b_1(c), \ldots, b_{p-1}(c)$ 
denote the periodic points 
$b_t, f_{c_t}(b_t), \cdots, f_{c_t}^{p-1}(b_t)$
respectively.

The next two lemmas rely on local dynamics of perturbed parabolic cycle.

\paragraph{Lemma G.}
{\it
Suppose that $q \ge 2$.
There exists a constant $K_{\mathrm{G}}>0$ with the following property:
For any $c \approx \chat$ and $z_0 \in J(f_c) \cap U_0$
such that $z_p$, $z_{2p}, \ldots, z_{qp} \in U_0$, 
we have 
$$
\left|\sum_{l=0}^{q-1}\frac{1}{Df_c^{lp}(z_j)}\right|
\le  K_{\mathrm{G}} 
\big(|c-\chat|+|z_j-b_j(c)|\big).
$$
for each $0 \le j \le p-1$.
 }
\medskip

\paragraph{Lemma H.}
{\it 
There exists a constant $K_{\mathrm{H}}>0$ with the following property:
For any $c \approx \chat$, 
$z_0 \in J(f_c) \cap V_0$,
and $m \in \N$ 
such that $z_{kpq} \in U_0$ for $0 \le k \le m-1$,
we have for each $0\le j\le p-1$ that
$$
\sum_{k=0}^{m-1} 
\frac{1}{|Df_c^{kpq}(z_j)|}
\le 
\frac{  K_{\mathrm{H}} }{|\,c-\chat\,|^{1/2}}
\qquad  \mbox{if }q=1,
$$
and
$$
\sum_{k=0}^{m-1} 
\frac{|c-\chat|+| f_c^{kpq}(z_j) -b_j(c)|}{ |Df_c^{kpq}(z_j)|}
\le    
\frac{K_{\mathrm{H}}}{|\,c-\chat\,|^{1-1/q}}  
   \qquad \mbox{if }q\ge 2.
$$
}
\medskip

The proofs of Lemmas G and H will be given later. 

\paragraph{Finite S-cycles.}
Set $f=f_c$. For the finite S-cycle 
$\sS=[0,M')$ with $M'=mpq+L$, we have 
\begin{align}
&\sum_{i=1}^{M'} \frac{1}{Df^i(z_0)} \nonumber \\
=&
\sum_{j=1}^{p} \frac{1}{Df^{j}(z_{0})}
\sum_{k=0}^{m-1} \frac{1}{Df^{kpq}(z_j)}
\sum_{l=0}^{q-1} 
\frac{1}{Df^{lp}(z_{kpq+j})}  \label{5P1} \\
&~~~~
+ 
\sum_{i=1}^{L-p} 
\frac{1}{Df^{mpq}(z_0)\, Df^{i}(z_{mpq})} 
+
\sum_{i=1}^{p} 
\frac{1}{Df^{mpq+(L-p)}(z_0)\, Df^{i}(z_{mpq+(L-p)})}. \nonumber 
\end{align}
When $1 \le j \le p$, we have $z_j \in f^j(U_0)$ 
and thus 
$$
|Df^j(z_{0})|  \ge  \xi^{j}
$$
in \eqref{5P1} by Remark \ref{rem_xi}.
Note that by the Koebe distortion theorem,
$$
  |Df^{mpq}(z_j)|\ge C_{2} \cdot \frac{\diam U_0}{\diam V_0} \ge \frac{C_{3}}{\nu}
$$
for some constants $C_{2}$ and $C_{3}$ independent of $c$ and $j\in \{0, 1, \ldots, p-1\}$.
Consequently,
\begin{eqnarray*}
\lefteqn{
\sum_{j=1}^{p} \frac{1}{|Df^{j}(z_{0})| }
\sum_{k=0}^{m-1} \frac{1}{|Df^{kpq}(z_j)|}
\left|\sum_{l=0}^{q-1} \frac{1}{Df^{lp}(z_{kpq+j})}\right|
} \\
&\le &\begin{cases}
       \displaystyle  \sum_{j=1}^{p} \frac{1}{\xi^{j}} \sum_{k=0}^{m-1} \frac{1}{|Df^{kpq}(z_j)|} 
       & \mbox{if } q=1 \\
\displaystyle  
\sum_{j=1}^{p} \frac{1}{\xi^{j}} \sum_{k=0}^{m-1} 
\frac{K_{\mathrm G} \, \paren{|c-\chat|+\left|z_{kpq+j}-b_j(c)\right|}}{|Df^{kpq}(z_j)|}  
& \mbox{if $q\ge 2$ (by using Lemma G)}
          \end{cases} \\
&\le& \frac{p}{\xi^p}\cdot \max\{1, K_{\mathrm G}\} 
\cdot \frac{K_{\mathrm H}}{|c-\chat|^{1-1/Q}}\qquad\mbox{(by Lemma H)}
\end{eqnarray*}
where $Q=\max\{2,\,q\}$.

When $n=mpq+i$ with $1 \le i \le L-p$, 
we have $z_{mpq},\,z_{mpq+i} \notin \cV(c)$ and thus
\begin{align*}
|Df^n(z_0)| &= 
|Df^{mpq}(z_0)|\,|Df^{i}(z_{mpq})|\\
&\ge |Df^{mpq}(z_0)| \cdot K_{\mathrm E} \cdot A^{i}. 
\end{align*}
Here the constant $A$ above is the same as that of Lemma F. 
Consequently,  
\begin{eqnarray*}
\sum_{i=1}^{L-p}  \frac{1}
                                      {     
                                        |Df^{mpq}(z_0)|\, |Df^{i} (z_{mpq})|
                                       }
&  \le & \sum_{i=1}^{L-p} \frac{1}
                                                 {    |Df^{mpq}(z_0)|
                                                  ~K_{\mathrm E}~A^i
                                                 } \\
&\le& \frac{\nu}{C_{3}~ K_{\mathrm E}}\frac{1}{A-1}.
\end{eqnarray*}

When $n=mpq+(L-p)+i$ with $1 \le i \le p$, 
$f(z_n) \in \cV(c)$ and thus
\begin{align*}
|Df^n(z)| &= 
|Df^{mpq+(L-p)}(z_0)|\,|Df^i(z_{mpq+(L-p)})|\\
&\ge |Df^{mpq}(z_0)| \cdot K_{\mathrm E} \cdot A^{L-p} \cdot \xi^{i}\\
&\ge|Df^{mpq}(z_0)| \cdot K_{\mathrm E} \cdot 1 \cdot \xi^{p}.
\end{align*}
Consequently,  
\begin{eqnarray*}
\sum_{i=1}^{p}  \frac{1}{|Df^{mpq+(L-p)}(z_0)|\, |Df^{i}(z_{mpq+(L-p)})|}
&\le & \sum_{i=1}^{p} \frac{1}{|Df^{mpq}(z_0)|~K_{\mathrm E}~\xi^p} \\
&\le& \frac{p \nu }{C_{3}~ K_{\mathrm E}~\xi^p}.
\end{eqnarray*}

By these estimates, when $M'<\infty$, we have:
\begin{eqnarray*}
\left|\sum_{i=1}^{M'} \frac{1}{Df^i(z_0)}\right| 
&\le&
\sum_{j=1}^{p} \frac{1}{|Df^{j}(z_{0})| }
 \sum_{k=0}^{m-1} \frac{1}{|Df^{kpq}(z_j)|}
\left|\sum_{l=0}^{q-1} \frac{1}{Df^{lp}(z_{kpq+j})}\right|
\\
&&\quad
+ 
\sum_{i=1}^{L-p} 
\frac{1}{|Df^{mpq}(z_0)|\, |Df^{i}(z_{mpq})|} \\
&&\quad
+
\sum_{i=1}^{p} 
\frac{1}{|Df^{mpq+(L-p)}(z_0)|\, |Df^{i}(z_{mpq+(L-p)})|} \\
&\le & 
\frac{p}{\xi^{p}}\cdot\max\{1,K_{\mathrm G}\}
\cdot\frac{K_{\mathrm H}}{|c-\chat|^{1-1/Q}}
+
\frac{\nu}{C_{3}~K_{\mathrm E}}
\paren{
\frac{1}{A-1}
+
\frac{p}{\xi^{p}}}.
\end{eqnarray*}
Hence by setting
$$
K_{\mathrm A}:=
\max\brac{1, K_{\mathrm G}}
 \cdot \frac{2p}{\xi^p}\cdot K_{\mathrm H},
$$
we obtain our desired estimate
$$
\left|\sum_{i=1}^{M'} \frac{1}{Df^i(z_0)}\right| 
\le \frac{K_{\mathrm A}}{ |\,c-\chat\,|^{1-1/Q}}
$$
when $|c-\chat|$ is sufficiently small.
Note that $K_{\mathrm A}$ does not depend on $m$ and $L$.

\paragraph{Infinite S-cycles, Type (I).}
If $M'=\infty$, then $L=\infty$ and one can easily check 
\begin{align*}
\left|\sum_{i=1}^{\infty} \frac{1}{Df^i(z_0)}\right|
\le&
\frac{p}{\xi^{p}}\cdot\max\{1,K_{\mathrm G}\}\cdot\frac{K_{\mathrm H}}{|c-\chat|^{1-1/Q}}
+
\frac{\nu}{C_{3}~K_{\mathrm E}}
\paren{
\frac{1}{A-1}
}\\
<& \frac{K_{\mathrm A}}{ |\,c-\chat\,|^{1-1/Q}}
\end{align*}
by the same argument as the case of finite S-cycles.

\paragraph{Infinite S-cycles, Type (II).}
Next we suppose that $\sS=[0,\infty)$ 
is an infinite S-cycle of Type (II).
As we have seen in Section 2
 (see in particular Propositions \ref{prop_near_parabolic} and \ref{prop_leaving_U_0}),
$z_0$ must be a repelling periodic point of period $p$ or $pq$
of $f_c$ contained in $V_0$ for $c =c_t \approx \chat~(t \in \Delta)$ 
by taking a smaller $T_0>0$ in the definition of $\Delta$ if necessary.
Note that Lemmas G and H are valid even in such a case.

\paragraph{Case 1 ($q=1$).}
Since the assumption of Lemma H is valid for any large $m$,
we have
$$
\sum_{k=0}^{\infty} 
\frac{1}{|Df_c^{kp}(z_j)|}
\le 
\frac{  K_{\mathrm{H}} }{|\,c-\chat\,|^{1/2}}
$$
for $q=1$ and $0 \le j \le p-1$. 
Hence for $c \approx \chat$,
we have 
\begin{align*}
\abs{\sum_{i =1}^\infty\frac{1}{Df^i(z_0)}}
&\le 
\sum_{k =0}^\infty \sum_{j=1}^{p}
\frac{1}{|Df^{kp}(z_j)||Df^{j}(z_0)|}\\
&\le 
\sum_{k =0}^\infty \sum_{j=1}^{p}
\frac{1}{|Df^{kp}(z_j)| \cdot \xi^j}
\le 
\frac{ p K_{\mathrm{H}} }
{ \xi^p |\,c-\chat\,|^{1/2}}.
\end{align*}

\paragraph{Case 2$^\pm$.}
As in Case 1, Lemma H implies 
$$
\sum_{k=0}^{\infty} 
\frac{|c-\chat|+|f_c^{kpq}(z_j) -b_j(c)|}{|Df_c^{kpq}(z_j)|}
\le 
\frac{  K_{\mathrm{H}} }{|\,c-\chat\,|^{1-1/q}}
$$
for $q \ge 2$ and $0 \le j \le p-1$. 
(In Case 2$^+$, $b_j(c)$ is a repelling periodic point 
and we have $f_c^{kpq}(z_j)=b_j(c)$ for any $k \ge 0$.)
Hence for $c \approx \chat$, we have
\begin{align*}
\abs{\sum_{i =1}^\infty\frac{1}{Df^i(z_0)}}
&\le
\sum_{j=1}^{p} \frac{1}{|Df^{j}(z_{0})| }
\sum_{k=0}^{\infty} \frac{1}{|Df^{kpq}(z_j)|}
\left|\sum_{l=0}^{q-1} \frac{1}{Df^{lp}(z_{kpq+j})}\right|
 \\
&\le 
\sum_{j=1}^{p} \frac{1}{\xi^{j}} 
\sum_{k=0}^{\infty} 
\frac{1}{|Df^{kpq}(z_j)|} 
\cdot
K_{\mathrm G}\big(|c-\chat|+\left|z_{kpq+j}-b_j(c)\right|\big) \\
&\le 
\frac{p}{\xi^p}\cdot K_{\mathrm G} \cdot \frac{K_{\mathrm H}}{|c-\chat|^{1-1/q}}
\end{align*}
by Lemma G.
In both Case 1 and Case 2$^\pm$,
we conclude that  
Lemma A is valid with the same constant $K_{\rm A} > 0$
defined as above. \QED

\section{Proof of Lemma B} \label{sec:Proof_B}
Without loss of generality, we may assume that either
\begin{itemize}
\item
$M < \infty$ and $z_{M} \in V_0$; or
\item
$M=\infty$.
\end{itemize}
For the first case, if $M \le p$, then 
$z_0 \in \bigcup_{k=0}^p f_c^{-k}(V_0)$
since $z_M \in V_0$. Hence 
$$
|Df_c^i(z_0)| \ge \xi^i \ge \xi^p 
$$ 
for $i \in [1,M]$. This implies 
$$
\sum_{i=0}^M\frac{1}{|Df_c^i(z_0)|} 
\le \frac{M}{\xi^p}
\le \frac{p}{\xi^p}.
$$
If $M > p$, then $z_{M-p} \in f_c^{-p}(V_0)$ and 
$z_i \notin \bigcup_{k=0}^p f_c^{-k}(V_0)$ 
for $i \in [0, M-p)$.
Then 
\begin{itemize}
\item
$|Df_c^i(z_0)| \ge K_{\mathrm E} A^i$ for $i \in [1, M-p]$ (by Lemma F)
\item
$|Df_c^{M-p+j}(z_0)| 
\ge K_{\mathrm E} A^{M-p} \cdot \xi^j 
\ge  K_{\mathrm E} \xi^p  $
for $j \in [1,p]$.
\end{itemize}
This implies:
$$
\sum_{i=1}^M \dfrac{1}{|Df_c^i(z_0)|}
\le 
\sum_{i=1}^{M-p}\dfrac{1}{K_{\mathrm E} A^i}
+
\sum_{j=1}^{p}\dfrac{1}{K_{\mathrm E} \xi^p}
\le
\frac{1}{K_{\mathrm E}}
\paren{
\frac{1}{A-1}
+
\frac{p}{\xi^p}
}.
$$
By defining
$$
K_{\mathrm B}:=\frac{1}{K_{\mathrm E}}
\paren{
\frac{1}{A-1}
+
\frac{p}{\xi^p}
}
$$
we have the desired estimate.

For the second case ($M=\infty$), the orbit never land on $\cV(c)$ 
and by Lemma F, we have
$$
|Df_c^i(z_0)| 
\ge K_{\mathrm E} A^i
$$
for any $i \in \N$. Hence 
$$
\sum_{i=1}^\infty \dfrac{1}{|Df_c^i(z_0)|}
\le 
\sum_{i=1}^{\infty}\dfrac{1}{K_{\mathrm E} A^i}
\le
\frac{1}{K_{\mathrm E}}
\cdot \frac{1}{A-1} < K_{\mathrm{B}}.
$$
\QED

\section{Proofs of
Propositions \ref{prop_c_and_epsilon}~and \ref{prop_near_parabolic}
}
\label{sec_parabolic}

This section is devoted to the proofs of
Propositions \ref{prop_c_and_epsilon} and
\ref{prop_near_parabolic},
which describe the local properties of parabolic periodic points and its perturbations 
along thick internal rays. 
The argument here relies on some well-known facts 
originally due to Douady and Hubbard \cite[Expos\'e XIV]{DH}
on perturbation of parabolic cycles. 
Here we will adopt Milnor's formulation in \cite[\S 4]{Mi1}.
(See also \cite[Thm. 1.1 and Thm. A.1]{T}.)

\paragraph{}
Let $\chat$ be a parabolic parameter 
and $\bhat$ be a parabolic periodic point of 
period $p$ of $f_\chat$. Let $\hat{\lambda}$ be the 
multiplier $Df_\chat^p (\bhat)$ of $\bhat$
that is a primitive $q$-th root of unity.

\begin{prop}[Lemma 4.5 of {\cite{Mi1}}]
\label{lem_Milnor}
There exist unique single valued functions $c(\lam)$
and $z(\lam)$ defined on a neighborhood of $\hat{\lam}$
so that $z(\lam)$ is a periodic point of period $p$
and multiplier $\lam$ for the map $f_{c(\lam)}$,
with $\chat=c(\hat{\lam})$ and $\bhat=z(\hat{\lam})$.
This function $c(\lam)$ has a single critical point at $\hat{\lam}$
when $q=1$ (Case 1) 
but is univalent when $q \ge 2$ (Case 2$^\pm$).  
\end{prop}

By applying the argument of \cite[Appendix A.4]{K3} 
one can find a holomorphic family of local coordinates 
$\zeta=\varphi_\lam(z)$ and $w=\psi_\lam(z)$
defined near $z(\lam)$ such that
$\varphi_\lam(z(\lam)) =\psi_\lam(z(\lam))=0$, 
\begin{align}
\varphi_{\lam} \circ f_{c(\lam)}^{p} \circ \varphi_{\lam}^{-1}(\zeta)
&= \lam \zeta +\zeta^{q+1} +O(\zeta^{2q+1}),
 ~\text{and} \label{eq_Milnor_varphi}
\\[.5em]
\psi_{\lam} \circ f_{c(\lam)}^{pq} \circ \psi_{\lam}^{-1}(w)
&= \lam^q w \,(1+w^q +O(w^{2q}))
\label{eq_Milnor_psi}
\end{align}
where 
$$
w= 
\paren{
\frac{1+\lambda^q+\lambda^{2q}+\cdots+\lambda^{(q-1)q}}
{\lam}}
^{1/q} \zeta. 
$$
(Note that the error terms in
\eqref{eq_Milnor_varphi} and 
\eqref{eq_Milnor_psi} are slightly refined 
compared with the similar coordinates given in 
\cite[Prop. 11.1]{DH} and \cite[Lem.4.2]{Mi1}.)
Then we take a sufficiently small $\hat{R}>0$
such that the domains of these local coordinates contain
$\hat{U}=\D(\bhat,\hat{R})$ for $\lam$ sufficiently close to $\hat{\lam}$. 
In particular, both $D\varphi^{-1}_{\lam}(0)$ and $D\psi^{-1}_{\lam}(0)$ 
are uniformly bounded away from zero.

\begin{prop}
\label{lem_q-roots}
For $\lam \neq \hat{\lam}$, one can find exactly $q$ non-zero fixed points of 
$\psi_{\lam} \circ f_{c(\lam)}^{pq} \circ \psi_{\lam}^{-1}(w)
= \lam^q w \,(1+w^q +O(w^{2q}))$
of the form $w=w_\lam(1+ o(1))$ 
with multiplier $1+q(1-\lam^q)+o(|\lam-\hat{\lam}|)$, where $w_\lam$ is a $q$-th root of $(1-\lambda^q)/\lambda^q$.
\end{prop}

\paragraph{Proof.}
Since the equation of the form 
$\lambda^q w \,(1+w^q +O(w^{2q}))=w$
is regarded as a perturbation of the $\lam=\hat{\lam}$ case, 
it has exactly $q$ non-zero roots by Hurwitz's theorem. 
To obtain the estimate of the solution, 
we apply Rouch\'e's theorem\footnote{
Let $F(w)=w^q-w_\lam^q 
= w^q-(1-\lam^q)/\lam^q$ and $G(w)=O(w^{2q})$. 
Let $\ell:=|\lam-\hat{\lam}|$.
Then consider the circle $|w-w_\lam|=\ell^s$ 
with $1/q < s < 1+1/q$.
For example, one can take $s:=(2q+1)/(2q)$.
Now one can check $|w_\lam| \asymp \ell^{1/q}$
and thus $|F(w)|\asymp \ell^{1-1/q+s}$ and $|G(w)|=O(\ell^{2})$
on this circle. Since $|F(w)|>|G(w)|$ for $\ell \ll 1$, 
$F(w)+G(w)=0$ and $F(w)=0$ has the same number of zeroes,
which is one. Hence the solution is of the form 
$w=w_\lam+O(\ell^s)=w_\lam(1+o(1))$.
}.
The multiplier comes from 
$D(\psi_{\lam} \circ f_{c(\lam)}^{pq} \circ \psi_{\lam}^{-1})
(w)= \lam^q (1+(q+1) w^q +O(w^{2q}))$ and $|w_\lam| \asymp|\lam-\hat{\lam}|^{1/q}$.
\QED

\paragraph{}
Now we fix a hyperbolic component $\X$ attached to $\chat$
with uniformization $\Phi=\Phi_\X:\Dbar \to \Xbar$,
and a Stolz angle $S(\Delta) \subset \D$ at $\hat{\lam}$
with $\Delta=\Delta(A_0,T_0)$. 
In what follows we will take a sufficiently 
small $T_0$ if necessary. 
Let us start with Case 2$^-$.

\paragraph{Proofs of Propositions \ref{prop_c_and_epsilon}
and \ref{prop_near_parabolic} for Case 2$^-$.}
We obtain Case 2$^-$ by letting $\lam=\lam_t:=(1-t/q) \hat{\lam}$ 
with $t \in \Delta$ in Proposition \ref{lem_Milnor}. 
Indeed, the hyperbolic component $\X$ is specified by 
the family of attracting fixed points $z(\lam)$ 
of $f_{c(\lam)}^p$ with $\lam=\lam_t$. 
Proposition \ref{lem_Milnor} implies that
$
|c_t -\chat| \asymp |\lam_t -\hat{\lam}| \asymp |t|
$ 
by setting $c_t=c(\lam_t)$. 
This proves Proposition \ref{prop_c_and_epsilon} for this case.

Items (1), (2) and (4) of Proposition \ref{prop_near_parabolic}
are straightforward by applying Proposition \ref{lem_q-roots} 
with
$b_t:=z(\lam_t)$, 
$\varphi_{t}:=\varphi_{\lam_t}$, and
$\psi_{t}:=\psi_{\lam_t}$
for $t \in \Delta$.

\paragraph{Proofs of Propositions \ref{prop_c_and_epsilon}
and \ref{prop_near_parabolic} for Case 1 and Case 2$^+$.}
In these cases, $q$ non-zero fixed points of $f_\lam^{pq}$
given in Proposition \ref{lem_q-roots} should be attracting 
with multiplier $1-qt$ for $t \in \Delta$,
and this specifies the hyperbolic component $\X$ in Case 1 or Case 2$^+$.
Hence by the equality $1-qt=1+q(1-\lam^q)+o(|\lam-\hat{\lam}|)$
we obtain $\lam=\lam_t:=(1+t/q+o(t))\hat{\lam}$ for $t \in \Delta$. 
By setting $c_t:=c(\lam_t)$, Proposition \ref{lem_Milnor} implies that
$
|c_t -\chat| \asymp |\lam_t -\hat{\lam}|^2 \asymp |t|^2
$
in Case 1, and 
$
|c_t -\chat| \asymp |\lam_t -\hat{\lam}| \asymp |t|
$
in Case 2$^+$.
This proves Proposition \ref{prop_c_and_epsilon} for these cases.

Items (1), (2) and (3) of Proposition \ref{prop_near_parabolic}
are straightforward by applying Proposition \ref{lem_q-roots} 
with
$b_t:=z(\lam_t)$, 
$\varphi_{t}:=\varphi_{\lam_t}$, and
$\psi_{t}:=\psi_{\lam_t}$
for $t \in \Delta$.
\QED

\section{Proof of Lemma G} \label{sec:Proof_G}

\paragraph{Preliminary.}
Lemma G comes from 
symmetric behavior of the local dynamics of $f_c^p$
near the periodic point $b(c)$.
Indeed, it holds without the radial access condition
and the assumption that $z_0 \in J(f_c)$.
To describe its principle in a generalized form,
we adopt Milnor's formulation  
as above.

Let $\lam$ range over a neighborhood $\hat{S}$ of $\hat{\lam}$
such that the holomorphic maps $c(\lam)$ and $z(\lam)$
in Proposition \ref{lem_Milnor} are defined.
Note that the map $\lam \mapsto c(\lam)$ gives an isomorphism
between $\hat{S}$ and a neighborhood of $\chat=c(\hat{\lam})$
since we are in the case of $q \ge 2$ (Case 2$^\pm$).
We claim:

\begin{prop}\label{prop_Lemma_O'}
Suppose that $q \ge 2$.
For any $\lambda \in \hat{S}$ and 
$z_0 \in \hat{U}=\D(\hat{b},\hat{R})$,
suppose that 
$z_{lp}:=f_{c(\lam)}^{lp}(z_0) \in \hat{U}$ 
for $0 \le l \le q$.
Then we have
$$
\left|\sum_{l=0}^{q-1}\frac{1}{Df_{c(\lam)}^{lp}(z_0)}\right|
=O(|c(\lam)-\chat|)+O(|z_0-z(\lam)|)
$$
by taking smaller $\hat{S}$ and $\hat{U}$ if necessary.
\end{prop}

\paragraph{Proof of Proposition \ref{prop_Lemma_O'}.}
Let $\e:=\lam/\hat{\lam}-1~ (\lam \in \hat{S})$ 
be an alternative parameter such that
$$
\lam=(1+\e)\hat{\lam} \in \hat{S}. 
$$
Since we are in the case of $q \ge 2$ (Case 2$^\pm$),
we have $\hat{\lam} \neq 1$ and $\lam \mapsto c(\lam)$ is univalent
(Proposition \ref{lem_Milnor}).
Hence we obtain
$|c(\lam)-\chat| \asymp |\lam-\lamhat| \asymp |\e|$ for $\lambda \in \hat{S}$.

Recall that we have a local coordinate $\zeta=\varphi_\lam(z)$
defined on $\hat{U}$
such that $\varphi_\lambda(z(\lam))=0$ and
$$
F_\lam(\zeta)
:=\varphi_\lam \circ f_{c(\lam)}^{p} \circ \varphi_\lam^{-1}(\zeta)
=\lam \zeta +\zeta^{q+1} +O(\zeta^{2q+1})
$$
as in \eqref{eq_Milnor_varphi}.
For a given $z_0 \in \hat{U}$,
let $\varphi_\lam^{-1}(\zeta_0)=z_0$ and 
$\zeta_{lp}:=\varphi_\lam(z_{lp})$.
 
Then it is easy to check that 
\begin{equation}
\zeta_{lp}=F_\lam^{l}(\zeta_0)
=\lambda^l\zeta_0+\lambda^{l-1}\left(1+\lambda^q+\cdots+\lambda^{(l-1)q}\right)\zeta_0^{q+1}+O(\zeta_0^{2q+1}) \label{zetalp}
\end{equation}
for any $1 \le l\le q$.
Since the local coordinate 
$\varphi_\lam(z)$ depends holomorphically on $\lam$,
we have
\begin{equation}
 \varphi_\lam^{-1}(\zeta)
 =z(\lam)+C_{4}\zeta+C_{5}\zeta^2+O(\zeta^3)  
 \label{Fact1implies}
\end{equation}
for some $C_{4}=C_{4}(\lam)$ and $C_{5}=C_{5}(\lam)$ with $C_{4}$ bounded away from zero for $\lam \in \hat{S}$.  

Let
$$
\mathcal{A}(z_0):=
\sum_{l=0}^{q-1} \frac{1}{Df_{c(\lam)}^{lp}(z_0)}.
$$
Now,
$$
Df_{c(\lam)}^{lp}(z_0)=D\varphi_\lam^{-1}(\zeta_{lp})~ DF^l_\lam(\zeta_0)~ D\varphi_\lam(z_0),
$$ 
and
$$
\mathcal{A}(z_0)
=1+\frac{   D\varphi^{-1}_\lam(\zeta_0)}{D\varphi_\lam^{-1}(\zeta_p)}\frac{1}{DF_\lam(\zeta_0)} +\frac{D\varphi_\lam^{-1}(\zeta_0)}{D\varphi_\lam^{-1}(\zeta_{2p})}\frac{1}{DF_\lam^{2}(\zeta_0)}
+\cdots
+\frac{D\varphi^{-1}_\lam(\zeta_0)}{D\varphi_\lam^{-1}(\zeta_{(q-1)p})}\frac{1}{DF_\lam^{q-1}(\zeta_0)}.
$$
By \eqref{zetalp} and \eqref{Fact1implies}, 
\begin{eqnarray*}
 D\varphi_\lam^{-1}(\zeta_{lp}) &=& C_{4}+2C_{5}\zeta_{lp}+O(\zeta^2_{lp}) \\
                                      &=& C_{4}+2C_{5}\lambda^l\zeta_0+O(\zeta_0^2),
\end{eqnarray*}
thus
\begin{eqnarray*}
\frac{D\varphi_\lam^{-1}(\zeta_0)}{D\varphi_\lam^{-1}(\zeta_{lp})}
&=& \frac{C_{4}+2C_{5}\zeta_0+O(\zeta_0^2)}{
C_{4}+2C_{5}\zeta_{lp}+O(\zeta_{lp}^2)
} \\
&=& 1+\frac{2C_{5}}{C_{4}} (1-\lambda^l)~\zeta_0+O(\zeta_0^2).
\end{eqnarray*}
Formula \eqref{zetalp} also gives
$$
 \frac{1}{DF_\lam^{l}(\zeta_0)} = \frac{1}{ 
                     \lambda^l} \left( 1+O(\zeta_0^2)
                                    \right).
$$
Therefore, 
\begin{eqnarray*}
 \mathcal{A}(z_0)  
 &=& 
1+
\sum_{l=1}^{q-1}
\bigg\{\left( 1+ \frac{2C_{5}}{C_{4}}(1-\lambda^l)\zeta_0+O(\zeta_0^2)
\right) \cdot 
\frac{1}{\lambda^l}\left( 1+O(\zeta_0^{2})
 \right)\bigg\} \\
&=&  1+\sum_{l=1}^{q-1}\frac{1}{\lambda^l}\left( 1+ \frac{2C_{5}}{C_{4}}(1-\lambda^l)\zeta_0+O(\zeta_0^2)
 \right)  \\
&=& 
\frac{1}{\lambda^q}
\left( 
\lambda+\lambda^2+\cdots+\lambda^q+\frac{2C_{5}}{C_{4}}
( \lambda+\lambda^2+\cdots+\lambda^{q-1}+\lambda^q-q\lambda^q)\zeta_0
\right)+O(\zeta_0^2)\\
&=&
\frac{\lambda+\lambda^2+\cdots+\lambda^q}{\lambda^q}
\left(1+ 
\frac{2C_{5}}{C_{4}}\zeta_0
\right)
-
\frac{2qC_{5}}{C_{4}} \zeta_0
+O(\zeta_0^2).
\end{eqnarray*}
Clearly, 
$$ 
\hat{\lambda}+\hat{\lambda}^2+\cdots+\hat{\lambda}^{q}=0.
$$
Thus, by using the alternative parameter $\e$, we have
\begin{align*}
\lambda+\lam^2+\cdots+\lambda^{q}
&~=~
\lefteqn{(1+\epsilon)\hat{\lambda}+(1+\epsilon)^2\hat{\lambda}^2+\cdots+
(1+\epsilon)^{q}\hat{\lambda}^{q}} \\
&~=~ 
(\hat{\lambda}+2\hat{\lambda}^2+\cdots+q\hat{\lambda}^{q})~\epsilon+O(\epsilon^2)\\
&~=~ \frac{q\hat{\lam}}{\hat{\lam}-1}\epsilon+O(\epsilon^2).
\end{align*}
Hence, 
\begin{eqnarray*}
 \mathcal{A}(z_0) 
&=& 
\frac{1}{(1+\epsilon)^{q}}
\brac{
\frac{q\hat{\lam}}{\hat{\lam}-1}\epsilon+O(\epsilon^2)
}
\left(1+ 
\frac{2C_{5}}{C_{4}}\zeta_0
\right)
-\frac{2q\,C_{5}}{C_{4}} \zeta_0
\paren{1+O(\zeta_0)} \\
&=& O(\e)+O(\zeta_0)
\end{eqnarray*}
where the implicit constants are 
independent of $\lam =(1+\e) \hat{\lam} \in \hat{S}$
and $\zeta_0 \in \varphi_\lam(\hat{U})$.

Since $C_{4}$ is uniformly bounded away from zero for $\lam \in \hat{S}$,
\eqref{Fact1implies} implies
$
|\zeta_0| \asymp |z_0-z(\lam)|.
$
Since $|\e| \asymp |c(\lam)-\chat|$, 
we conclude that 
$$
| \mathcal{A}(z_0) | = O(|c(\lam)-\chat|)+O(|z_0-z(\lam)|).
$$ 
\QED

\paragraph{Proof of Lemma G.}
When $j=0$, we apply Proposition \ref{prop_Lemma_O'}
by letting 
$\lam:=\lam_t$,
$c=c_t:=c(\lam_t)$, $b_0(c)=b_t:=z(\lam_t)$ 
as in the proof of statements (3) and (4) of 
Proposition \ref{prop_near_parabolic} given in the previous section.
(Note that we have $U_0 \subset \hat{U}$ and $S(\Delta) \subset \hat{S}$ by taking a smaller $T_0$.)
Then we immediately obtain 
$$
\left|\sum_{l=0}^{q-1}\frac{1}{Df_c^{lp}(z_0)}\right|
=O(|c-\chat|)+O(|z_0 -b_0(c)|).
$$

When $1 \le j \le p-1$, 
we can apply the same argument as Proposition \ref{prop_Lemma_O'}
by replacing $\hat{U}$ with $\hat{U}_j:=f_c^j(\hat{U})$ 
and $b_0(c)$ with $b_j(c)$.
More precisely,  
having the local coordinate $\varphi_0=\varphi_{\lam_t}$ 
on $\hat{U}$ for $c=c_t=c(\lam_t)$, 
we can define a local coordinate 
$\varphi_{j}:=\varphi_c\circ f_c^{-j}|_{\hat{U}_j}$ on $\hat{U}_j$ 
for each $1 \le j\le p-1$ 
by a branch such that $\varphi_{j}(z_{lp+j})=\zeta_{lp}$ 
for $0 \le l \le q$.
Since $\varphi_{j}^{-1}$ is of the form
$$
 \varphi^{-1}_{j}(\zeta)=b_j(c)+C_{4,j}\zeta+C_{5,j}\zeta^2+O(\zeta^3)  
$$
where $C_{4,j}$ is a constant bounded away from zero,
we have $|z_j -b_j(c)| \asymp |\zeta_0|$ and
the same argument as Proposition \ref{prop_Lemma_O'} yields
$$
\left|\sum_{l=0}^{q-1}\frac{1}{Df_c^{lp}(z_j)}\right|
=O(|c-\chat|)+O(|z_j -b_j(c)|).
$$
Hence there exists a constant $K_{\mathrm G}>0$
independent of $c=c_t \approx \chat$ and $z_0 \in U_0$
such that 
$$
\left|\sum_{l=0}^{q-1}\frac{1}{Df_c^{lp}(z_j)}\right|
\le K_{\mathrm G}\big(|c-\chat|+|z_j -b_j(c)|\big)
$$
for $0 \le j \le p-1$.
\QED

\section{Proof of Proposition \ref{prop_leaving_U_0}} \label{sec:Proof_of_Leaving}
This section is devoted to the proof of 
Proposition \ref{prop_leaving_U_0}.

\paragraph{Branched coordinates near infinity.}
Let $c=c_t=c(\lam_t) \in \X\cup\{\chat\}$ with 
$\lam_t \in S(\Delta) \cup \{\hat{\lam}\}$ and $t \in \Delta \cup \{0\}$.
It is convenient to use a branched coordinate 
$\Psi_t: \psi_{t}(\hat{U}) \to \Cbar_W$ given by
$W=\Psi_t(w):=-\lam_t^{q^2}/(q\, w^q)$
where $\Cbar_W$ denotes the Riemann sphere in $W$-coordinate.
Let $W=\Phi_t(z):=\Psi_t \circ  \psi_{t}(z)$.
We may assume that 
$\Phi_t(\hat{U})$ is always contained in 
$\brac{W \in \Cbar_W \st |W| \ge \hat{r}}$ 
for some constant $\hat{r}>0$ 
independent of $t \in \Delta$,
with $\hat{r}\hat{R}^q \asymp 1$
by taking a smaller $\hat{R}$ if necessary.

In this coordinate we observe the map $f_{c_t}^{pq}$ 
as $G:=
\Phi_{t} \circ f_{c_t}^{pq} \circ \Phi_{t}^{-1}$
(taking an appropriate branch of $\Phi_{t}^{-1}$)
of the form
$$
G(W)=G_t(W)= \tau \,W +1+O(1/W),
$$
where 
$$
\tau:= \lam_t^{-q^2}
=\left\{\begin{array}{ll} 
1 - qt + O(t^2) & \text{for Case 1 and Case 2$^+$},     
\\[.5em] 
1 + qt + O(t^2) & \text{for Case 2$^-$}  
\end{array}\right.
$$ 
by Proposition \ref{prop_near_parabolic}.
For $t \in \Delta$, let
$$
B:=\frac{1}{1-\tau} = \pm \frac{1}{qt}+O(1).
$$
Note that $B$ is a fixed point of $W \mapsto \tau W +1$ and
$G(W)-B=\tau(W-B)+O(1/W)$.
By Rouch\'e's theorem, 
there is a fixed point 
$B^\ast$ of $G$ of the form $B^\ast=B+O(1)=O(1/t)$ with multiplier $\tau^\ast:=DG(B^\ast)=\tau+O(t^2)$. 

\begin{figure}[htbp]
\begin{center}
\includegraphics
{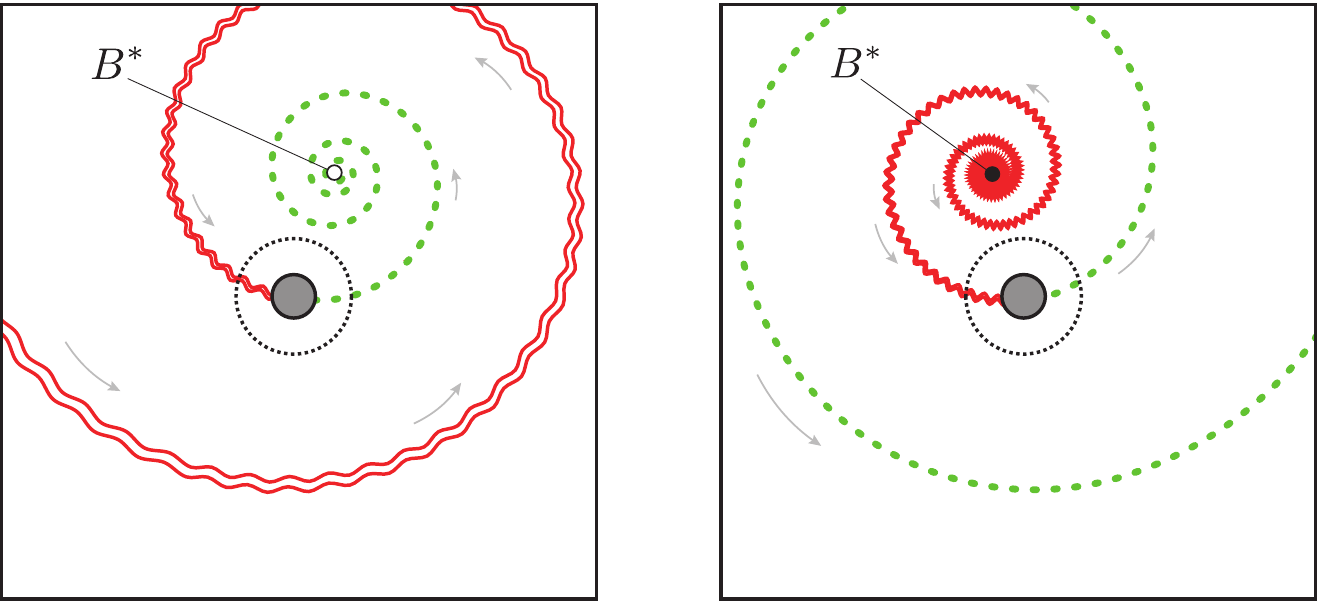}
\caption{The image of the Julia set (red) by $\Phi_t$ 
in $W$-coordinate for Cases 1 and 2$^+$ (left) and 
Case 2$^-$ (right). 
The origin of the $W$-plane is the center of the dotted (or black) circle. The radius of the black circle equals $\hat{r}$. 
The annulus between the dotted circle and the black one encloses 
the region where the dynamics by $G$ is relatively close to the translation $W \mapsto W +1$
(see Proposition \ref{prop_G_near_origin}).} \label{fig_G}
\end{center}
\end{figure}

Figure \ref{fig_G} illustrates the Julia set in $\Cbar_W$.  

\paragraph{Dynamics of $G$ near infinity.}
Let 
$$
d_0:=\frac{\cos A_0}{4} > 0,
$$ 
where $A_0 \in [0, \pi/2)$ is given in the definition of $\Delta$
(hence it is determined by the given thickness of the 
thick internal ray).
 
\begin{prop}
\label{prop_G_near_infinity}
When $|W|>1/(q d_0 |t|)$, we have 
\begin{align*}
|G(W)|
&<|W|  \quad \text{in Cases 1 and 2$^+$, and}     
\\ 
|G(W)|&>|W|  \quad \text{in Case 2$^-$}
\end{align*}
by taking sufficiently small $\hat{R}$ and $T_0$
in the definitions of $\hat{U}$ and $\Delta$.
\end{prop}

\paragraph{Proof.}
Since $G(W)/W=\tau +W^{-1} + O(W^{-2})
= 1 \mp qt +O(t^2) +W^{-1} + O(W^{-2})$,
we have
$$
\log \abs{\frac{G(W)}{W}}
=\Re \log \frac{G(W)}{W} 
=\Re \paren{\mp qt +O(t^2) +W^{-1} + O(W^{-2})}.
$$
Now suppose that $|W|>1/(q d_0 |t|)$.
By taking a smaller $T_0$ if necessary,
for any $t \in \Delta=\Delta(A_0, T_0)$ we have
\begin{align*}
&\Re (q t) \ge q |t| \cos A_0 = 4qd_0 |t|,\\
&\Re (W^{-1}) \le |W^{-1}| \le q d_0 |t|, \quad\text{and}
\quad
\Re \left(O(t^2) + O(W^{-2})\right) \le q d_0 |t|.
\end{align*}
Hence 
$$
\log \abs{\frac{G(W)}{W}} \le (-4+1+1)qd_0 |t| <0
$$
in Cases 1 and 2$^+$, and 
$$
\log \abs{\frac{G(W)}{W}} \ge (4-1-1)qd_0 |t| >0
$$
in Case 2$^-$.
\QED

\paragraph{}
We first prove Proposition \ref{prop_leaving_U_0} for Case 2$^-$,
then prove it for Case 1 and Case 2$^+$.

\paragraph{Proof of Proposition \ref{prop_leaving_U_0}: Case 2$^-$.}
In this case $\tau=1+qt+O(t^2)$ and thus $B=-1/(qt)+O(1)$.
Note that we may assume $|B| > 7 \hat{r}$ for any $t \in \Delta$
by taking a smaller $T_0$ if necessary.
We first claim:

\begin{prop}\label{prop_2-_1}
In Case 2$^-$, for any $W$ with $|W| \ge \hat{r}$ and $|W-B| \ge |B|/2$,
we have 
$$
|G(W)-B| \ge |W-B| +d_0
$$
by taking sufficiently small $\hat{R}$ and $T_0$
in the definitions of $\hat{U}$ and $\Delta$.
\end{prop}

\paragraph{Proof.}
Since $G(W)-B=\tau(W-B) +O(W^{-1})$ 
and $|W-B| \ge |B|/2$,
$$
|G(W)-B|=|\tau| \, |W-B| +|O(W^{-1})|
\ge 
|W-B| + \frac{|\tau|-1}{2|\tau-1|} +|O(W^{-1})|.
$$
Taking smaller $\hat{R}$ (hence a larger $\hat{r}$)
and $T_0$ if necessary, we obtain
$$
\frac{|\tau|-1}{2|\tau-1|}
=\frac{\Re (q t) +O(|t^2|) }{2|qt+O(|t^2|)|}
\ge 
\frac{q |t|\cos A_0 \cdot (1 +O(|t|))}{2q|t| \cdot (1 +O(|t|))}
\ge 
\frac{\cos A_0}{2} \cdot \frac{3}{4}
$$
and 
$|O(W^{-1})|\le (\cos A_0)/8$. Hence 
$$
|G(W)-B| \ge |W-B| 
+  \paren{\frac{3}{8}-\frac{1}{8}}\cos A_0
= |W-B| + d_0.
$$
\QED

\begin{prop}\label{prop_2-_2} 
In Case 2$^-$, there exists a unique linearizing coordinate
$
L: \D(B, 3|B|/4) \to \C
$
of the repelling fixed point $B^\ast$ of $G$ 
such that $L(B^\ast)=0$, $DL(B^\ast)=1$, 
and if $G(W) \in \D(B, 3|B|/4)$, 
$$
L \circ G (W)
=\tau^\ast L(W)
$$
where $\tau^\ast=DG(B^\ast)$.
\end{prop}

\paragraph{Proof.}
By the previous proposition, 
the disk $\D(B, 3|B|/4)$ is compactly contained in
$G(\D(B, 3|B|/4))$. Hence the univalent branch of $G^{-1}$
on $G(\D(B, 3|B|/4))$ is strictly contracting. 
By the Riemann mapping theorem and the Schwarz lemma, 
there exists a unique attracting fixed point 
of $G^{-1}$ in $\D(B, 3|B|/4)$, which must be $B^\ast$. 
Hence there is a unique linearizing coordinate 
$L$ near $B^\ast$ that satisfies 
$L(B^\ast)=0$, $DL(B^\ast)=1$, and the relation
$L \circ G^{-1}(W)
=(\tau^\ast)^{-1} L(W)$.
We obtain the desired linearizing coordinate
by extending $L$ to $\D(B, 3|B|/4)$ by this relation.
\QED 

\paragraph{}
Now we are ready to finish the proof of Proposition \ref{prop_leaving_U_0} in Case 2$^-$.
Take any 
$z_0 \in J(f_{c_t}) \cap U_0$ for a given $t \in \Delta$.
Let $W_0 := \Phi_t (z_0)$ and $W_k:=G^k(W)$ for $k \ge 0$.

If $|W_0-B|\le |B|/2$, 
Proposition \ref{prop_2-_2} implies that 
$|W_k-B|> |B|/2$ for some $k \ge 1$ unless $W_0=B^\ast$.
If $W_0=B^\ast$, then $z_0$ 
is a repelling fixed point of $f_c^{pq}$ in $U_0$.
If $W_0 \neq B^\ast$, 
by Proposition \ref{prop_2-_1},
either $|W_{k+j} -B| \ge |W_k-B|+jd_0
\to \infty $ 
as $j \to \infty$ or $|W_{k+j}| \le \hat{r}$
for some $j \ge 0$. 
The former implies that the orbit of $z_0$ 
is attracted to an attracting fixed point of $f_c^{pq}$
(see Proposition \ref{prop_G_near_infinity}) 
and thus a contradiction.
The latter implies that $z_{pq(k+j)}$ is not
contained in $U_0$ any more.
\QED  

\paragraph{Proof of Proposition \ref{prop_leaving_U_0}: Cases 1 and  2$^+$.}
In this case $\tau=1 - qt+O(t^2)$ and thus $B=1/(qt)+O(1)$.
The proof is analogous to Case 2$^-$ and 
we only give an outline.

Since $G^{-1}(W)-B=\tau^{-1}(W-B) +O(W^{-1})$, we have
\begin{equation}\label{eq_2+_1}
|G^{-1}(W)-B| \ge |W-B| +d_0
\end{equation}
when $|W| \ge \hat{r}$ and $|W-B| \ge |B|/2$.
The fixed point $B^\ast$ of $G^{-1}$ is repelling 
and we can find a linearizing coordinate $L$
of $G^{-1}$ in $\D(B, 3|B|/4)$.

Let $z_0 \in J(f_{c_t}) \cap U_0$ for a given $t \in \Delta$.
Let $W_0 := \Phi_t(z_0)$
and $W_k:=G^k(W_0)$ for $k \ge 0$.
If $W_0=\infty$, then $z_0$ is the repelling fixed point 
of $f_{c_t}^{pq}$ (see Proposition \ref{prop_G_near_infinity}). 
Suppose that $W_0 \neq \infty$.  
Then $W_0$ is not contained in $\D(B, |B|/2)$,
otherwise $W_k$ tends to $B^\ast$ as $k \to \infty$.
By \eqref{eq_2+_1}, 
we have $|W_0-B| \ge|W_k-B| + kd_0 $
and thus
either $|W_{k}-B| < 3|B|/4$
or $|W_{k}| \le \hat{r}$
for some $k \ge 0$. 
The former implies that $z_{pqk}$
is contained in the attracting basin 
of an attracting fixed point of $f_c^{pq}$,
and thus a contradiction.
The latter implies that $z_{pqk}$ is not
contained in $U_0$ any more.
\QED

\section{Proof of Lemma H} \label{sec:Proof_H}
In this section we give a proof of Lemma H.
The proof employs the branched coordinates 
and the dynamics of the form $G(W)=\tau W+1+O(1/W)$
as in the previous section.

\paragraph{Dynamics of $G$ near the origin.}
Since $\tau \to 1$ as $t \in \Delta$ tends to $0$,
the dynamics of $G$ (relatively) near the origin becomes closer to 
the translation $W \mapsto W+1$:

\begin{prop} \label{prop_G_near_origin}
For any $t \in \Delta$ and 
any $W$ with $\hat{r} \le |W| \le 1/(5 q |t|)$, 
we have $|G(W)-(W+1)| \le 1/2$
by taking sufficiently small $\hat{R}$ and $T_0$
in the definitions of $\hat{U}$ and $\Delta$.
In particular, we have 
\begin{equation}
1/2 \le \Re G(W) -\Re W \le 3/2
\quad\text{and}\quad
|\arg (W-G(W))| \le \pi/6.
\label{eq_G_near_origin}
\end{equation}
\end{prop}

\paragraph{Proof.}
We have $|(\tau-1)W| =|(\mp qt +O(t^2))W| \le 1/4$
when $|W| \le (5 q |t|)^{-1}$ by taking a smaller $T_0$.
We also have 
$|G(W)-(\tau W+1)|=|O(1/W)| \le 1/4$
when $|W| \ge \hat{r}$ 
by taking a smaller $\hat{R}$ (hence a larger $\hat{r}$).
It follows that $|G(W)-(W+1)| = |(\tau-1)W+O(1/W)| \le 1/2$
and \eqref{eq_G_near_origin} is an immediate consequence.
\QED

\paragraph{Some additional conditions 
on $U_0$ in $W$-coordinate.}
In what follows we suppose that $z_0 \in J(f_{c_t}) \cap U_0$,
and there exists an $m \in \N$ such that 
$z_{(m-1)pq} \in U_0$ but $z_{mpq} \notin U_0$.
Let $W_k:=\Phi_t(z_{kpq})$.
Then we may assume that 
$W_{m-1} \in \Phi_t(U_0)$ but $W_{m} \notin \Phi_t(U_0)$.

By taking smaller and appropriate $\hat{R}$ and $R_0$
in the definitions of $\hat{U}$ and $U_0$
(and taking a smaller $T_0$ if necessary), 
we may assume that the set
$\Phi_t(\partial U_0)$ is contained in the annulus
$$
\hat{A}:=\brac{W \in \Cbar_W \st 5 \hat{r} < |W| <  6 \hat{r}}
$$ 
(where $\hat{r} \asymp \hat{R}^{-q}$ is introduced in the previous section)
for any $t \in \Delta$.
Indeed, we can apply the Koebe distortion theorem
to the family of local coordinates $\{\psi_{t}\}$
to control the shape (eccentricity) of
the image $\Phi_t(\partial U_0)$. 
Because $\Re W_{m}-\Re W_{m-1} \le 3/2$
by Proposition \ref{prop_G_near_origin}, 
we may assume in addition that $W_m$ is contained 
in this $\hat{A}$.

\begin{figure}[htbp]
\begin{center}
\includegraphics[width=.6\textwidth]{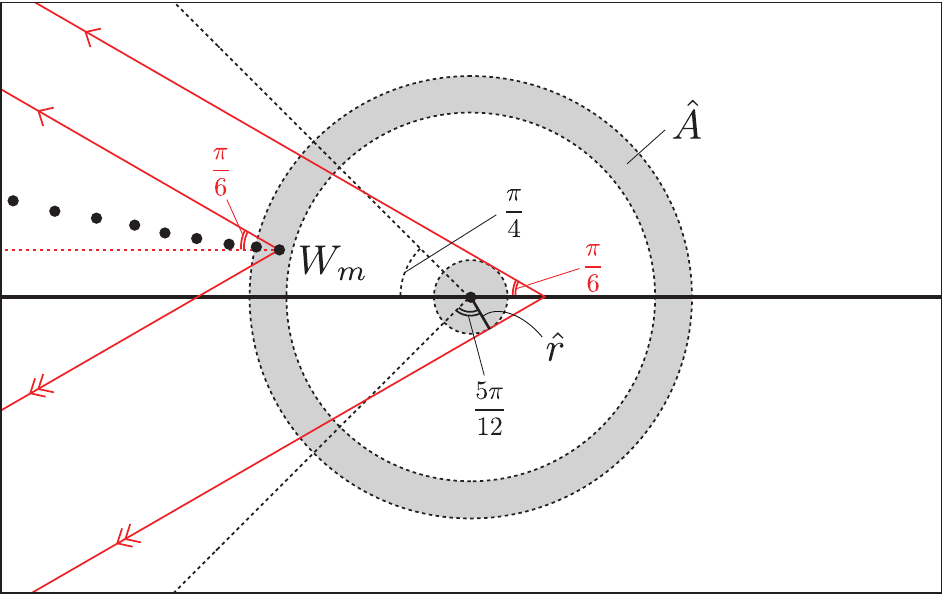}
\caption{The orbit leaving the domain $\Phi_t(U_0)$ of $G$.}
\label{fig_W_m}
\end{center}
\end{figure}

Now we claim that $3 \pi/4 \le \arg W_m \le 5\pi/4$,
equivalently, $|\arg (-W_m)| \le \pi/4$.
Indeed, since $z_0$ is contained in the Julia set,
the same argument as Proposition \ref{prop_leaving_U_0}
yields that the orbit of $W_0$ by $G$ 
must leave the domain of $G$
(by taking a larger $\hat{r}$ if necessary). 
In other words, there exists an integer $k >m$ such that 
$|W_k|>\hat{r}$ and $|W_{k+1}| \le \hat{r}$.
Since $|\arg (W_k-W_{m})| \le \pi/6$
by the proposition above, the condition
$|W_m| > 5\hat{r}$ implies $|\arg (-W_m)| \le \pi/4$.
(Here we have used $\cos (5 \pi/12) > 1/5$. 
See Figure \ref{fig_W_m}.)

\paragraph{}
The lemma below is a very important conclusion from the radial access condition:

\paragraph{Lemma I.} 
{\it 
There exists a constant $K_{\mathrm I}\in (0,  1/(5q)]$ such that
for any $t \in \Delta$, $z_0 \in J(f_{c_t})\cap U_0$
and $k \ge 0$, the condition
$|W_k| \le K_{\mathrm I}/|t|$ implies $|\arg (-W_k)| \le \pi/4$. Moreover, for such a $W_k$, we have $|W_{k+1}|<|W_k|$. 
}

\paragraph{}

Hence we can only find $W_k$ near the negative real axis
in the disk $\D(0, K_{\mathrm I}/|t|)$.  

\paragraph{Proof.}
By the assumption on $|W_m|$ as above, we may assume that there exists  $j \in \N$
such that 
$$
W_{m-j},\,W_{m-j+1}, \cdots, W_m 
\in \brac{W \in \Cbar_W 
\st 
|W| \le (5 q |t|)^{-1}, \, 
|\arg (-W)| \le \pi/4}
$$
and that $|W_{m-j}| \ge  (6 q |t|)^{-1}$.
Note that the condition $|\arg (-W)| \le \pi/4$ implies 
$
|W|/\sqrt{2} \le -\Re W \le |W|.
$

By Proposition \ref{prop_G_near_origin}, we have
$
\Re W_{m} -\Re W_{m-j} \le 3j/2
$
and thus 
$$
3j/2 \ge -|W_m| + |W_{m-j}|/\sqrt{2} \ge -|W_m| + (6\sqrt{2}q|t|)^{-1}.
$$
Since the assumption $5\hat{r}<|W_m|<6\hat{r}$ implies $|W_m| \asymp 1$,
we conclude that $j \ge C_6/|t|$ for some constant $C_6\asymp 1$
independent of $t \in \Delta$ (taking a smaller $T_0$ 
if necessary). 

In Case 2$^-$, as in the proof of Proposition  \ref{prop_leaving_U_0},  there is $k'>j\ge1$ for which $|W_{m-k'}-B|>|B|/2$ (otherwise $W_0=B^*$). We have
$$
|W_m-B| \ge |W_{m-k'}-B| + k'd_0 
> |W_{m-k'}-B| + C_6d_0/|t|  
$$
by Proposition \ref{prop_2-_1}.
Hence 
$$
|W_{m-k'}| \ge |B|-|W_{m-k'}-B|
> |B|-|W_m-B|+C_6d_0/|t|
\ge -|W_m|+C_6d_0/|t|.
$$
Since $|W_m| \asymp 1$, there exists a constant $K_{\mathrm I}\asymp 1$
such that $|W_{m-k'}|> K_{\mathrm I}/|t|$ for any such $k'$  and 
$t \in \Delta$ (again taking a smaller $T_0$ if necessary). 
Hence if we have $|W_k| \le K_{\mathrm I}/|t|$ for some $k$,
then $m-j \le k \le m$ and thus $|\arg (-W_k)| \le \pi/4$. 
Since $|W_k|\le K_{\mathrm I}/|t|$ also implies $|W_k|\le 1/(5q|t|)$,
 we conclude $K_{\mathrm I}\le 1/(5q)$.

The proof for  Cases 1 and 2$^+$ is analogous:
By \eqref{eq_2+_1}, we have 
$$
|W_{m-k'}-B| \ge |W_{m}-B| + k'd_0 
> |W_{m}-B| + C_6d_0/|t|  
$$
instead, and this implies 
$$
|W_{m-k'}| \ge |W_{m-k'}-B|-|B|
 > |W_{m}-B| -|B| + C_6d_0/|t| 
\ge -|W_m|+ C_6d_0/|t|. 
$$
Then we repeat the same argument as above.

Now, we show the condition $|W_k|\le K_{\mathrm I}/|t|$ implies $|W_{k+1}|<|W_k|$. By Proposition \ref{prop_G_near_origin}, we get $|W_{k+1}-(W_k+1)|\le 1/2$ and thus 
$$ 
\left| \frac{W_{k+1}}{W_k}\right|\le \left|1+\frac{1}{W_k}\right|+\frac{1}{2|W_k|}.
$$
Since $|\arg (-W_k)|\le \pi /4$ implies $3\pi /4\le \arg (1/W_k)\le 5\pi/4$, we obtain
$$
\left|1+\frac{1}{W_k}\right| \le 
\sqrt{
\paren{
1-\frac{1}{\sqrt{2}~|W_k|}}^2+
\paren{\frac{1}{\sqrt{2}~|W_k|}}
^2}
$$
and thus
\begin{eqnarray*}
\left| \frac{W_{k+1}}{W_k}\right|
&\le&  1- \frac{1}{\sqrt{2}~W_k}+\frac{1}{2~|W_k|}+O\left(\frac{1}{|W_k|^2}\right) \\
&=& 1- \frac{\sqrt{2}-1}{2~|W_k|}+O\left(\frac{1}{|W_k|^2}\right) \\
&<& 1
\end{eqnarray*}
by taking a sufficiently large $\hat{r}$ if necessary.
\QED

\medskip

\begin{remark}
 Once we have $|W_k|\le K_{\mathrm I}/|t|$, then $|W_j|\le K_{\mathrm I}/|t|$ for $k\le j\le m$.
\end{remark}

\begin{prop}\label{prop_DG_k}
For any $t \in \Delta$, any $z_0 \in J(f_{c_t}) \cap U_0$
that is not a fixed point of $f_{c_t}^{pq}$,
and any $k$ with $0 \le k \le m$, we have 
$$
|DG^k(W_0)| \asymp |\tau|^k 
$$
in Case 1 or Case 2$^+$, and 
$$
|DG^k(W_0)| \asymp |\tau^\ast|^k 
$$ 
in Case 2$^-$, where $W_0=\Phi_t(z_0)$. 
The implicit constants are independent of $t, z_0$, and $k$.
\end{prop}

\paragraph{Proof for Case 1 and Case 2$^+$.}
Since $DG(W)=\tau +O(W^{-2})$, 
\begin{equation}
DG^k(W_0)=\tau^k \, \prod_{j=0}^{k-1}\paren{1+O({W_j}^{-2})}. \label{DGk=}
\end{equation}
Hence it is enough to show that $\sum_{j=0}^{k-1}|W_j|^{-2}$
is uniformly bounded in $t$, $W_0$, and $k$. 

As in the proof of Proposition   \ref{prop_leaving_U_0}, $\D(B, 3|B|/4)$ is contained in the attracting basin of $B^*$. Therefore, we may always assume $|W_j|\ge \hat{r}$ and $|W_j-B|>|B|/2$ for $j=0, 1, \ldots, k$. 
By \eqref{eq_2+_1}, 
$|W_j-B|$ is strictly decreasing in $j$. 
This implies that if $|W_{k_1}-B|< 4 |B|$ for some $0\le k_1\le k$, 
then $|W_j-B|< 4|B|$ for all $k_1\le j\le k$. 
On the other hand, by Lemma I, if $|W_{k_2}|\le K_{\mathrm I}/|t|$ for some $0\le k_2\le k$, then $|W_j|\le K_{\mathrm I}/|t|$ and $|\arg (-W_j)|\le \pi/4$ for all  $k_2\le j\le k$ $(\le m)$. Consequently,  we can assume the location of 
 $W_0, \ldots, W_k$  to be 
\begin{enumerate}[(i$^+$)]
\item
$|W_j-B| \ge 4 |B|$  for $j=0, \ldots, k_1-1$,
\item
$|W_j-B| < 4|B|$ and $|W_j| > K_{\mathrm I}/|t|$ for $j=k_1, \ldots, k_2-1$,
\item
$|W_j-B| < 4 |B|$ and $|W_j| \le K_{\mathrm I}/|t|$ for $j=k_2, \ldots, k$.
\end{enumerate}
 
For the case (i$^+$), by   \eqref{eq_2+_1} again we have
$$
|W_j| \ge |W_j-B| -|B| 
\ge |W_{k_1-1}-B| -|B| +(k_1-1-j)d_0 
\ge |B| +(k_1-1-j)d_0. 
$$
Hence 
\[
\displaystyle \sum_{j=0}^{k_1-1}|W_j|^{-2}
=O\paren{\sum_{j=0}^{k_1-2}\frac{1}{(k_1-1-j)^2}}+\frac{1}{|B|^2}
=O\paren{\sum_{j=1}^{\infty}\frac{1}{j^2}}+O\paren{|t|^2}
=O(1).
\]

The condition $|W_j-B|< 4|B|$ implies 
$|W_j| < 5 |B| = 5/(q|t|)+O(1)$.
Hence the case (ii$^+$) implies 
 $|W_j| \asymp |t|^{-1}$ and thus $|W_j|^{-2} \asymp |t|^2$. From  \eqref{eq_2+_1}, we obtain 
$$
4|B| > |W_{k_1}-B| \ge |W_{k_2-1}-B| + (k_2-1-k_1)d_0 \ge (k_2-1-k_1) d_0
$$
and thus $k_2-1-k_1 =O(|t|^{-1})$. 
This gives
\[ \sum_{j=k_1
             }^{k_2-1}|W_j|^{-2}=O(|t|). 
\]

 In the  (iii$^+$) case, we have  $|\arg (-W_j)| \le \pi/4$  and thus 
$|W_j|/\sqrt{2} \le -\Re W_j \le |W_j|$.
By Proposition \ref{prop_G_near_origin},
we have 
$
 -\Re W_j \ge -\Re W_k+(k-j)/2.
$
Hence,  we get $|W_j| \ge |W_k|/\sqrt{2} +(k-j)/2 
\ge (k-j)/2$, and 
\[
\displaystyle \sum_{j=k_2
             }^{k-1}|W_j|^{-2}
=O\paren{\sum_{j=k_2}^{k-1}\frac{1}{(k-j)^2}}
=O\paren{\sum_{j=1}^{\infty}\frac{1}{j^2}}
=O(1).
\]

\paragraph{Proof for Case 2$^-$.}
First note  by Proposition \ref{prop_G_near_infinity}
that $|W_j|\le 1/(qd_0|t|)$ for $0\le j\le k$
(otherwise $W_j$ must be attracted by $\infty$). Second, by Proposition \ref{prop_2-_1}, if $|W_{k_1}-B|\ge |B|/2$ for some  least $k_1\le k$, then $|W_j-B|$ is strictly increasing for $j$ with $k_1\le j\le k$.
 Hence,
 $|W_j-B|\ge |B|/2$ for such a $j$. Furthermore, 
if $|W_{k_2}|\le K_{\mathrm I}/|t|$ for some $0\le k_2\le k$, then as in the Case 1 and Case 2$^+$,  the location of 
 $W_0, \ldots, W_k$  must be 
\begin{enumerate}[(i$^-$)]
\item
$|W_j-B| < |B|/2$  for  $j=0, \ldots, k_1-1$,
\item
$|W_j-B| \ge |B|/2$ and $|W_j| > K_{\mathrm I}/|t|$ for $j=k_1, \ldots, k_2-1$,
\item
$|W_j-B| \ge |B|/2$ and $|W_j| \le K_{\mathrm I}/|t|$ for   $j=k_2, \ldots, k$.
\end{enumerate}
 
Now, 
\begin{eqnarray*}
\lefteqn{
\frac{DG^k(W_0)}{(\tau^\ast)^k}}\\
&=&\frac{DG^{k_1-1}(W_0)}{(\tau^\ast)^{k_1-1}}  \cdot \frac{DG^{k-k_1+1}(W_{k_1-1})}{(\tau^\ast)^{k-k_1+1}} \\
&\asymp& 1\cdot  \frac{DG^{k-k_1+1}(W_{k_1-1})}{(\tau^\ast)^{k-k_1+1}} \qquad\mbox{(by Proposition \ref{prop_2-_2})} \\
&\asymp&   \paren{\frac{\tau}{\tau^\ast}}^{k-k_1+1}
\prod_{j=k_1-1}^{k-1}\paren{1+O({W_j}^{-2})} \qquad\mbox{(as \eqref{DGk=})}\\
&\asymp &  \paren{1+O(t^2)}^{k-k_1+1}\cdot 
\prod_{j=k_1-1}^{k-1}\paren{1+O({W_j}^{-2})} \qquad\mbox{(since $\displaystyle\frac{\tau}{\tau^\ast}=1+O(t^2)$.)} 
\end{eqnarray*}

Because $|W_{k_1-1}| \ge |B|-|W_{k_1-1}-B|> |B|-|B|/2=|B|/2=O(1/|t|)$ and because $|W_{k_1-1}|\le 1/(qd_0|t|)$, we have $|W_{k_1-1}^{-2}|=O(|t|^{2})$.
 
By Proposition \ref{prop_2-_1}, we know 
$$
|W_k-B| \ge |W_{k_1}-B| + (k-k_1)d_0 \ge (k-k_1) d_0
$$
and thus $k-k_1  \le (|W_k| + |B|)/d_0 =O(|t|^{-1})$.
This implies that $(1+O(t^2))^{k-k_1+2}=1+O(t)$.

It remains to obtain an  estimate of $\prod_{j=k_1}^{k-1}\paren{1+O(W_j^{-2})}$. 
Again, it is enough to obtain an estimate by showing that $\sum_{j=k_1}^{k-1}|W_j|^{-2}$ is uniformly bounded in $k$ and $W_0$. 

For the (ii$^-$) case, since $|W_j|<1/(qd_0|t|)$,   
we obtain $|W_j| \asymp |t|^{-1}$.
Therefore, 
\[ \sum_{j=k_1}^{k_2-1}|W_j|^{-2}=O(|t|) 
\qquad \mbox{(since $k-k_1 = O(|t|^{-1})$ so is $k_2-k_1$.)}
\]
The  (iii$^-$) case follows exactly as the (iii$^+$) case. 

The proof of Proposition \ref{prop_DG_k} is complete.
\QED

\paragraph{Lemma H in the branched coordinates.}
To give estimates for the sums of the form
$$
\sum_{k=0}^{m-1} 
\frac{1}{|Df_c^{kpq}(z_j)|}
\quad\text{and}\quad
\sum_{k=0}^{m-1} 
\frac{|c-\chat|+| f_c^{kpq}(z_j) -b_j(c)|}{ |Df_c^{kpq}(z_j)|}
$$
with $0\le j\le p-1$, we rewrite them in $W$-coordinates.
It is enough to consider the case of $j=0$,
since we may apply the same argument as in the proof of 
Lemma G.

Let $W=\Phi_t(z)=\Psi_t \circ \psi_{t}(z)$
for $z \in \hat{U}$ such that 
$W_k:=\Phi_t(z_{kpq})$ for each $k$.
Note that $|D\psi_t(z)| \asymp 1$ for any $z \in \hat{U}$,
where the implicit constant is independent of $t \in \Delta$.
Since $D\Psi_t(w) =\lam_t^{q^2}/w^{q+1}$, we have 
$$
|D\Phi_t(z)|=|D\Psi_t(w)\cdot D\psi_t(z)| 
\asymp |w|^{-(q+1)} \asymp |W|^{1+1/q}. 
$$
By the chain rule, we obtain
\begin{align*}
|Df_{c_t}^{kpq}(z_0)|
&=|D(\Phi_t^{-1} \circ G^k \circ \Phi_t)(z_0)|
=|D\Phi_t^{-1}(W_k)\cdot D G^k(W_0) \cdot D\Phi_t(z_0)|\\
&\asymp |W_k|^{-1-1/q}\cdot|DG^k(W_0)|\cdot|W_0|^{1+1/q}\\
&= \abs{\frac{W_0}{W_k}}^{1+1/q}\cdot|DG^k(W_0)|.
\end{align*}
We let
$$
\cS_1
:=
\sum_{k=0}^{m-1} 
\abs{\frac{W_k}{W_0}}^{1+1/q}\cdot\frac{1}{|DG^k(W_0)|}
\asymp
\sum_{k=0}^{m-1} 
\frac{1}{|Df_{c_t}^{kpq}(z_0)|},
$$
and 
$$
\cS_2:=
\frac{1}{|W_0|^{1/q}}
\sum_{k=0}^{m-1}  
\abs{\frac{W_k}{W_0}}\cdot\frac{1}{|DG^k(W_0)|}
\asymp
\sum_{k=0}^{m-1} 
\frac{|f_{c_t}^{kpq}(z_0) -b_t|}{ |Df_{c_t}^{kpq}(z_0)|},
$$
where we used the fact that  
$|f_c^{kpq}(z_0) -b_0(c)|=|z_{kpq} -b_t|\asymp |W_k|^{-1/q}$.
Now Lemma H is reduced to the estimates
$$
\cS_1 = O\paren{\frac{1}{|t|}}
\quad\text{and}
\quad
\cS_2 = O\paren{\frac{1}{|t|^{1-1/q}}}. 
$$
Indeed, Proposition \ref{prop_c_and_epsilon} implies
$$
\sum_{k=0}^{m-1} 
\frac{1}{|Df_{c_t}^{kpq}(z_0)|}
\asymp \cS_1 
= O\paren{\frac{1}{|t|}}
=
O\paren{\frac{1}{\sqrt{|c_t-\chat|}}}
$$ 
in Case 1, and
\begin{align*}
\sum_{k=0}^{m-1} 
\frac{|c_t-\chat|+|f_{c_t}^{kpq}(z_0) -b_t|}{ |Df_{c_t}^{kpq}(z_0)|}
&\asymp
|c_t-\chat|\cS_1+\cS_2 \\
&=O\paren{ |t|\cdot \frac{1}{|t|} + \frac{1}{|t|^{1-1/q}}}
=O\paren{\frac{1}{|c_t-\chat|^{1-1/q}}}
\end{align*}
in Case 2.

\paragraph{Proof of Lemma H for  Case $2^-$.}
We will use tha fact that 
$\tau^\ast=1+qt+O(t^2)$ and 
$$
|\tau^\ast|= 1+\Re (qt) + O(|t^2|) \ge 1 + q d_0 |t|
$$
for $t \in \Delta$ if we take a smaller $T_0$ if necessary.

First, suppose that $|W_0| \ge K_{\mathrm I}/|t|$.
Since it must be $|W_k| \le 1/( q d_0 |t|)$
for any $k=0,1,\ldots, m$ 
by Proposition \ref{prop_G_near_infinity}
(otherwise $W_k$ is attracted by $\infty$),
we have $|W_0| \asymp |t|^{-1}$
and $|W_k/W_0|=O(1)$. 
Hence by Proposition \ref{prop_DG_k}, we have 
$$
\cS_1 
=
\sum_{k=0}^{m-1} 
\abs{\frac{W_k}{W_0}}^{1+1/q}\cdot\frac{1}{|DG^k(W_0)|}
=
\sum_{k=0}^{m-1} O(1) \cdot \frac{1}{|\tau^\ast|^k}
=
O\paren{
\sum_{k=0}^{m-1} \frac{1}{(1+qd_0|t|)^k}}
=O\paren{\frac{1}{|t|}}
$$
and
$$
\cS_2
=
\frac{1}{|W_0|^{1/q}}
\sum_{k=0}^{m-1} 
\abs{\frac{W_k}{W_0}}\cdot\frac{1}{|DG^k(W_0)|}
= 
|t|^{1/q}
\sum_{k=0}^{m-1} O(1) \cdot \frac{1}{|\tau^\ast|^k}
=O\paren{\frac{1}{|t|^{1-1/q}}}.
$$

Next, suppose that $|W_0| \le K_{\mathrm I}/|t|$.
By Proposition 
\ref{prop_G_near_origin}
and Lemma I,
we have $|\arg (-W_k)| \le \pi/4$
and thus 
$|W_k|/\sqrt{2} \le -\Re W_k \le |W_k|$ 
for any $k=0,1,\ldots, m$.
By Proposition
\ref{prop_G_near_origin} again we have
$\Re W_0 + k/2 \le \Re W_k$. 
Hence
\begin{equation}
0 \le |W_k|/\sqrt{2} 
\le -\Re W_k \le -\Re W_0 - k/2
\le |W_0|-k/2  \label{W0k/2}
\end{equation}
and this implies 
$$
\abs{\frac{W_k}{W_0}}
\le \sqrt{2} \cdot \frac{|W_0|-k/2}{|W_0|}=O(1).
$$
In particular, by letting $k=m$ in \eqref{W0k/2} we have 
$m \le 2|W_0|=O(|W_0|)$. 

Since $|W_0| \le K_{\mathrm I}/|t|$, 
we have
$$
\cS_1 
=
\sum_{k=0}^{m-1} O(1) \cdot \frac{1}{|\tau^\ast|^k}
\le
\sum_{k=0}^{m-1} O(1) \cdot 1^k
=O(m) = O(|W_0|)=O\paren{\frac{1}{|t|}}
$$
and
\begin{align*}
\cS_2
&= 
\frac{1}{|W_0|^{1/q}}
\sum_{k=0}^{m-1} 
O(1) \cdot \frac{1}{|\tau^\ast|^k}
\le 
\frac{1}{|W_0|^{1/q}}
\sum_{k=0}^{m-1} 
O(1) \cdot 1^k\\
&=\frac{O(m)}{|W_0|^{1/q}}
= O(|\,W_0|^{1-1/q}|)
=O\paren{\frac{1}{|t|^{1-1/q}}}.
\end{align*}
This completes the proof for Case 2$^-$.

\paragraph{Preliminary for Case 1 and Case $2^+$.}
We will use the fact that 
$\tau^{-1}=1 +qt+O(t^2)$ and 
$$
|\tau|^{-1}= 1 + \Re (qt) + O(|t^2|) \ge 1 + q d_0 |t|
$$
for $t \in \Delta$ if we take a smaller $T_0$ if necessary.

It is very convenient to use {\it Ueda's modulus} defined 
as follows (see \cite{U}): For $W \neq \infty$, 
$$
N(W):=|W-B|-|B|.
$$
Then one can easily check that 
\begin{itemize}
\item $|W| \ge N(W)$; and
\item $N(W) \ge |W|/3$ when $|W-B| \ge 4 |B|$.
\end{itemize}
Indeed, $|W| \ge N(W)$ is just the triangle 
inequality; and when $|W-B| \ge 4 |B|$, we have 
$|W| \ge 3|B|$ and thus $N(W)/|W|=|1-B/W|-|B/W| \ge 1/3$.

Here is another useful fact:
\begin{prop}\label{prop_N}
In Case 1 and Case 2$^+$,
we have 
$$
N(G^{-1}(W)) \ge |\tau|^{-1} N(W)+d_0
$$
by taking sufficiently small $\hat{R}$ and $T_0$
in the definitions of $\hat{U}$ and $\Delta$.
\end{prop}

\paragraph{Proof.}
Since $|G^{-1}(W)-B| =|\tau|^{-1}|W-B|+O(|W^{-1}|)$,
we have 
$$
N(G^{-1}(W))
=|\tau|^{-1}N(W) + \frac{|\tau|^{-1}-1}{|1-\tau|}+O(|W^{-1}|).
$$
By taking a smaller $T_0$ if necessary, 
we have 
$$
\frac{|\tau|^{-1}-1}{|1-\tau|} 
= \frac{\Re t}{|t|}(1+O(|t|)) 
\ge \frac{\cos A_0}{2} = 2 d_0
$$ 
for any $t \in \Delta$. 
By taking a smaller $\hat{R}$
(hence a larger $\hat{r}$) if necessary, 
we have $O(|W^{-1}|) \le (\cos A_0)/4=d_0$
for any $|W| \ge \hat{r}$.
Hence we have the desired inequality.
\QED

\begin{prop}\label{prop_Lem_P_1}
In Case 1 and Case 2$^+$, we have
$$
\abs{\frac{W_k}{W_0}}=O(|\tau|^k)
$$
for any $k=0, 1, \ldots, m$.
\end{prop}

\paragraph{Proof.}
We adopt the argument of Proposition \ref{prop_DG_k}
and assume that 
\begin{enumerate}[(i$^+$)]
\item
$|W_j-B| \ge 4 |B|$  for $j=0, \ldots, k_1-1$,
\item
$|W_j-B| < 4|B|$ and $|W_j| > K_{\mathrm I}/|t|$ for $j=k_1, \ldots, k_2-1$,
\item
$|W_j-B| < 4 |B|$ and $|W_j| \le K_{\mathrm I}/|t|$ for $j=k_2, \ldots, k$
\end{enumerate}
for some $k_1$ and $k_2$.
Now we consider the product
\begin{equation}\label{eq_W_kW_0}
\abs{\frac{W_k}{W_0}}
= 
\abs{\frac{W_{k_1-1}}{W_{0}}}
\cdot 
\abs{\frac{W_{k_2}}{W_{k_1-1}}}
\cdot 
\abs{\frac{W_k}{W_{k_2}}}.
\end{equation}
In the case of (i$^+$), 
we can apply Proposition \ref{prop_N} and
it follows that
$$
|W_0| \ge N(W_0) \ge |\tau|^{-(k_1-1)}N(W_{k_1-1})
\ge |\tau|^{-k_1+1} |W_{k_1-1}|/3.
$$
Hence $|W_{k_1-1}/W_0| \le 3 \,|\tau|^{k_1-1}$.
Since $|W_{k_1-1}| \ge  3|B|=3/(q|t|)+O(1)$ and 
$|W_{k_2}| \le K_{\mathrm I}/|t|$,
 we have $|W_{k_2}/W_{k_1-1}| =O(1)$. 
Since $|W_j|$ is decreasing in the case of (iii$^+$) by Lemma I, 
we have $|W_{k}/W_{k_2}| \le 1$. 
Hence by \eqref{eq_W_kW_0}, we obtain $|W_{k}/W_0|=O(|\tau|^{k_1-1})$.

To show the estimate $|W_{k}/W_0|=O(|\tau|^{k})$, 
it is enough to show that 
$|\tau|^{k-(k_1-1)}=O(1)$.
Because $|W_{k_1}-B| \le 4 |B|$
and hence $|W_{k_1}|\le 5|B|$,
by \eqref{eq_2+_1} we have 
$$
|W_{k_1}-B| \ge |W_k-B| +(k-k_1)\,d_0
\ge (k-k_1)\,d_0
$$
and thus 
$k-k_1 \le (|W_{k_1}|+|B|)/d_0 =O(|t|^{-1})$.
Hence 
$$
|\tau|^{k-(k_1-1)} 
=(1-\Re(qt)+O(t^2))^{k-k_1}\cdot |\tau|=O(1)
$$
and this completes the proof.
\QED

\paragraph{Proof of Lemma H for Case 1 and Case $2^+$.}
By Propositions  \ref{prop_DG_k} and 
\ref{prop_Lem_P_1}, we have
$$
\cS_1 
=
\sum_{k=0}^{m-1} 
\abs{\frac{W_k}{W_0}}^{1+1/q}\cdot\frac{1}{|DG^k(W_0)|}
=
\sum_{k=0}^{m-1} O(|\tau|^{k(1+1/q)}) 
\cdot \frac{1}{|\tau|^k}
=O\paren{
\sum_{k=0}^{m-1} |\tau|^{k/q}}.
$$
Since $|\tau|^{1/q}\le 1-d_0|t|$
for $t \in \Delta$ (by taking a smaller $T_0$
if necessary),
we obtain $\cS_1=O(|t|^{-1})$.

Similarly we have
$$
\cS_2
=
\frac{1}{|W_0|^{1/q}}
\sum_{k=0}^{m-1} 
\abs{\frac{W_k}{W_0}}\cdot\frac{1}{|DG^k(W_0)|}
=
\frac{1}{|W_0|^{1/q}}
\sum_{k=0}^{m-1} 
O(|\tau|^{k})\cdot\frac{1}{|\tau|^k}
=
\frac{O(m)}{|W_0|^{1/q}}.
$$

First we suppose that $|W_0 - B| < 4 |B|$.
Then we have  $|W_0| < 5 |B| =O(|t|^{-1})$.
By \eqref{eq_2+_1},
$$
|W_0| \ge |W_0-B| - |B| 
\ge |W_m-B| - |B| +m d_0 
\ge -|W_m|+md_0.
$$
By the assumption that $W_m \in \hat{A}$,
we have $|W_m| < 6 \hat{r}$ and thus $md_0 \le |W_0| + 6\hat{r}$. 
We obtain $m =O(|W_0|)$ since $5\hat{r} < |W_0| = O(|t|^{-1})$.
Hence 
$$
\cS_2=\frac{O(m)}{|W_0|^{1/q}}=
O(|W_0|^{1-1/q}) = O\paren{\frac{1}{|t|^{1-1/q}}}.
$$

Next we suppose that $|W_0-B|\ge 4 |B|$.
There exists a $k\ge 0$ such that 
\begin{itemize}
\item
$|W_j-B| \ge 4 |B|$ for $0 \le j \le k$; and
\item
$|W_j-B| < 4 |B|$ for $k \le j \le m$.
\end{itemize}
(We may regard $k$ and $m$ as $k_1-1$ and $k$ in the proof of Proposition \ref{prop_Lem_P_1} respectively.)
Since $|W_{k}/W_0| \le 3 |\tau|^{k}$
as in the proof of Proposition \ref{prop_Lem_P_1},
and $|W_{k}| \ge 3|B| > 1/(q|t|)$
for $t \in \Delta$,
we have $ 3 |\tau|^{k} \ge  1/(q|t W_0|)$.
Hence 
$$
1< |\tau|^{-k} \le 3q|t W_0|, 
\quad\mbox{or equivalently}\quad
0< k \log |\tau|^{-1} \le \log(3q|t W_0|).
$$
This implies that 
$$
0 < k (\Re (qt) +O(t^2))  \le \log(3q|t W_0|).
$$
Since $\Re (qt) +O(t^2) \ge (q|t|\cos A_0)/2$ 
for $t \in \Delta$ (by taking a smaller $T_0$ if necessary),
we conclude
$$
k=O\paren{\frac{\log (3q|t W_0|)}{|t|}}.
$$
As in the proof of Proposition \ref{prop_Lem_P_1}, 
we have $m-k=O(1/|t|)$.
Since $|W_0| \ge 3 |B| \ge 2/(q|t|)$
(for $t \in \Delta$ by taking a smaller $T_0$ if necessary), 
we obtain $\log (3q|t W_0|) \ge \log 6 >1$
and thus 
$$
m=k+(m-k)
=O\paren{\frac{\log (3q|t W_0|)}{|t|}}+O\paren{\frac{1}{|t|}}
=O\paren{\frac{\log (3q|t W_0|)}{|t|}}.
$$
Hence we conclude that 
$$
\cS_2=\frac{O(m)}{|W_0|^{1/q}}=
O\paren{\frac{\log (3q|t W_0|)}{|t W_0|^{1/q}}}
\cdot \frac{1}{|t|^{1-1/q}}
 = O\paren{\frac{1}{|t|^{1-1/q}}},
$$
where we used the fact that 
$x^{-1/q}\log (3q x)$ is bounded if $3qx \ge  6$. \QED

\section{Proof of Lemma D}\label{sec:Proof_of_Lemma_D}
In this section we give a proof of Lemma D. 
Since the proof heavily relies on Lemma I,
we restate Lemma D in terms of the parameter $t \in \Delta \cup \{0\}$. (Note that 
we have $\chat=c_0$).

\begin{prop}\label{prop_Lemma_D_in_t}
There exist constants 
$K_{\mathrm D}>0$ and 
$\nu_0 \in (0, R_0)$
such that 
for any $\nu \in (0,\nu_0)$, 
there exists a constant 
$T_0=T_0(\nu) \in (0, 2 \cos A_0)$
such that for any
$t \in \Delta=\Delta(A_0,T_0)$, 
$\zeta \in P(f_{c_0})$, 
and $z \in J(f_{c_t})-\cV(c_t)$,
$$
|\zeta-z| \ge K_{\mathrm D} \,\nu.
$$
\end{prop}

Without loss of generality, we may assume that $K_{\mathrm D} \le 1$.

\paragraph{Proof.}
We prove it by contradiction. 
Suppose that for any
$K_{\mathrm D} >0$ and 
$\nu_0 \in (0, R_0)$,
there exists a $\nu \in (0, \nu_0)$
such that for any 
$T_0 \in (0, 2 \cos A_0)$,
there exist 
$t \in \Delta=\Delta(A_0,T_0)$, 
$\zeta \in P(f_{c_0})$, and $z \in J(f_{c_t})-\cV(c_t)$
such that 
$$
|\zeta-z|< K_{\mathrm D} \,\nu.
$$
For instance, let $K_{\mathrm D} :=1/n$, $\nu_0:=R_0/n$ for integer $n \ge 1$.
Then for any sufficiently large $n$, we may assume that $R_0/n <1$ 
and 
there exists a $\nu_n \in (0, R_0/n)$ such that  
for $T_0:= 2 \nu_n^{2q} \cos A_0 $
we can find sequences 
$t_n \in \Delta=\Delta(A_0,T_0)$,
$\zeta_n \in P(f_{c_0})$, and
$z_n \in J(f_{c_{n}})-\cV(c_{n})$ with $c_n=c_{t_n}$
such that 
$$
|\zeta_n-z_n|< \frac{1}{n} \cdot \nu_n.
$$
(Now, $\mathcal{V}(c_n)=\bigcup_{k=0}^{p-1}f_{c_n}^{-k}(V_0)$ and $V_0=\D(\bhat, \nu_n)$.)
We may assume that $\zeta_n$ and $z_n$
tends to the same limit $\hat{\zeta} \in P(f_{c_0})$
by passing through a subsequence, 
since $P(f_{c_0})$ is compact and the Julia sets $J(f_{c_n})$
are uniformly bounded. 
Moreover, there exists a $j\ge0$ 
such that $f_{c_0}^j(\hat{\zeta})$ is contained in 
$\D(\bhat, R_0/2) \subset U_0=\D(\bhat, R_0) $. 
By replacing $\zeta_n$ and $z_n$ by
$f_{c_0}^j(\zeta_n)$ and $f_{c_n}^j(z_n)$ respectively, 
we may assume that $\zeta_n \in U_0$ 
and $z_n \in U_0 - V_0$ for sufficiently large $n$.

In $w$-coordinate $w=\psi_0(z)$ given 
in Proposition \ref{prop_near_parabolic} for $t=0$,
the local dynamics of $f_{c_0}^{pq}$ in $U_0$ can be observed as
 $w \mapsto w(1+w^q+O(w^{2q}))$
and the orbit of $\psi_0(\zeta_n)$ 
by this map tends to $0$
tangentially to the attracting direction 
(that is, the set of $w$'s with $\Im (w^q)=0$ and $\Re (w^q)<0$).
This implies that (by taking sufficiently large $j$ if necessary) 
$\psi_0(\zeta_n)$ is contained in the set
$X:=\{w \in \psi_{0}(U_0) \st |\arg (-w^q)| \le \pi/4\}$.

\begin{figure}[htbp]
\begin{center}
\includegraphics{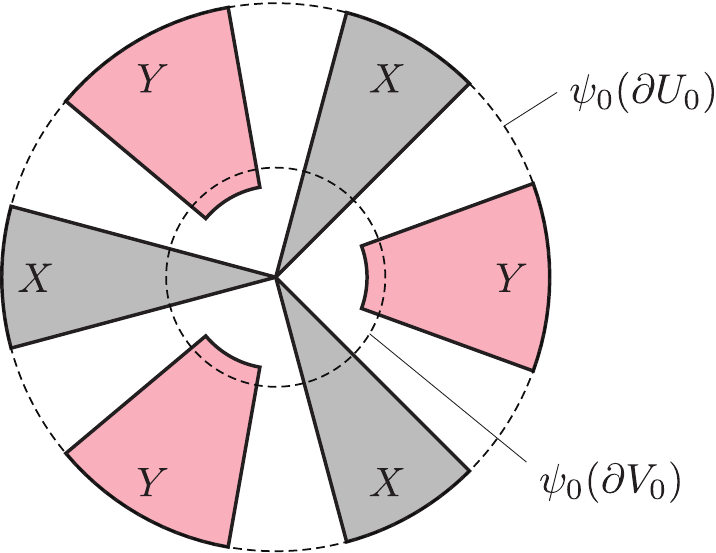}
\end{center}
\caption{The sets $X$ and $Y$ in $w$-coordinate for $q=3$.}
\label{fig_XY}
\end{figure}

Next we consider $z_n \in U_0 -V_0$ in the local coordinate $w=\psi_{t_n}(z)$. 
Since $\psi_{t_n}$ is univalent and $\nu_n < |z_n-\bhat| <R_0$,
there exists a constant $C_{7}>0$ independent of $n$
such that $w_n:=\psi_{t_n}(z_n)$ satisfies $|w_n| \ge C_{7}\nu_n$.
In $W$-coordinate $W=\Psi_{t_n}(w) =-\lam_{t_n}^{q^2}/(q w^q)$, 
we have $\arg(-W) = -\arg (w^q) +O(t_n)$ since 
$\lam_{t_n}^{q^2}=1 +O(t_n)$. Now, 
for $W_n=\Psi_{t_n}(w_n)=\Psi_{t_n} \circ \psi_{t_n}(z_n)$
with sufficiently large $n$, we obtain
\[
  |W_n|\le \left| \frac{1+O(t_n)}{qC_{7}^q \nu_n^q}\right|\le 
\left| 
            \frac{\sqrt{2\cos A_0}}
                   {qC_{7}^q}  \frac{1+O(|t_n|)}
                   {           \sqrt{|t_n|}
                   }
    \right|\le \frac{K_{\mathrm I}}{|t_n|}.
\]
Therefore, by Lemma I, 
we have $|\arg(-W_n)| \le \pi/4$ 
since $z_n \in J(f_{c_n}) \cap U_0$.
Hence in $w$-coordinate we have 
$|\arg (w_n^q)| =|\arg(-W_n)|+O(|t_n|) \le \pi/4+O(\nu_n^{2q})$ for 
$w_n=\psi_{t_n}(z_n)$. 

We also have $|\psi_{t_n}(z_n) -\psi_0(z_n)|=O(|t_n|)=O(\nu_n^{2q})$,
since the map $(t,z) \mapsto \psi_t(z)$ 
is holomorphic in both $t$ and $z$. 
Hence $\psi_0(z_n)$ is contained in 
$Y:=\{w \in \psi_{0}(U_0) 
\st |\arg (w^q)| \le \pi/3,\, |w| \ge C_{7} \nu_n \cdot (3/4)\}$
for sufficiently large $n$.
It follows that 
$|\psi_{0}(\zeta_n) -\psi_{0}(z_n)|$
is larger than the distance between the sets $X$ and $Y$,
which is comparable with the radius $\nu_n$ of $V_0$
 (see Figure \ref{fig_XY}).
However, the univalence of $\psi_0$ implies that $|\psi_{0}(\zeta_n) -\psi_{0}(z_n)|$ is not comparable with $\nu_n$ by
$|\psi_{0}(\zeta_n) -\psi_{0}(z_n)|
\asymp |\zeta_n - z_n|< \nu_n/n$. 
This is a contradiction. 
\QED

\section{Proof of Lemma C} \label{sec:Proof_Of_Lemma_C}
We will show that $|Df_c^{M'-M}(z_M)| \ge \kappa_{\mathrm{C}}/\nu$
for some constant $\kappa_{\mathrm{C}}$ 
that does not depend on $c=c_t$ with $t \in \Delta$.
By choosing $\nu$ sufficiently small, 
we have $\Lam:=\kappa_{\mathrm{C}}/\nu>1$. 

As in the proof of Lemma A, we assume that $M=0$ 
and $M'=mpq +L$ for which 
$z_0 \in V_0$ (equivalently, $|z_{0} -\bhat| < \nu$), 
$z_{(m-1)pq} \in U_0$,
$z_{mpq} \notin U_0$ (equivalently, 
$|z_{mpq} -\bhat| \ge R_0$),
and $z_{M'} \in V_0$. 
By the chain rule we have
\begin{equation}
|Df_c^{M'}(z_0)|
=|Df_c^{mpq}(z_0)|\cdot |Df_c^{L}(z_{mpq})|.
\label{eq_C_6-0}
\end{equation}

First let us give an estimate of $|Df_c^{mpq}(z_0)|$.
For $c =c_t~(t \in \Delta)$, 
let $b(c)$ be the periodic point $b_t$ of period $p$
given in Proposition \ref{prop_near_parabolic}
with $b(c) \to \bhat$ as $c =c_t$ tends to $\chat$.
We may assume that $|b(c)-\bhat| \le \nu$
for $c \approx \chat$, 
and by the Koebe distortion Theorem, we have 
$$
|Df_c^{mpq}(z_0)| \asymp 
\frac{|z_{mpq}-b(c)|}{|z_0-b(c)|}
\ge \frac{|z_{mpq}-\bhat|-|b(c)-\bhat|}{|z_0-\bhat|+|b(c)-\bhat|}
\ge \frac{R_0-\nu}{\nu + \nu}.
$$
This implies that $|Df_c^{mpq}(z_0)| \ge C_8/\nu$
for some constant $C_8>0$ independent of $c\approx \chat$, $z_0$ and $\nu \ll 1$. 

Next we give an estimate of the form 
$|Df_c^{L}(z_{mpq})| \ge C_9$, 
 where $C_9$ is a constant 
independent of $c \approx \chat$, $\nu \ll 1$, 
and $z_0 \in J(f_c)$.
(Then by (\ref{eq_C_6-0}) the proof is done by setting 
$\kappa_{\mathrm C} :=C_8 C_9$.)

As in the proof of Lemma A, 
by shrinking $U_0$ if necessary, 
we may always assume that $L >p$. 
Since $z_{M'} \in V_0$, 
we have $z_{M'-p} \in f_c^{-p}(V_0)$
and thus $|Df_c^p(z_{M'-p})| \ge \xi^p$.
By Proposition \ref{previous_Lemma_W} and inequality \eqref{ineqInF} in Lemma F we obtain 
\begin{align*}
|Df_c^{L}(z_{mpq})|
=&|Df_c^{L-p}(z_{mpq})||Df_c^p(z_{M'-p})|\\
\ge& 
\frac{\rho(z_{mpq})}{\rho(z_{M'-p})} A^{L-p} \cdot \xi^p \\
\ge& 
\rho(z_{mpq})\cdot 
\dist(z_{M'-p}, P(f_\chat)) \cdot
\xi^p.
\end{align*}
Since $z_{mpq} \notin U_0$ 
and $z_{mpq}$ has a definite distance from 
the parabolic cycle, 
the distance between $z_{mpq}$ and the postcritical set 
$P(f_\chat)$ is larger  than 
a positive constant independent of $c \approx \chat$,
$\nu \ll 1$, and $z_0 \in J(f_c)$.
(Indeed, we may apply the proof of Lemma D 
to obtain $\dist(z_{mpq}, P(f_\chat)) \ge K_{\mathrm D} R_0$
by taking a smaller $R_0$ if necessary.) 
Hence we always have $\rho(z_{mpq}) \asymp 1$.
Finally we show that 
$\dist(z_{M'-p}, P(f_\chat))$
has a uniform distance away from zero.
Since $z_{M'} \in V_0$, 
$\dist(z_{M'-p}, P(f_\chat))$ can be close to $0$ 
only if $z_{M'-p}$ belongs to 
the connected component $V_0'$ of $f_c^{-p}(V_0)$ contained in $U_0$.
(See Figure \ref{fig_S-cycle}.)
However, the local dynamics of $f_c^p$ on 
$U_0$ does not allow any point $z \in J(f_c)$ satisfying both
$z \in V_0'-V_0$ and $f_c^p(z) \in V_0$.
(More precisely, we may apply 
Proposition \ref{prop_G_near_origin}
and Lemma I 
 to exclude such a $z$ in the Julia set.)
Hence $z_{M'-p}$ belongs to a connected component of 
$f_c^{-p}(V_0)$ that has a definite distance from $P(f_\chat)$ 
and $\dist(z_{M'-p}, P(f_\chat))$ is bounded below by
a positive constant independent of $c \approx \chat$,
$\nu \ll 1$, and $z_0 \in J(f_c)$.
In conclusion, $|Df_c^{L}(z_{mpq})|$ is bounded away from zero by a constant
independent of $c \approx \chat$, $\nu \ll 1$, 
and $z_0 \in J(f_c)$. \QED

\section{Proofs of Theorems \ref{thm_HolderContinuity}~and \ref{thm_H2P}} \label{sec:together}
In this section we prove Theorem \ref{thm_HolderContinuity} and Theorem \ref{thm_H2P} together. 

\paragraph{Proofs of Theorems \ref{thm_HolderContinuity} and 
\ref{thm_H2P}.}
Let $\sigma_0$ be any point in the thick internal ray 
$\cI(\theta,\delta)$ and $\ell:=|\sigma_0-\chat|$. 
For each $n=0,1,\ldots$, let
$$
E_n:= \brac{c \in \cI(\theta,\delta) \st 
\frac{\ell}{2^{n+1}}\le |c-\chat| \le \frac{\ell}{2^n}
}.
$$
Note that $\mathrm{diam}\,E_n \asymp \ell/2^n$ and that 
for any $c \in E_n$ we have $|c-\chat| \asymp \ell/2^n$. 

We take a $\sigma_n \in E_n$ for each $n \ge 1$. 
Then we can join $\sigma_n$ and 
$\sigma_{n+1}$ by a piecewise smooth path $\gamma_n$
of length compatible with $\ell/2^n$ contained in $E_n \cup E_{n+1}$.
By the Main Theorem we have
$$
\abs{\frac{d}{dc}z(c)} \le \frac{K}{|c-\chat|^{1-1/Q}}
\asymp \paren{\frac{\ell}{2^n}}^{-1+1/Q} 
$$
for any $c \in \gamma_n$ and $n \ge 0$. Hence
$$
|z(\sigma_{n+1}) - z(\sigma_n)| = \abs{\int_{\gamma_n}\frac{d}{dc}z(c) \,dc}
= O\paren{\paren{\frac{\ell}{2^n}}^{-1+1/Q} \cdot \frac{\ell}{2^n} }
= O\paren{\frac{1}{2^{n/Q}}}\cdot \ell^{1/Q}.
$$
It follows that the sequence $\{z(\sigma_n)\}_{n \ge 0}$ is Cauchy,
and we denote the limit by $z(\chat)$.
(One can easily check that the limit does not depend on the 
choice of the sequences $\{\sigma_n\}$ and $\{\gamma_n\}$.)
Moreover, we have
$$
|z(\sigma_0) -z(\sigma_n)| 
\le O\paren{1+\frac{1}{2^{1/Q}}+\cdots
+ \frac{1}{2^{(n-1)/Q}}} 
\cdot \ell^{1/Q}
=O(\ell^{1/Q})
$$
and thus $|z(\sigma_0)-z(\chat)| =O(\ell^{1/Q})=O(|\sigma_0-\chat|^{1/Q})$,
where the implicit constants depend only on $\chat$ and 
the thickness $\delta$ of the thick internal ray.
This proves  \eqref{eq_thm1.2} of Theorem \ref{thm_HolderContinuity}.

\paragraph{Semiconjugacy.}
Next we show that $z(\chat)$ belongs to the Julia set $J(f_\chat)$.
For each $z_\ast \in J(f_{\sigma})$ and 
its motion $z(c)=h_c(z_\ast)=H(c,z_\ast)$ along 
the thick internal ray $\cI(\theta, \delta)$, 
we define $h_\chat(z_\ast)$ by the limit $z(\chat)$ 
given as above. 
Since $h_c:J(f_{\sigma}) \to J(f_c)$ 
is continuous and the convergence of 
$h_c$ to $h_\chat$ as $c \in \cI(\theta, \delta)$ tends to $\chat$
is uniform, $h_\chat:J(f_{\sigma}) \to h_\chat(J(f_{\sigma}))$ 
is continuous as well.
Hence by $f_c \circ h_c = h_c \circ f_{\sigma}$ 
we obtain $f_\chat \circ h_\chat = h_\chat \circ f_{\sigma}$.
In particular, this implies that 
$f_\chat \circ h_\chat(J(f_{\sigma}))=h_\chat(J(f_{\sigma}))$
and thus the image $h_\chat(J(f_{\sigma}))$ 
is a compact forward invariant set contained in the filled Julia set 
of $f_\chat$.
Suppose that there exists a $z_\ast \in J(f_\sigma)$
such that $h_\chat(z_\ast)=z(\chat) \in h_\chat(J(f_{\sigma}))$ belongs to the Fatou set 
of $f_\chat$. 
It actually belongs to the basin of attraction of $\hat{b}$, 
and thus there exists an $l$ such that $f_{\chat}^{l +kpq}(z(\chat))$
is contained in $V_0$ for any $k \ge 0$. 
Since $h_\chat$ is continuous, the image of any nearby point of $z_*$ under $h_\chat$ also belongs to the same basin. In addition, since the points that are not eventually periodic in the dynamics of $f_\sigma$ are dense in $J(f_{\sigma})$, 
we may assume that the orbit $f_\sigma^n(z_\ast)~(n \in \N_0)$ 
of $z_\ast$ never lands on the periodic points.
By $f_{\chat}^{l +kpq}(z(\chat)) =h_\chat(f_\sigma^{l +kpq}(z_\ast))$
and uniform convergence of $h_c$ to $h_\chat$,
we have $h_c(f_\sigma^{l +kpq}(z_\ast)) \in V_0$
for any $k$ and $c \approx \chat$.
Since $h_c(f_\sigma^{l +kpq}(z_\ast)) \in J(f_c)$, 
Proposition \ref{prop_leaving_U_0} implies that 
$h_c(f_\sigma^{l +kpq}(z_\ast))$ is a repelling periodic point
contained in $U_0$.
This is impossible because $z_\ast$ is not eventually periodic.
This completes the proof of Theorem \ref{thm_HolderContinuity}.

To confirm that $h_\chat$ is a semiconjugacy, 
we show surjectivity of $h_{\chat}:J(f_{\sigma}) \to J(f_\chat)$:
First we take any repelling periodic point $\hat{x} \in J(f_\chat)$.
Since there is a holomorphic family $x(c)$ of repelling periodic points 
for $c$ sufficiently close to $\chat$ such that $\hat{x}=x(\chat)$, 
we have some $z_0 \in J(f_{\sigma})$ with 
$h_c(z_0)=x(c)$ for any $c \in \cI(\theta, \delta)$.
In particular, we have $h_\chat(z_0)=\hat{x}$. 
Next we take any $w \in J(f_{\chat})$ and a sequence of repelling 
periodic points $\hat{x}_n$ of $f_{\chat}$ that converges to $w$ as $n \to \infty$.
(Such a sequence exists since repelling periodic points are dense in the Julia set.)
Let $z_n \in J(f_{\sigma})$ be the repelling periodic point with
$h_\chat(z_n)=\hat{x}_n$. 
Then any accumulation point $y$ of the sequence $z_n$ 
satisfies $h_\chat(y)=w$ by continuity.

Finally we check properties (1) -- (3) of Theorem \ref{thm_H2P}.
 Property (3) is an immediate consequence of 
 \eqref{eq_thm1.2} in Theorem \ref{thm_HolderContinuity}. 
(In particular, we obtain Corollary \ref{cor_HausdorffConvergence}.)
To show (1) and (2), 
let $\eta_c:=h_\chat \circ h_c^{-1}:J(f_c) \to J(f_\chat)$ 
for $c \approx \chat$ such that $\eta_c$ 
is a semiconjugacy between $f_c$ and $f_\chat$
and satisfies $|\eta_c(z)-z| \le K'|c-\chat|^{1/Q}$ by (3).
We assume that $c$ is sufficiently close to $\chat$ 
such that $K'|c-\chat|^{1/Q} < K_{\mathrm D}\,\nu/4$,
where the constant $K_{\mathrm D} \le 1$ is given in Lemma D. 
Suppose that $\eta_c(z)=\eta_c(z')$ for some distinct points 
$z, \, z' \in J(f_c)$.
Let $z_n:=f_c^n(z)$ and $z_n':=f_c^n(z')$ for $n \in \N$.
Then $\eta_c \circ f_c = f_\chat \circ \eta_c$ implies 
$\eta_c(z_n)=\eta_c(z_n')$ for any $n \in \N$,
and thus 
\begin{equation}\label{eq_pf_thm_H2P}
|z_n-z_n'| \le |\eta_c(z_n)-z_n|+|\eta_c(z_n')-z_n'|
=2K'|c-\chat|^{1/Q} < K_{\mathrm D} \nu/2 < \nu/2.
\end{equation}
Suppose that the orbit $z_n~(n \in \N)$ 
never lands on repelling periodic points of 
$f_c$ in $\hat{U} =\D(\bhat, \hat{R})$ 
described in assertion (3) or (4) of 
Proposition \ref{prop_near_parabolic}. 
By Proposition \ref{prop_leaving_U_0}, 
the orbit $z_n'~(n \in \N)$ must behave in the same way. 
Let $\{n_k\}_{k \in \N}$ be a subsequence 
such that $z_n, z_n' \in J(f_c)-\cV(c)$ when $n=n_k$ for each $k$.
By Lemma D, $\D(z_n, K_{\mathrm D} \nu)$ does not contain any point in 
the postcritical set $P(f_\chat)$
when $n$ ranges over the subsequence $\{n_k\}_{k \in \N}$.
Hence we have a univalent branch $g_n$ of $f_c^{-n}$ 
defined on $\D(z_n, K_{\mathrm D} \nu)$
that sends $z_n$ and $z_n'$ to $z$ and $z'$.
However, since $f_c$ is hyperbolic,
the Koebe distortion theorem (see  
\cite[Theorem 5-3]{Ah}) implies
$$
|z-z'| \le \diam g_n(\D(z_n, K_{\mathrm D} \nu/2))
\asymp \frac{\nu}{|Df_c^n(z)|} \to 0
$$
as $n=n_k \to \infty$. 
This contradicts the assumption $z \neq z'$.
Hence we may assume that the orbit $z_n~(n \in \N)$ 
lands on a repelling fixed point of 
$f_c^{pq}$ in $V_0 =\D(\bhat, \nu)$.
Let $z_N$ be such a repelling fixed point.

\paragraph{Case 1 or Case 2$^+$.} 
By \eqref{eq_pf_thm_H2P} and Proposition \ref{prop_leaving_U_0},
$z_N'$ must be the only repelling fixed point of $f_c^{pq}$ in $V_0$ 
given in assertion (3) of Proposition \ref{prop_near_parabolic},
and thus $z_N=z_N'$. 
Since $z \neq z'$, we must have $z_n \neq z_n'=-z_n$ 
for some $n <N$.
However, this is impossible 
because the Julia set $J(f_\chat)$ 
and the critical point 0 have a definite distance,
and the same holds for $J(f_c)$ by \eqref{eq_H2P} for $c \approx \chat$.
(Indeed, $J(f_c)$ converges to $J(f_\chat)$ as $c$ tends to $\chat$
along a thick internal ray. See Corollary \ref{cor_HausdorffConvergence}.)
Hence we conclude that $\eta_c$ is injective. 
Since $h_\chat=\eta_c \circ h_c$ and $h_c$ is a conjugacy,
property (1) holds. 

\paragraph{Case 2$^-$.}
By the same argument as above, 
$z_N'$ must be one of 
the $q$ repelling fixed points of 
$f_c^{pq}$ in $V_0$ given in assertion (4) of 
Proposition \ref{prop_near_parabolic}. 
Hence we conclude that 
$\eta_c(z_N)$ is a parabolic periodic point $\bhat$
and property (2) holds by the relation $h_\chat=\eta_c \circ h_c$.
\QED

\section*{Acknowledgments}
Chen was partly supported by MOST 108-2115-M-001-005 and 109-2115-M-001-006. 
Kawahira was partly supported by JSPS KAKENHI Grants  numbers
16K05193 and 19K03535. 
They thank the hospitality of Academia Sinica, Nagoya University, 
RIMS in Kyoto University, and Tokyo Institute of Technology
where parts of this research were carried out.

\vspace{1cm}

{~}\\
Yi-Chiuan Chen \\
Institute of Mathematics\\
Academia Sinica\\
Taipei 106319, Taiwan\\
YCChen@math.sinica.edu.tw

\vspace{.5cm}

{~}\\
Tomoki Kawahira\\
Graduate School of Economics\\ 
Hitotsubashi University\\
Tokyo 186-8602, Japan\\
t.kawahira@r.hit-u.ac.jp

\end{document}